\documentclass[a4paper,11pt]{article}

\usepackage[T1]{fontenc}
\usepackage{lmodern}

\usepackage[latin1]{inputenc}
\usepackage{amsmath}
\usepackage{amsthm}
\usepackage{amssymb}
\usepackage{mathrsfs}
\usepackage{graphicx}
\usepackage[all]{xy}
\usepackage{hyperref}

\usepackage{abstract} 

\newtheorem{thm}{Theorem}[section]
\newtheorem{cor}[thm]{Corollary}
\newtheorem{claim}[thm]{Claim}
\newtheorem{fact}[thm]{Fact}
\newtheorem{lemma}[thm]{Lemma}
\newtheorem{prop}[thm]{Proposition}

\theoremstyle{definition}
\newtheorem{definition}[thm]{Definition}
\newtheorem{ex}[thm]{Example}
\newtheorem{remark}[thm]{Remark}
\newtheorem{question}[thm]{Question}

\title{Contracting isometries of CAT(0) cube complexes and acylindrical hyperbolicity of diagram groups}
\date{\today}
\author{Anthony Genevois}

\begin{document}

\maketitle

\begin{abstract}
The main technical result of this paper is to characterize the contracting isometries of a CAT(0) cube complex without any assumption on its local finiteness. Afterwards, we introduce the combinatorial boundary of a CAT(0) cube complex, and we show that contracting isometries are strongly related to isolated points at infinity, when the complex is locally finite. This boundary turns out to appear naturally in the context of Guba and Sapir's diagram groups, and we apply our main criterion to determine precisely when an element of a diagram group induces a contracting isometry on the associated Farley cube complex. As a consequence, in some specific case, we are able to deduce a criterion to determine precisely when a diagram group is acylindrically hyperbolic.
\end{abstract}

\tableofcontents

\section{Introduction}

Given a metric space $X$, an isometry $g \in X$ is \emph{contracting} if
\begin{itemize}
	\item $g$ is \emph{loxodromic}, ie., there exists $x_0 \in X$ such that $n \mapsto g^n \cdot x_0$ defines a quasi-isometric embedding $\mathbb{Z} \to X$;
	\item if $C_g= \{ g^n \cdot x_0 \mid n \in \mathbb{Z} \}$, then the diameter of the nearest-point projection of any ball disjoint from $C_g$ onto $C_g$ is uniformly bounded.
\end{itemize}
For instance, any loxodromic isometry of a Gromov-hyperbolic space is contracting. In fact, if a group $G$ acts by isometries on a metric space, the existence of a contracting isometry seems to confer to $G$ a ``hyperbolic behaviour''. To make this statement more precise, one possibility is to introduce \emph{acylindrically hyperbolic groups} as defined in \cite{arXiv:1304.1246}.

\begin{definition}
Let $G$ be a group acting on a metric space $X$. The action $G \curvearrowright X$ is \emph{acylindrical} if, for every $d \geq 0$, there exist some constants $R,N \geq 0$ such that, for every $x,y \in X$,
\begin{center}
$d(x,y) \geq R \Rightarrow \# \{ g \in G \mid d(x,gx),d(y,gy) \leq d \} \leq N$.
\end{center}
A group is \emph{acylindrically hyperbolic} if it admits a \emph{non-elementary} (ie., with an infinite limit set) acylindrical action on a Gromov-hyperbolic space.
\end{definition}

\noindent
Acylindrically hyperbolic groups may be thought of as a generalisation of relatively hyperbolic groups. See for example \cite[Appendix]{arXiv:1304.1246} and references therein for more details. The link between contracting isometries and acylindrically hyperbolic groups is made explicit by the following result, which is a consequence of \cite[Theorem H]{arXiv:1006.1939} and \cite[Theorem 1.2]{arXiv:1304.1246}:

\begin{thm}\label{BFO}
If a group acts by isometries on a geodesic space with a WPD contracting isometry, then it is either virtually cyclic or acylindrically hyperbolic.
\end{thm}

\noindent
We do not give the precise definition of a WPD isometry. The only thing to know in our paper is that any isometry turns out to be WPD if the action of our group is properly discontinuous. 

\medskip \noindent
In this paper, we focus on a specific class of geodesic spaces: the CAT(0) cube complexes, ie., simply connected cellular complexes obtained by gluing cubes of different dimensions by isometries along their faces so that the link of every vertex is a simplicial flag complex. See Section 2 for more details. The first question we are interested in is: when is an isometry of a CAT(0) cube complex contracting? The answer we give to this question is the main technical result of our paper. Our criterion essentially deals with hyperplane configurations in the set $\mathcal{H}(\gamma)$ of the hyperlanes intersecting a combinatorial axis $\gamma$ of a given loxodromic isometry $g$. We refer to Section 3 for precise definitions.

\begin{thm}\label{main1}
Let $X$ be a CAT(0) cube complex and $g \in \mathrm{Isom}(X)$ an isometry with a combinatorial axis $\gamma$. The following statements are equivalent:
\begin{itemize}
	\item[(i)] $g$ is a contracting isometry;
	\item[(ii)] joins of hyperplanes $(\mathcal{H}, \mathcal{V})$ with $\mathcal{H} \subset \mathcal{H}(\gamma)$ are uniformly thin;
	\item[(iii)] there exists $C \geq 1$ such that:
		\begin{itemize}
			\item $\mathcal{H}(\gamma)$ does not contain $C$ pairwise transverse hyperplanes,
			\item any grid of hyperplanes $(\mathcal{H}, \mathcal{V})$ with $\mathcal{H} \subset \mathcal{H}(\gamma)$ is $C$-thin;
		\end{itemize}
	\item[(iv)] $g$ skewers a pair of well-separated hyperplanes.
\end{itemize}
\end{thm}

\begin{remark}
The equivalence $(i) \Leftrightarrow (iv)$ generalises a similar criterion established in \cite[Theorem 4.2]{MR3339446}, where the cube complex is supposed uniformly locally finite.
\end{remark}

\noindent
In view of our application to diagram groups, it turned out to be natural to introduce a new boundary of a CAT(0) cube complex; in the definition below, if $r$ is a combinatorial ray, $\mathcal{H}(r)$ denotes the set of the hyperplanes which intersect $r$.

\begin{definition}
Let $X$ be a CAT(0) cube complex. Its \emph{combinatorial boundary} is the poset $(\partial^c X,\prec)$, where $\partial^c X$ is the set of the combinatorial rays modulo the relation: $r_1 \sim r_2$ if $\mathcal{H}(r_1) \underset{a}{=} \mathcal{H}(r_2)$; and where the partial order $\prec$ is defined by: $r_1 \prec r_2$ whenever $\mathcal{H}(r_1) \underset{a}{\subset} \mathcal{H}(r_2)$.
\end{definition}

\noindent
Here, we used the following notation: given two sets $X,Y$, we say that $X$ and $Y$ are \emph{almost-equal}, denoted by $X \underset{a}{=} Y$, if the symmetric difference between $X$ and $Y$ is finite; we say that $X$ is \emph{almost-included} into $Y$, denoted by $X \underset{a}{\subset} Y$, if $X$ is almost-equal to a subset of $Y$. In fact, this boundary is not completely new, since it admits strong relations with other boundaries; see Appendix A for more details. 

\medskip \noindent
A point in the combinatorial boundary is \emph{isolated} if it is not comparable with any other point. Now, thanks to Theorem \ref{main1}, it is possible to read at infinity when an isometry is contracting:

\begin{thm}\label{main2}
Let $X$ be a locally finite CAT(0) cube complex and $g\in \mathrm{Isom}(X)$ an isometry with a combinatorial axis $\gamma$. Then $g$ is a contracting isometry if and only if $\gamma(+ \infty)$ is isolated in $\partial^cX$. 
\end{thm}

\noindent
Thus, the existence of a contracting isometry implies the existence of an isolated point in the combinatorial boundary. Of course, the converse cannot hold without an additional hypothesis on the action, but in some specific cases we are able to prove partial converses. For instance:

\begin{thm}\label{main3}
Let $G$ be a group acting on a locally finite CAT(0) cube complex $X$ with finitely many orbits of hyperplanes. Then $G \curvearrowright X$ contains a contracting isometry if and only if $\partial^c X$ has an isolated vertex.
\end{thm}

\begin{thm}\label{main4}
Let $G$ be a countable group acting on a locally finite CAT(0) cube complex $X$. Suppose that the action $G \curvearrowright X$ is minimal (ie., $X$ does not contain a proper $G$-invariant combinatorially convex subcomplex) and $G$ does not fix a point of $\partial^cX$. Then $G$ contains a contracting isometry if and only if $\partial^c X$ contains an isolated point.
\end{thm}

\noindent
As mentionned above, the main reason we introduce combinatorial boundaries of CAT(0) cube complexes is to apply our criteria to Guba and Sapir's diagram groups. Loosely speaking, diagram groups are ``two-dimensional free groups'': in the same way that free groups are defined by concatenating and reducing words, diagram groups are defined by concatenating and reducing some two-dimensional objets, called \emph{semigroup diagrams}. See Section \ref{section:preliminarydiag} for precise definitions. Although these two classes of groups turn out to be quite different, the previous analogy can be pushed further. On the one hand, free groups act on their canonical Cayley graphs, which are simplicial trees; on the other hand, diagram groups act on the natural Cayley graphs of their associated groupoids, which are CAT(0) cube complexes, called \emph{Farley complexes}. Moreover, in the same way that the boundary of a free group may be thought of as a set of infinite reduced words, the combinatorial boundary of a Farley cube complex may be thought of as a set of infinite reduced diagrams. See Section \ref{section:Fboundary} for a precise description. 

\medskip \noindent
If $g$ is an element of a diagram group, which is \emph{absolutely reduced} (ie., the product $g^n$ is reduced for every $n \geq 1$), let $g^{\infty}$ denote the infinite diagram obtained by concatenating infinitely many copies of $g$. Then, thanks to Theorem \ref{main2} and a precise description of the combinatorial boundaries of Farley cube complexes, we are able to deduce the following criterion:

\begin{thm}\label{main5}
Let $G$ be a diagram group and $g \in G \backslash \{ 1 \}$ be absolutely reduced. Then $g$ is a contracting isometry of the associated Farley complex if and only if the following two conditions are satisfied:
\begin{itemize}
	\item $g^{\infty}$ does not contain any infinite proper prefix;
	\item for any infinite reduced diagram $\Delta$ containing $g^{\infty}$ as a prefix, all but finitely many cells of $\Delta$ belong to $g^{\infty}$.
\end{itemize}
\end{thm}

\noindent
Of course, although it is sufficient for elements with few cells, this criterion may be difficult to apply in practice, because we cannot draw the whole $g^{\infty}$. This why we give a more ``algorithmic'' criterion in Section 5.3.

\medskip \noindent
Thus, we are able to determine precisely when a given element of a diagram group is a contracting isometry. In general, if there exist such isometries, it is not too difficult to find one of them. Otherwise, it is possible to apply Theorem \ref{main3} or Theorem \ref{main4}; we emphasize that Theorem \ref{main3} may be particularly useful since diagram groups often act on their cube complexes with finitely many orbits of hyperplanes. Conversely, we are able to state that no contracting isometry exists if the combinatorial boundary does not contain isolated points, see for instance Example \ref{Thompson}.

\medskip \noindent
Therefore, we get powerful tools to pursue the cubical study of negatively-curved properties of diagram groups we initialized in \cite{arXiv:1505.02053}. To be precise, the question we are interested in is:

\begin{question}\label{Q1}
When is a diagram group acylindrically hyperbolic?
\end{question}

\noindent
Thanks to Theorem \ref{BFO}, we are able to deduce that a diagram group is acylindrically hyperbolic if it contains a contracting isometry and if it is not cyclic. Section \ref{section:examples} provides families of acylindrically hyperbolic diagram groups which we think to be of interest. On the other hand, a given group may be represented as a diagram group in different ways, and an acylindrically hyperbolic group (eg., a non-abelian free group) may be represented as a diagram group so that there is no contracting isometry on the associated cube complex. We do not know if it is always possible to a find a ``good'' representation. 

\medskip \noindent
Nevertheless, a good class of diagram groups, where the problem mentionned above does not occur, corresponds to the case where the action on the associated complex is cocompact. Using the notation introduced below, they correspond exactly to the diagram groups $D(\mathcal{P},w)$ where the class of $w$ modulo the semigroup $\mathcal{P}$ is finite. We call them \emph{cocompact diagram groups}. Focusing on this class of groups, we prove (see Theorem \ref{decompositioncor} for a precise statement):

\begin{thm}\label{main7}
A cocompact diagram group decomposes naturally as a direct product of a finitely generated free abelian group and finitely many acylindrically hyperbolic diagram groups.
\end{thm}

\noindent
On the other hand, acylindrically hyperbolic groups cannot split as a direct product of two infinite groups (see \cite[Corollary 7.3]{arXiv:1304.1246}), so we deduce a complete answer of Question \ref{Q1} in the cocompact case:

\begin{cor}
A cocompact diagram group is acylindrically hyperbolic if and only if it is not cyclic and it does not split as a non-trivial direct product.
\end{cor}

\noindent
Notice that this statement is no longer true without the cocompact assumption. Indeed, Thompson's group $F$ and the lamplighter group $\mathbb{Z} \wr \mathbb{Z}$ are diagram groups which are not acylindrically hyperbolic (since they do not contain a non-abelian free group) and they do not split non-trivially as a direct product. Another consequence of Theorem \ref{main7} is that a cocompact diagram group either is free abelian or has a quotient which is acylindrically hyperbolic. Because acylindrically hyperbolic groups are \emph{SQ-universal} (ie., any countable group can be embedded into a quotient of the given group; see \cite[Theorem 8.1]{arXiv:1304.1246}), we get the following dichotomy:

\begin{cor}
A cocompact diagram group is either free abelian or SQ-universal.
\end{cor}

\noindent
In our opinion, cocompact diagram groups turn out to be quite similar to (finitely generated) right-angled Artin groups. In this context, Theorem \ref{main7} should be compared to the similar statement, but more general, \cite[Theorem 5.2]{BehrstockCharney} (it must be noticed that, although the statement is correct, the proof contains a mistake; see \cite[Remark 6.21]{arXiv:1310.6289}). Compare also \cite[Lemma 5.1]{BehrstockCharney} with Proposition \ref{decompositionprop}. Notice that there exist right-angled Artin groups which are not (cocompact) diagram groups (see \cite[Theorem 30]{MR1725439}), and conversely there exist cocompact diagram groups which are not right-angled Artin groups (see Example \ref{nonRAAG}).

\medskip \noindent
The paper is organized as follows. In Section 2, we give the prerequisites on CAT(0) cube complexes needed in the rest of the paper. Section 3 is essentially dedicated to the proof of Theorem \ref{main1}, and in Section 4, we introduce combinatorial boundaries of CAT(0) cube complexes and we prove Theorem \ref{main2}, as well as Theorem \ref{main3} and Theorem \ref{main4}. Finally, in Section 5, we introduce diagram groups and we apply the results of the previous sections to deduce the various statements mentionned above. We added an appendix at the end of this paper to compare the combinatorial boundary with other known boundaries of CAT(0) cube complexes.

\paragraph{Acknowledgement.} I am grateful to Jean L\'ecureux for having suggested me a link between the combinatorial boundary and the Roller boundary, and to my advisor, Peter Ha\"{i}ssinsky, for all our discussions.

\section{Preliminaries}

A \textit{cube complex} is a CW complex constructed by gluing together cubes of arbitrary (finite) dimension by isometries along their faces. Furthermore, it is \textit{nonpositively curved} if the link of any of its vertices is a simplicial \textit{flag} complex (ie., $n+1$ vertices span a $n$-simplex if and only if they are pairwise adjacent), and \textit{CAT(0)} if it is nonpositively curved and simply-connected. See \cite[page 111]{MR1744486} for more information.

Alternatively, CAT(0) cube complexes may be described by their 1-skeletons. Indeed, Chepoi notices in \cite{mediangraphs} that the class of graphs appearing as 1-skeletons of CAT(0) cube complexes coincides with the class of \textit{median graphs}, which we now define.

Let $\Gamma$ be a graph. If $x,y,z \in \Gamma$ are three vertices, a vertex $m$ is called a \textit{median point of $x,y,z$} whenever
\begin{center}
$d(x,y)=d(x,m)+d(m,y)$, $d(x,z)=d(x,m)+d(m,z)$, $d(y,z)=d(y,m)+d(m,z)$.
\end{center}
Notice that, for every geodesics $[x,m]$, $[y,m]$ and $[z,m]$, the concatenations $[x,m] \cup [m,y]$, $[x,m] \cup [m,z]$ and $[y,m] \cup [m,z]$ are also geodesics; furthermore, if $[x,y]$, $[y,z]$ and $[x,z]$ are geodesics, then any vertex of $[x,y] \cap [y,z] \cap [x,z]$ is a median point of $x,y,z$.

The graph $\Gamma$ is \textit{median} if every triple $(x,y,z)$ of pairwise distinct vertices admits a unique median point, denoted by $m(x,y,z)$.

\begin{thm}\emph{\cite[Theorem 6.1]{mediangraphs}} 
A graph is median if and only if it is the 1-skeleton of a CAT(0) cube complex.
\end{thm}

A fundamental feature of cube complexes is the notion of \textit{hyperplane}. Let $X$ be a nonpositively curved cube complex. Formally, a \textit{hyperplane} $J$ is an equivalence class of edges, where two edges $e$ and $f$ are equivalent whenever there exists a sequence of edges $e=e_0,e_1,\ldots, e_{n-1},e_n=f$ where $e_i$ and $e_{i+1}$ are parallel sides of some square in $X$. Notice that a hyperplane is uniquely determined by one of its edges, so if $e \in J$ we say that $J$ is the \textit{hyperplane dual to $e$}. Geometrically, a hyperplane $J$ is rather thought of as the union of the \textit{midcubes} transverse to the edges belonging to $J$.
\begin{center}
\includegraphics[scale=0.4]{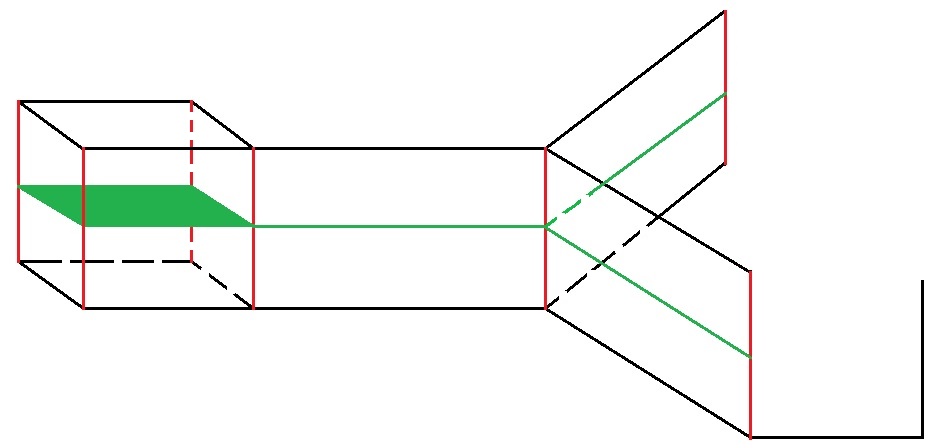}
\end{center}
The \textit{neighborhood} $N(J)$ of a hyperplane $J$ is the smallest subcomplex of $X$ containing $J$, i.e., the union of the cubes intersecting $J$. In the following, $\partial N(J)$ will denote the union of the cubes of $X$ contained in $N(J)$ but not intersecting $J$, and $X \backslash \backslash J= \left( X \backslash N(J) \right) \cup \partial N(J)$. Notice that $N(J)$ and $X \backslash \backslash J$ are subcomplexes of $X$.

\begin{thm}\emph{\cite[Theorem 4.10]{MR1347406}}
Let $X$ be a CAT(0) cube complex and $J$ a hyperplane. Then $X \backslash \backslash J$ has exactly two connected components.
\end{thm}

\noindent
The two connected components of $X \backslash \backslash J$ will be refered to as the \textit{halfspaces} associated to the hyperplane $J$.

\paragraph{Distances $\ell_p$.} 
There exist several natural metrics on a CAT(0) cube complex. For example, for any $p \in (0,+ \infty)$, the $\ell_p$-norm defined on each cube can be extended to a distance defined on the whole complex, the \emph{$\ell_p$-metric}. Usually, the $\ell_1$-metric is referred to as the \emph{combinatorial distance} and the $\ell_2$-metric as the \emph{CAT(0) distance}. Indeed, a CAT(0) cube complex endowed with its CAT(0) distance turns out to be a CAT(0) space \cite[Theorem C.9]{Leary}, and the combinatorial distance between two vertices corresponds to the graph metric associated to the 1-skeleton $X^{(1)}$. In particular, \textit{combinatorial geodesics} are edge-paths of minimal length, and a subcomplex is \textit{combinatorially convex} if it contains any combinatorial geodesic between two of its points.

In fact, the combinatorial metric and the hyperplanes are strongly linked together: the combinatorial distance between two vertices corresponds exactly to the number of hyperplanes separating them \cite[Theorem 2.7]{MR2413337}, and

\begin{thm}\emph{\cite[Corollary 2.16]{MR2413337}}
Let $X$ be a CAT(0) cube complex and $J$ a hyperplane. The two components of $X \backslash \backslash J$ are combinatorially convex, as well as the components of $\partial N(J)$.
\end{thm}

\noindent
This result is particularly useful when it is combined with the following well-known Helly property; see for example \cite[Theorem 2.2]{Roller}.

\begin{thm}
If $\mathcal{C}$ is a finite collection of pairwise intersecting combinatorially subcomplexes, then the intersection $\bigcap\limits_{C \in \mathcal{C}} C$ is non-empty.
\end{thm}

\paragraph{Combinatorial projection.}
In CAT(0) spaces, and so in particular in CAT(0) cube complexes with respect to the CAT(0) distance, the existence of a well-defined projection onto a given convex subspace provides a useful tool. Similarly, with respect to the combinatorial distance, it is possible to introduce a \emph{combinatorial projection} onto a combinatorially convex subcomplex, defined by the following result.

\begin{prop}\label{projection}
\emph{\cite[Lemma 1.2.3]{arXiv:1505.02053}} Let $X$ be a CAT(0) cube complex, $C \subset X$ be a combinatorially convex subspace and $x \in X \backslash C$ be a vertex. Then there exists a unique vertex $y \in C$ minimizing the distance to $x$. Moreover, for any vertex of $C$, there exists a combinatorial geodesic from it to $x$ passing through $y$.
\end{prop}

\noindent
The following result makes precise how the combinatorial projection behaves with respect to the hyperplanes.

\begin{prop}\label{hyperplanseparantcor2}
Let $X$ be a CAT(0) cube complex, $C$ a combinatorially convex subspace, $p : X \to C$ the combinatorial projection onto $C$ and $x,y \in X$ two vertices. The hyperplanes separating $p(x)$ and $p(y)$ are precisely the hyperplanes separating $x$ and $y$ which intersect $C$. In particular, $d(p(x),p(y)) \leq d(x,y)$.
\end{prop}

\noindent
\textbf{Proof of Proposition \ref{hyperplanseparantcor2}.} A hyperplane separating $p(x)$ and $p(y)$ separates $x$ and $y$ according to \cite[Lemma 2.10]{coningoff}. Conversely, let $J$ be a hyperplane separating $x$ and $y$ and intersecting $C$. Notice that, according to Lemma \ref{hyperplan séparant} below, if $J$ separates $x$ and $p(x)$, or $y$ and $p(y)$, necessarily $J$ must be disjoint from $C$. Therefore, $J$ has to separate $p(x)$ and $p(y)$. $\square$

\begin{lemma}\label{hyperplan séparant}\emph{\cite[Lemma 2.8]{coningoff}}
Let $X$ be a CAT(0) cube complex and $N \subset X$ a combinatorially convex subspace. Let $p : X \to N$ denote the combinatorial projection onto $N$. Then every hyperplane separating $x$ and $p(x)$ separates $x$ and $N$.
\end{lemma}

\noindent
The following lemma will be particularly useful in this paper.

\begin{lemma}\label{inclusion}
Let $X$ be a CAT(0) cube complex and $C_1 \subset C_2$ two subcomplexes with $C_2$ combinatorially convex. Let $p_2 : X \to C_2$ denote the combinatorial projection onto $C_2$, and $p_1$ the nearest-point projection onto $C_1$, which associates to any vertex the set of the vertices of $C_1$ which minimize the distance from it. Then $p_1 \circ p_2=p_1$. 
\end{lemma}

\noindent
\textbf{Proof.} Let $x \in X$. If $x \in C_2$, clearly $p_1(p_2(x))=p_1(x)$, so we suppose that $x \notin C_2$. If $z \in C_1$, then according to Proposition \ref{projection}, there exists a combinatorial geodesic between $x$ and $z$ passing through $p_2(x)$. Thus, 
\begin{center}
$d(x,z)=d(x,p_2(x))+d(p_2(x),z)$. 
\end{center}
In particular, if $y \in p_1(x)$,
\begin{center}
$d(x,C_1)=d(x,y)=d(x,p_2(x))+d(p_2(x),y)$.
\end{center}
The previous two equalities give:
\begin{center}
$d(p_2(x),z)-d(p_2(x),y)= d(x,z)-d(x,C_1) \geq 0$.
\end{center}
Therefore, $d(p_2(x),y)=d(p_2(x),C_1)$, ie., $y \in p_1(p_2(x))$. We have proved $p_1(x) \subset p_1(p_2(x))$. 

\medskip \noindent
It is worth noticing that we have also proved that $d(x,C_1)=d(x,p_2(x))+d(p_2(x),C_1)$. Thus, for every $y \in p_1(p_2(x))$, once again because there exists a combinatorial geodesic between $x$ and $y$ passing through $p_2(x)$ according to Proposition \ref{projection},
\begin{center}
$d(x,y)=d(x,p_2(x))+d(p_2(x),y)=d(x,p_2(x))+d(p_2(x),C_1)=d(x,C_1)$,
\end{center}
ie., $y \in p_1(x)$. We have proved that $p_1(p_2(x)) \subset p_1(x)$, concluding the proof. $\square$

\medskip \noindent
We conclude with a purely technical lemma which will be used in Section 3.3.

\begin{lemma}\label{d(z,N)}
Let $X$ be a CAT(0) cube complex, $N \subset X$ a combinatorially convex subspace and $x,y \notin N$ two vertices. Fix a combinatorial geodesic $[x,y]$ and choose a vertex $z \in [x,y]$. Then $d(z,N) \leq d(x,N)+d(y,N)$.
\end{lemma}

\noindent
\textbf{Proof.} For convenience, if $A,B \subset X$ are two sets of vertices, let $\Delta(A,B)$ denote the number of hyperplanes separating $A$ and $B$. According to Lemma \ref{hyperplan séparant}, $\Delta(\{a \},N)=d(a,N)$ for every vertex $a \in X$. Notice that, because $z \in [x,y]$ implies that no hyperplane can separate $z$ and $\{x,y\} \cup N$, we have 
\begin{center}
$d(z,N)= \Delta(z,N)= \Delta( \{x,y,z \},N)+ \Delta(\{x,z\},N \cup \{y\})+\Delta(\{y,z \},N \cup \{ x\})$,
\end{center}
and
\begin{center}
$d(x,N)=\Delta(x,N)= \Delta(x, N \cup \{y,z \})+\Delta(\{x,z \}, N \cup \{ y\}) + \Delta(\{x,y,z\} , N )$,
\end{center}
and
\begin{center}
$d(y,N)=\Delta(y,N)= \Delta(y, N \cup \{x,z \})+\Delta(\{y,z \}, N \cup \{ x\}) + \Delta(\{x,y,z\} , N )$.
\end{center}
Now, it is clear that $d(z,N) \leq d(x,N)+d(y,N)$. $\square$

\paragraph{Disc diagrams.} 
A fundamental tool to study CAT(0) cube complexes is the theory of \emph{disc diagrams}. For example, they were extensively used by Sageev in \cite{MR1347406} and by Wise in \cite{Wise}. The rest of this section is dedicated to basic definitions and properties of disc diagrams.

\begin{definition}
Let $X$ be a nonpostively curved cube complex. A \emph{disc diagram} is a continuous combinatorial map $D \to X$, where $D$ is a finite contractible square complex with a fixed topological embedding into $\mathbb{S}^2$; notice that $D$ may be \emph{non-degenerate}, ie., homeomorphic to a disc, or may be \emph{degenerate}. In particular, the complement of $D$ in $\mathbb{S}^2$ is a $2$-cell, whose attaching map will be refered to as the \emph{boundary path} $\partial D \to X$ of $D \to X$; it is a combinatorial path. The \emph{area} of $D \to X$, denoted by $\mathrm{Area}(D)$, corresponds to the number of squares of $D$.
\end{definition}

Given a combinatorial closed path $P \to X$, we say that a disc diagram $D \to X$ is \emph{bounded} by $P \to X$ if there exists an isomorphism $P \to \partial D$ such that the following diagram is commutative:
\begin{displaymath}
\xymatrix{ \partial D \ar[rr] & & X \\ P \ar[u] \ar[urr] & & }
\end{displaymath}
According to a classical argument due to van Kampen \cite{vanKampen} (see also \cite[Lemma 2.17]{McCammondWiseCancellation}), there exists a disc diagram bounded by a given combinatorial closed path if and only if this path is null-homotopic. Thus, if $X$ is a CAT(0) cube complex, then any combinatorial closed path bounds a disc diagram.

\medskip 
As a square complex, a disc diagram contains hyperplanes: they are called \emph{dual curves}. Equivalently, they correspond to the connected components of the reciprocal images of the hyperplanes of $X$. Given a disc diagram $D \to X$, a \emph{nogon} is a dual curve homeomorphic to a circle; a \emph{monogon} is a subpath, of a self-intersecting dual curve, homeomorphic to a circle; an \emph{oscugon} is a subpath of a dual curve whose endpoints are the midpoints of two adjacent edges; a \emph{bigon} is a pair of dual curves intersecting into two different squares.
\begin{figure}[ht!]
\includegraphics[scale=0.4]{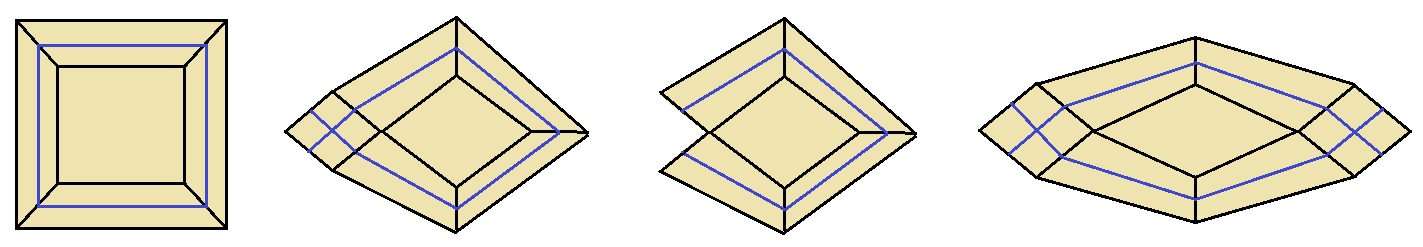}
\caption{From left to right: a nogon, a monogon, an oscugon and a bigon.}
\end{figure}

\begin{thm}\label{Wise}\emph{\cite[Lemma 2.2]{Wise}}
Let $X$ be a nonpositively curved cube complex and $D \to X$ a disc diagram. If $D$ contains a nogon, a monogon, a bigon or an oscugon, then there exists a new disc diagram $D' \to X$ such that:
\begin{itemize}
	\item[(i)] $D'$ is bounded by $\partial D$,
	\item[(ii)] $\mathrm{Area}(D') \leq \mathrm{Area}(D) -2$.
\end{itemize}
\end{thm}

Let $X$ be a CAT(0) cube complex. A \emph{cycle of subcomplexes} $\mathcal{C}$ is a sequence of subcomplexes $C_1, \ldots, C_r$ such that $C_1 \cap C_r \neq \emptyset$ and $C_i \cap C_{i+1} \neq \emptyset$ for every $1 \leq i \leq r-1$. A disc diagram $D \to X$ is \emph{bounded} by $\mathcal{C}$ if $\partial D \to X$ can be written as the concatenation of $r$ combinatorial geodesics $P_1, \ldots, P_r \to X$ such that $P_i \subset C_i$ for every $1 \leq i \leq r$. The \emph{complexity} of such a disc diagram is defined by the couple $(\mathrm{Area}(D), \mathrm{length}(\partial D))$, and a disc diagram bounded by $\mathcal{C}$ will be of \emph{minimal complexity} if its complexity is minimal with respect to the lexicographic order among all the possible disc diagrams bounded by $\mathcal{C}$ (allowing modifications of the paths $P_i$). It is worth noticing that such a disc diagram does not exist if our subcomplexes contain no combinatorial geodesics. On the other hand, if our subcomplexes are combinatorially geodesic, then a disc diagram always exists.

\medskip 
Our main technical result on disc diagrams is the following, which we proved in \cite[Theorem 2.13]{coningoff}.

\begin{thm}\label{disc diagram}
Let $X$ be a CAT(0) cube complex, $\mathcal{C}=(C_1, \ldots, C_r)$ a cycle of subcomplexes, and $D \to X$ a disc diagram bounded by $\mathcal{C}$. For convenience, write $\partial D$ as the concatenation of $r$ combinatorial geodesics $P_1, \ldots, P_r \to X$ with $P_i \subset C_i$ for every $1 \leq i\leq r$. If the complexity of $D \to X$ is minimal, then:
\begin{itemize}
	\item[(i)] if $C_i$ is combinatorially convex, two dual curves intersecting $P_i$ are disjoint;
	\item[(ii)] if $C_i$ and $C_{i+1}$ are combinatorially convex, no dual curve intersects both $P_i$ and $P_{i+1}$.
\end{itemize}
\end{thm}

\noindent
In general, a disc diagram $D \to X$ is not an embedding. However, the proposition below, which we proved in \cite[Proposition 2.15]{coningoff}, determines precisely when it is an isometric embedding.

\begin{prop}\label{disc embedding}
Let $X$ be a CAT(0) cube complex and $D \to X$ a disc diagram which does not contain any bigon. With respect to the combinatorial metrics, $\varphi : D \to X$ is an isometric embedding if and only if every hyperplane of $X$ induces at most one dual curve of $D$. 
\end{prop}

\noindent
We mention two particular cases below.

\begin{cor}\label{embedding}
Let $\mathcal{C}=(C_1,C_2,C_3,C_4)$ be a cycle of four subcomplexes. Suppose that $C_2,C_3,C_4$ are combinatorially convex subcomplexes. If $D \to X$ is a disc diagram of minimal complexity bounded by $\mathcal{C}$, then $D \to X$ is an isometric embedding.
\end{cor}

\noindent
\textbf{Proof.} First, we write $\partial D \to X$ as the concatenation of four combinatorial geodesics $P_1,P_2,P_3,P_4 \to X$, where $P_i \subset C_i$ for every $1 \leq i \leq 4$. Suppose that there exist two dual curves $c_1,c_2$ of $D$ induced by the same hyperplane $J$. Because a combinatorial geodesic cannot intersect a hyperplane twice, the four endpoints of $c_1$ and $c_2$ belong to four different sides of $\partial D$. On the other hand, because $C_2,C_3,C_4$ are combinatorially convex, it follows from Theorem \ref{disc diagram} that any dual curve intersecting $P_3$ must intersect $P_1$; say that $c_1$ intersects $P_1$ and $P_3$. Therefore, $c_2$ must intersect $P_2$ and $P_4$, and we deduce that $c_1$ and $c_2$ are transverse. But this implies that $J$ self-intersects in $X$, which is impossible. $\square$

\begin{cor}\label{flat rectangle}\emph{\cite[Corollary 2.17]{coningoff}}
Let $X$ be a CAT(0) cube complex and $\mathcal{C}$ a cycle of four combinatorially convex subcomplexes. If $D \to X$ is a disc diagram of minimal complexity bounded by $\mathcal{C}$, then $D$ is combinatorially isometric to a rectangle $[0,a] \times [0,b] \subset \mathbb{R}^2$ and $D \to X$ is an isometric embedding.
\end{cor}

\paragraph{Contracting isometries.} 
We gave in the introduction the definition of a contracting isometry. It turns out that an isometry is contracting if and only if its axis is contracting, in the following sense:

\begin{definition}
Let $(S,d)$ be a metric space, $Y \subset C$ a subspace and $B \geq 0$ a constant. We say that $Y$ is \emph{$B$-contracting} if, for every ball disjoint from $Y$, the diameter of its nearest-point projection onto $Y$ is at most $B$. A subspace will be \emph{contracting} if it is $B$-contracting for some $B \geq 0$.
\end{definition}

\noindent
In fact, in the context of CAT(0) cube complexes, we may slightly modify this definition in the following way:

\begin{lemma}\label{contractant critère}
Let $X$ be a CAT(0) cube complex, $S \subset X$ a combinatorially convex subspace and $L \geq 0$ a constant. Then $S$ is contracting if and only if there exists $C \geq 0$ such that, for all vertices $x,y \in X$ satisfying $d(x,y)<d(x,S)-L$, the projection of $\{x,y \}$ onto $S$ has diameter at most $C$.
\end{lemma}

\noindent
\textbf{Proof.} The implication is clear: if $S$ is $B$-contracting, $\{x,y \}$ is included into the ball $B(x,d(x,S))$ whose projection onto $S$ has diameter at most $B$.

\medskip \noindent
Conversely, suppose that $S$ is not contracting, ie., for all $n \geq 0$, there exists a ball $B(x_n,r_n)$, with $r_n < d(x_n,S)$, whose projection onto $S$ has diameter at least $n$. Let $p : X \to S$ denote the combinatorial projection onto $S$. If for all $y \in B(x_n,r_n)$ we had $d(p(y),p(x_n)) < n/2$, then the projection of $B(x_n,r_n)$ onto $S$ would have diameter at most $n$. Therefore, there exists $y_n \in B(x_n,r_n)$ satisfying $d(p(x_n),p(y_n)) \geq n/2$. In particular, the projection of $\{ x_n,y_n \}$ onto $S$ has diameter at least $n/2$, with $d(x_n,y_n)<d(x_n,S)$. 

\medskip \noindent
If $d(x_n,y_n) < d(x_n,S)-L$, there is nothing to prove. Ohterwise, two cases may happen. If $d(x_n,S) \leq L$, then, because $p$ is $1$-Lipschitz according to Proposition \ref{hyperplanseparantcor2}, the projection of $B(x_n,r_n)$ onto $S$ has diameter at most $2L$. Without loss of generality, we may suppose $n > 2L$, so that this case cannot occur. From now on, suppose $d(x_n,y_n) \geq d(x_n,S)-L$. Then, along a combinatorial geodesic $[x_n,y_n]$ between $x_n$ and $y_n$, there exists a vertex $z_n$ such that $d(z_n,y_n) \leq L$ and $d(x_n,z_n) \leq d(x_n,S)-L$ because we saw that $d(x_n,y_n)<d(x_n,S)$. Now,
\begin{center}
$d(p(x_n),p(z_n)) \geq d(p(x_n),p(y_n))- d(p(y_n),p(z_n)) \geq \frac{n}{2}-L$,
\end{center}
so that the projection of $\{x_n,z_n \}$ onto $S$ has diameter at least $\frac{n}{2}-L$. Because this diameter may be chosen arbitrarily large, this concludes the proof. $\square$

\medskip \noindent
However, the previous criterion only holds for combinatorially convex subspaces. Of course, it is always possible to look at the combinatorial convex hull of the subspace which is considered, and the following lemma essentially states when it is possible to deduce that the initial subspace is contracting.

\begin{lemma}\label{convexhull}
Let $X$ be a CAT(0) cube complex, $S \subset X$ a subspace and $N \supset S$ a combinatorially convex subspace. Suppose that the Hausdorff distance $d_H(S,N)$ between $S$ and $N$ is finite. Then $S$ is contracting if and only if $N$ is contracting as well.
\end{lemma}

\noindent
\textbf{Proof.} We know from Lemma \ref{inclusion} that the projection onto $S$ is the composition of the projection of $X$ onto $N$ with the projection of $N$ onto $S$. Consequently, the Hausdorff distance between the projections onto $N$ and onto $S$ is at most $d_H(N,S)$. Our lemma follows immediately. $\square$

\medskip \noindent
Typically, the previous lemma applies to quasiconvex combinatorial geodesics:

\begin{lemma}\label{d_H finie}
Let $X$ be a CAT(0) cube complex and $\gamma$ an infinite combinatorial geodesic. Let $N(\gamma)$ denote the combinatorial convex hull of $\gamma$. Then $\gamma$ is quasiconvex if and only if the Hausdorff distance between $\gamma$ and $N(\gamma)$ is finite.
\end{lemma}

\noindent
\textbf{Proof.} We claim that (the 1-skeleton of) $N(\gamma)$ is the union $\mathcal{L}$ of the combinatorial geodesics whose endpoints are on $\gamma$. First, it is clear that any such geodesic must be included into $N(\gamma)$, whence $\mathcal{L} \subset N(\gamma)$. Conversely, given some $x \in N(\gamma)$, we want to prove that there exists a geodesic between two vertices of $\gamma$ passing through $x$. Fix some $y \in \gamma$. Notice that, because $x$ belongs to $N(\gamma)$, no hyperplane separates $x$ from $\gamma$, so that any hyperplane separating $x$ and $y$ must intersect $\gamma$. As a consequence, for some $n \geq 1$ sufficiently large, every hyperplane separating $x$ and $y$ separates $\gamma(-n)$ and $\gamma(n)$; furthermore, we may suppose without loss of generality that $y$ belongs to the subsegment of $\gamma$ between $\gamma(-n)$ and $\gamma(n)$. We claim that the concatenation of a geodesic from $\gamma(-n)$ to $x$ with a geodesic from $x$ to $\gamma(n)$ defines a geodesic between $\gamma(-n)$ and $\gamma(n)$, which proves that $x \in \mathcal{L}$. Indeed, if this concatenation defines a path which is not a geodesic, then there must exist a hyperplane $J$ intersecting it twice. Notice that $J$ must separate $x$ from $\{\gamma(n), \gamma(-n) \}$, so that we deduce from the convexity of half-spaces that $J$ must separate $x$ and $y$. By our choice of $n$, we deduce that $J$ separates $\gamma(-n)$ and $\gamma(n)$. Finally, we obtain a contradiction since we already know that $J$ separates $x$ from $\gamma(-n)$ and $\gamma(n)$. Thus, we have proved that $N(\gamma) \subset \mathcal{L}$. 

\medskip \noindent
It follows that $\gamma$ is quasiconvex if and only if the Hausdorff distance between $\gamma$ and $N(\gamma)$ is finite. $\square$

\medskip \noindent
It is worth noticing that any contracting geodesic is quasiconvex. This is a particular case of \cite[Lemma 3.3]{MR3175245}; we include a proof for completeness.

\begin{lemma}\label{Sultan}
Let $X$ be a geodesic metric space. Any contracting geodesic of $X$ is quasiconvex.
\end{lemma}

\noindent
\textbf{Proof.} Let $\varphi$ be a geodesic between two points of $\gamma$. Let $\varphi(t)$ be a point of $\varphi$. We want to prove that $d(\varphi(t),\gamma) \leq 11B$.

\medskip \noindent
If $d(\varphi(t), \gamma) < 3B$, there is nothing to prove. Otherwise, $\varphi(t)$ belongs to a maximal subsegment $[\varphi(r), \varphi(s)]$ of $\varphi$ outside the $3B$-neighborhood of $\gamma$, ie. $d(\varphi(r),\gamma),d(\varphi(s),\gamma) \leq 3B$ but $d(\varphi(p),\gamma) \geq 3B$ for every $p \in [r,s]$. Let
\begin{equation} \label{eq:S1}
r=t_0<t_1< \cdots < t_{k-1} < t_k = s
\end{equation}
with $t_{i+1}-t_i=2B$ for $1 \leq i \leq k-2$ and $t_k-t_{k-1} \leq 2B$. Observe that
\begin{center}
$\displaystyle d(\varphi(r),\varphi(s)) = s-r= \sum\limits_{i=0}^{k-1} (t_{i+1}-t_i) \geq 2(k-1)B$,
\end{center}
but on the other hand, if $p_i$ belongs to the projection of $\varphi(t_i)$ on $\gamma$,
\begin{center}
$d(\varphi(r),\varphi(s)) \leq d(\varphi(r),p_0)+d(p_0,p_1)+ \cdots + d(p_{k-1},p_k)+ d(p_k,\varphi(s))$.
\end{center}
Noticing that $d(\varphi(t_i),\varphi(t_{i+1}))= t_{i+1}-t_i \leq 2B<3B \leq d(\varphi(t_i),\gamma)$, we deduce from the fact that $\gamma$ is $B$-contracting that
\begin{equation} \label{eq:S2}
d(\varphi(r), \varphi(s)) \leq 6B+kB.
\end{equation}
Finally, combining \ref{eq:S1} and \ref{eq:S2}, we get $k \leq 8$. Thus,
\begin{center}
$\displaystyle d(\varphi(r),\varphi(s)) = \sum\limits_{i=0}^{k-1} (t_{i+1}-t_i) \leq 2kB \leq 16B$.
\end{center}
Now, since $d(\varphi(t),\varphi(r))$ or $d(\varphi(t),\varphi(s))$, say the first one, is bounded above by $d(\varphi(r),\varphi(s))/2 \leq 8B$, we deduce
\begin{center}
$d(\varphi(t),\gamma) \leq d(\varphi(t),\varphi(r))+ d(\varphi(r),\gamma) \leq 8B+3B=11B$.
\end{center}
This concludes the proof. $\square$

\section{Contracting isometries from hyperplane configurations}

\subsection{Quasiconvex geodesics I}

\noindent
Given a subcomplex $Y$ of some CAT(0) cube complex, we denote by $\mathcal{H}(Y)$ the set of hyperplanes intersecting $Y$. In this section, the question we are interested in is to determine when a combinatorial geodesic $\gamma$ is quasiconvex just from the set $\mathcal{H}(\gamma)$. We begin by introducing some definitions.

\begin{definition}
A \emph{facing triple} in a CAT(0) cube complex is the data of three hyperplanes such that no two of them are separated by the third hyperplane.
\end{definition}

\begin{definition}
A \emph{join of hyperplanes} is the data of two families $(\mathcal{H}=(H_1, \ldots, H_r),\mathcal{V}=(V_1, \ldots, V_s))$ of hyperplanes which do not contain facing triple and such that any hyperplane of $\mathcal{H}$ is transverse to any hyperplane of $\mathcal{V}$. If $H_i$ (resp. $V_j$) separates $H_{i-1}$ and $H_{i+1}$ (resp. $V_{j-1}$ and $V_{j+1}$) for every $2 \leq i \leq r-1$ (resp. $2 \leq j \leq s-1$), we say that $(\mathcal{H},\mathcal{V})$ is a \emph{grid of hyperplanes}.\\
If $(\mathcal{H}, \mathcal{V})$ is a join or a grid of hyperplanes satisfying $\min( \# \mathcal{H}, \# \mathcal{V}) \leq C$, we say that $(\mathcal{H}, \mathcal{V})$ is \emph{$C$-thin}.
\end{definition}

\noindent
Our main criterion is:

\begin{prop}\label{prop-qc}
Let $X$ be a CAT(0) cube complex and $\gamma$ an infinite combinatorial geodesic. The following statements are equivalent:  
\begin{itemize}
	\item[(i)] $\gamma$ is quasiconvex;
	\item[(ii)] there exists $C \geq 1$ satisfying:
	\begin{itemize}
		\item $\mathcal{H}(\gamma)$ does not contain $C$ pairwise transverse hyperplanes,
		\item any grid of hyperplanes $(\mathcal{H},\mathcal{V})$ with $\mathcal{H}, \mathcal{V} \subset \mathcal{H}(\gamma)$ is $C$-thin;
	\end{itemize}
	\item[(iii)] there exists some constant $C \geq 1$ such that any join of hyperplanes $(\mathcal{H},\mathcal{V})$ satisfying $\mathcal{H}, \mathcal{V} \subset \mathcal{H}(\gamma)$ is $C$-thin.
\end{itemize}
\end{prop}

\noindent
It will be a direct consequence of the following lemma and a criterion of hyperbolicity we proved in \cite{coningoff}:

\begin{lemma}\label{lem-qc}
An infinite combinatorial geodesic $\gamma$ is quasiconvex if and only if its combinatorial convex hull is a hyperbolic subcomplex.
\end{lemma}

\noindent
\textbf{Proof.} According to Lemma \ref{d_H finie}, $\gamma$ is quasiconvex if and only if the Hausdorff distance between $\gamma$ and $N(\gamma)$ is finite. In particular, if $\gamma$ is quasiconvex, then $N(\gamma)$ is quasi-isometric to a line, and a fortiori is hyperbolic. Conversely, if $N(\gamma)$ is a hyperbolic subcomplex, its bigons are uniformly thin. This implies the quasiconvexity of $\gamma$. $\square$

\begin{thm}\label{critère d'hyperbolicité}\emph{\cite[Theorem 3.3]{coningoff}}
A CAT(0) cube complex is hyperbolic if and only if it is finite-dimensional and its grids of hyperplanes of $X$ are uniformly thin.
\end{thm}

\noindent
\textbf{Proof of Proposition \ref{prop-qc}.} From Lemma \ref{lem-qc} and Theorem \ref{critère d'hyperbolicité}, we deduce the equivalence $(i) \Leftrightarrow (ii)$. Then, because a grid of hyperplanes or a collection of pairwise transverse hyperplanes gives rise to a join of hyperplanes, $(iii)$ clearly implies $(ii)$. 

\medskip \noindent
To conclude, we want to prove $(ii) \Rightarrow (iii)$. Let $C$ denote the constant given by $(ii)$. Let $(\mathcal{H},\mathcal{V})$ be a join of hyperplanes satisfying $\mathcal{H}, \mathcal{V} \subset \mathcal{H}(\gamma)$. Suppose by contradiction that $\# \mathcal{H},\# \mathcal{V} \geq \mathrm{Ram}(C)$. Because $\mathcal{H}(\gamma)$ does not contain $C$ pairwise transverse hyperplanes, we deduce that $\mathcal{H}$ and $\mathcal{V}$ each contain a subfamily of $C$ pairwise disjoint hyperplanes, say $\mathcal{H}'$ and $\mathcal{V}'$ respectively. Notice that any hyperplane of $\mathcal{H}'$ is transverse to any hyperplane of $\mathcal{V}'$. Since the hyperplanes of $\mathcal{H}'$ and $\mathcal{V}'$ are all intersected by $\gamma$, we conclude that $(\mathcal{H}', \mathcal{V}')$ defines a $(C,C)$-grid of hyperplanes in $X$, contradicting the definition of $C$. Therefore, $\min(p,q) < \mathrm{Ram}(C)$. $\square$

\subsection{Contracting convex subcomplexes I}

\noindent
In the previous section, we showed how to recognize the quasiconvexity of a combinatorial geodesic from the set of the hyperplanes it intersects. Therefore, according to Proposition \ref{prop-qc}, we reduced the problem of determining when a combinatorial geodesic is contracting to the problem of determining when a combinatorially convex subcomplex is contracting. This is the main criterion of this section:

\begin{thm}\label{contracting subcomplexes}
Let $X$ be a CAT(0) cube complex and $Y \subset X$ a combinatorially convex subcomplex. Then $Y$ is contracting if and only if there exists $C \geq 0$ such that any join of hyperplanes $(\mathcal{H},\mathcal{V})$ with $\mathcal{H} \subset \mathcal{H}(Y)$ and $\mathcal{V} \cap \mathcal{H}(Y) = \emptyset$ satisfies $\min ( \# \mathcal{H}, \# \mathcal{V}) \leq C$.
\end{thm}

\noindent
\textbf{Proof.} Suppose that $Y$ is $B$-contracting for some $B \geq 0$, and let $(\mathcal{H},\mathcal{V})$ be a join of hyperplanes satisfying $\mathcal{H} \subset \mathcal{H}(Y)$ and $\mathcal{V} \cap \mathcal{H}(Y) = \emptyset$. Because any hyperplane of $\mathcal{H}$ intersects $Y$ and $\mathcal{H}$ does not contain facing triples, there exist two vertices $x_-,x_+ \in Y$ separated by the hyperplanes of $\mathcal{H}$. For convenience, let $H^{\pm}$ denote the halfspaces delimited by some $H \in \mathcal{H}$ which contains $x_{\pm}$; and, for any $V \in \mathcal{V}$, let $V^+$ denote the halfspace delimited by $V$ which does not contain $Y$. Now, we set
\begin{center}
$\mathfrak{H}^{\pm} = \bigcap\limits_{H \in \mathcal{H}} H^{\pm}$ and $\mathfrak{V}= \bigcap\limits_{V \in \mathcal{V}} V^+$.
\end{center}
If $\mathcal{C}$ is the cycle of subcomplexes $(\mathfrak{H}^-,\mathfrak{V},\mathfrak{H}^+,Y)$, then Corollary \ref{flat rectangle} states that a disc diagram of minimal complexity bounded by $\mathcal{C}$ defines a flat rectangle $F \subset X$. Let $x,y \in F \cap \mathfrak{V}$ be two vertices satisfying $d(x,y)= \min(\# \mathcal{V},\# \mathcal{H})-1$. Then, because $Y$ and $\{x,y\}$ are separated by the hyperplanes of $\mathcal{V}$, we have $d(x,y) < d(x,Y)$. Moreover, if $p : X \to Y$ denotes the combinatorial projection onto $Y$, because $x$ and $y$ are separated by $\min(\# \mathcal{H},\# \mathcal{V})-1$ hyperplanes which intersect $Y$, it follows from Proposition \ref{hyperplanseparantcor2} that $d(p(x),p(y)) \geq \min(\#\mathcal{H},\# \mathcal{V})-1$. Finally, since $Y$ is $B$-contracting, we conclude that
\begin{center}
$ \min ( \# \mathcal{H}, \# \mathcal{V}) \leq 1+d(p(x), p(y)) \leq 1+B$.
\end{center}
Conversely, suppose that there exists $C \geq 0$ such that any join of hyperplanes $(\mathcal{H},\mathcal{V})$ with $\mathcal{H} \subset \mathcal{H}(Y)$ and $\mathcal{V} \cap \mathcal{H}(Y) = \emptyset$ satisfies $\min ( \# \mathcal{H}, \# \mathcal{V}) \leq C$. We want to prove that $Y$ is contracting by applying Lemma \ref{contractant critère}. Let $x,y \in X$ be two distinct vertices satisfying $d(x,y)<d(x,Y)-C$ and let $p : X \to Y$ denote the combinatorial projection onto $Y$. If $\mathcal{H}$ denotes the set of the hyperplanes separating $\{x,y\}$ and $Y$, and $\mathcal{V}$ the set of the hyperplanes separating $p(x)$ and $p(y)$, then $(\mathcal{H},\mathcal{V})$ defines a join of hyperplanes: because any hyperplane of $\mathcal{H}$ separates $\{x,y\}$ and $\{ p(x),p(y) \}$, and hyperplane of $\mathcal{V}$ separates $\{x,p(x) \}$ and $\{ y,p(y) \}$ according to Proposition \ref{hyperplanseparantcor2} and Lemma \ref{hyperplan séparant}, we deduce that any hyperplane of $\mathcal{H}$ is transverse to any hyperplane of $\mathcal{V}$. Because the assumption $d(x,y)<d(x,Y)-C$ implies $\# \mathcal{H} > C$, we conclude that
\begin{center}
$d(p(x),p(y)) \leq \# \mathcal{V} \leq C$.
\end{center}
Therefore, $Y$ is contracting. $\square$

\medskip \noindent
Combining Proposition \ref{prop-qc} and Theorem \ref{contracting subcomplexes}, we get:

\begin{cor}\label{géodésiques contractantes}
Let $X$ be a CAT(0) cube complex and $\gamma$ an infinite combinatorial geodesic. Then $\gamma$ is contracting if and only if there exists a constant $C \geq 1$ such that any join of hyperplanes $(\mathcal{H},\mathcal{V})$ satisfying $\mathcal{H} \subset \mathcal{H}(\gamma)$ is necessarily $C$-thin.
\end{cor}

\noindent
\textbf{Proof.} Suppose that $\gamma$ is contracting. In particular, $\gamma$ is quasiconvex so, according to Proposition \ref{prop-qc}, there exists a constant $C_1$ such that any join of hyperplanes $(\mathcal{H},\mathcal{V})$ satisfying $\mathcal{H}, \mathcal{V} \subset \mathcal{H}(\gamma)$ is $C_1$-thin. Then, because $N(\gamma)$ is also contracting according to Lemma \ref{d_H finie} and Lemma \ref{convexhull}, we deduce from Theorem \ref{contracting subcomplexes} that there exists a constant $C_2$ such that any join of hyperplanes $(\mathcal{H}, \mathcal{V})$ satisfying $\mathcal{H} \subset \mathcal{H}(\gamma)$ and $\mathcal{V} \cap \mathcal{H}(\gamma)=\emptyset$ is $C_2$-thin. Now, let $(\mathcal{H},\mathcal{V})$ be any join of hyperplanes satisfying $\mathcal{H} \subset \mathcal{H}(\gamma)$. If we set $\mathcal{V}_1= \mathcal{V} \cap \mathcal{H}(\gamma)$ and $\mathcal{V}_2= \mathcal{V}\backslash \mathcal{V}_1$, then $(\mathcal{H},\mathcal{V}_1)$ is $C_1$-thin and $(\mathcal{H},\mathcal{V}_2)$ is $C_2$-thin. Thus,
\begin{center}
$\min( \# \mathcal{H}, \#\mathcal{V}) \leq \min( \# \mathcal{H}, \#\mathcal{V}_1) + \min( \# \mathcal{H}, \#\mathcal{V}_2) \leq C_1+C_2$,
\end{center}
ie., $(\mathcal{H}, \mathcal{V})$ is $(C_1+C_2)$-thin.

\medskip \noindent
Conversely, suppose that there exists a constant $C \geq 1$ such that any join of hyperplanes $(\mathcal{H},\mathcal{V})$ satisfying $\mathcal{H} \subset \mathcal{H}(\gamma)$ is necessarily $C$-thin. According to Proposition \ref{prop-qc}, $\gamma$ is quasiconvex, so that it is sufficient to prove that $N(\gamma)$ is contracting to conclude that $\gamma$ is contracting as well according to Lemma \ref{convexhull}. Finally, it follows from Theorem \ref{contracting subcomplexes} that $N(\gamma)$ is contracting. $\square$

\subsection{Well-separated hyperplanes}

\noindent
Separation properties of hyperplanes play a central role in the study of contracting isometries. \emph{Strongly separated} hyperplanes were introduced in \cite{BehrstockCharney} in order to characterize the rank-one isometries of right-angled Artin groups; they were also crucial in the proof of the Rank Rigidity Theorem for CAT(0) cube complexes in \cite{MR2827012}. In \cite{MR3339446}, Charney and Sultan used contracting isometries to distinguish quasi-isometrically some cubulable groups, and they introduced \emph{$k$-separated} hyperplanes in order to characterize contracting isometries in uniformly locally finite CAT(0) cube complexes \cite[Theorem 4.2]{MR3339446}. In this section, we introduce \emph{well-separated} hyperplanes in order to generalize this charaterization to CAT(0) cube complexes without any local finiteness condition.

\begin{definition}
Let $J$ and $B$ be two disjoint hyperplanes in some CAT(0) cube complex and $L \geq 0$. We say that $J$ and $H$ are \textit{$L$-well separated} if any family of hyperplanes transverse to both $J$ and $H$, which does not contain any facing triple, has cardinality at most $L$. Two hyperplanes are \textit{well-separated} if they are $L$-well-separated for some $L \geq 1$.
\end{definition}

\begin{thm}\label{géodésiques contractantes 2}
Let $\gamma$ be an infinite combinatorial geodesic. Then $\gamma$ is contracting if and only if there exist constants $r,L \geq 1$, hyperplanes $\{ H_i,i \in \mathbb{Z} \}$, and vertices $\{x_i \in \gamma \cap N(H_i), i \in \mathbb{Z} \}$ such that, for every $i \in \mathbb{Z}$:
\begin{itemize}
	\item[$\bullet$] $d(x_i,x_{i+1}) \leq r$,
	\item[$\bullet$] the hyperplanes $H_i$ and $H_{i+1}$ are $L$-well separated.
\end{itemize}
\end{thm}

\noindent
The following lemma is a combinatorial analogue to \cite[Lemma 4.3]{MR3339446}, although our proof is essentially different. This is the key technical lemma needed to prove Proposition \ref{géodésiques quasi-convexes 2}, which is in turn the first step toward the proof of Theorem \ref{géodésiques contractantes 2}. 

\begin{lemma}
Let $\gamma$ be an infinite quasiconvex combinatorial geodesic. There exists a constant $C$ depending only on $\gamma$ such that, if two vertices $x ,y \in \gamma$ satisfy $d(x,y)> C$, then there exists a hyperplane $J$ separating $x$ and $y$ such that the projection of $N(J)$ onto $\gamma$ is included into $[x,y] \subset \gamma$.
\end{lemma}

\noindent
\textbf{Proof.} Let $J$ be a hyperplane separating $x$ and $y$ whose projection onto $\gamma$ is not included into $[x,y] \subset \gamma$. Say that this projection contains a vertex at the right of $y$; the case where it contains a vertex at the left of $x$ is completely symmetric. Thus, we have the following configuration:
\begin{center}
\includegraphics[scale=0.6]{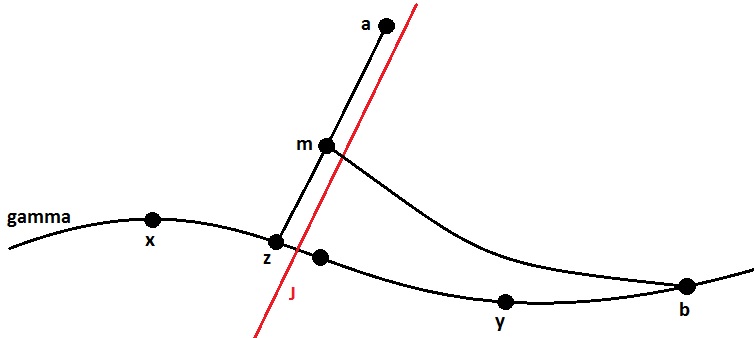}
\end{center}
where $b \notin [x,y]$ is a projection onto $\gamma$ of some vertex $a \in N(J)$, $z \in \gamma$ is a vertex adjacent to $J$, and $m=m(a,b,z)$ the median point of $\{a,b,z \}$. Because $m$ belongs to a combinatorial geodesic between $z$ and $b$, and that these two points belong to $\gamma$, necessarily $m \in N(\gamma)$. On the other hand, $m$ belongs to a combinatorial geodesic between $a$ and $b$, and by definition $b$ is a vertex of $\gamma$ minimizing the distance to $a$, so $d(m,b)=d(m, \gamma)$; because, by the quasiconvexity of $\gamma$, the Hausdorff distance between $\gamma$ and $N(\gamma)$ is finite, say $L$, we deduce that $d(m,b)=d(m,\gamma) \leq L$ since $m$ belongs to a combinatorial geodesic between $b$ and $z$ (which implies that $m \in N(\gamma)$). Using Lemma \ref{d(z,N)}, we get
\begin{center}
$d(y,N(J)) \leq d(z,N(J)) + d(b,N(J)) \leq d(b,m) \leq L$.
\end{center}
Thus, $J$ intersects the ball $B(y,L)$. 

\medskip \noindent
Let $\mathcal{H}$ be the set of hyperplanes separating $x$ and $y$, and intersecting $B(y,L)$. Of course, because any hyperplane of $\mathcal{H}$ intersects $\gamma$, the collection $\mathcal{H}$ contains no facing triple. On the other hand, it follows from Proposition \ref{prop-qc} that $\dim N(\gamma)$ is finite, so that if $\# \mathcal{H} \geq \mathrm{Ram}(s)$ for some $s \geq \dim N(\gamma)$ then the collection $\mathcal{H}$ must contain at least $s$ pairwise disjoint hyperplanes. Since this collection of hyperplanes does not contain facing triples and intersects the ball $B(y,L)$, we deduce that $s \leq 2L$, hence $\# \mathcal{H} \leq \mathrm{Ram}( \max( \dim N(\gamma),2L))$. 

\medskip \noindent
Finally, we have proved that there are at most $2 \mathrm{Ram}( \max( \dim N(\gamma),2L))$ hyperplanes separating $x$ and $y$ such that the projections of their neighborhoods onto $\gamma$ are not included into $[x,y] \subset \gamma$. It clearly follows that $d(x,y) > 2\mathrm{Ram}( \max( \dim N(\gamma),2L))$ implies that there is a hyperplane separating $x$ and $y$ such that the projection of its neighborhood onto $\gamma$ is included into $[x,y] \subset \gamma$. $\square$

\begin{prop}\label{géodésiques quasi-convexes 2}
Let $\gamma$ be an infinite combinatorial geodesic. If $\gamma$ is quasiconvex then there exist constants $r,L \geq 1$, hyperplanes $\{ H_i,i \in \mathbb{Z} \}$, and vertices $\{x_i \in \gamma \cap N(H_i), i \in \mathbb{Z} \}$ such that, for every $i \in \mathbb{Z}$:
\begin{itemize}
	\item[$\bullet$] $d(x_i,x_{i+1}) \leq r$,
	\item[$\bullet$] the hyperplanes $H_i$ and $H_{i+1}$ are disjoint.
\end{itemize}
\end{prop}

\noindent
\textbf{Proof.} Because $\gamma$ is quasiconvex, we may apply the previous lemma: let $C$ be the constant it gives. Let $\ldots, y_{-1},y_0,y_1 , \ldots \in \gamma$ be vertices along $\gamma$ satisfying $d(y_i,y_{i+1})= C+1$ for all $i \in \mathbb{Z}$. According to the previous lemma, for every $k \in \mathbb{Z}$, there exists a hyperplane $J_k$ separating $y_{2k}$ and $y_{2k+1}$ whose projection onto $\gamma$ is included into $[y_{2k},y_{2k+1}] \subset \gamma$; let $x_k$ be one of the two vertices in $\gamma \cap N(J_k)$. Notice that
\begin{center}
$d(x_i,x_{i+1}) \leq d(y_{2i},y_{2i+3}) \leq d(y_{2i},y_{2i+1}) + d(y_{2i+1},y_{2i+2}) + d(y_{2i+2},y_{2i+3}) = 3(C+1)$.
\end{center}
Finally, it is clear that, for every $k \in \mathbb{Z}$, the hyperplanes $J_k$ and $J_{k+1}$ are disjoint: if it was not the case, there would exist a vertex of $N(J_k) \cap N(J_{k+1})$ whose projection onto $\gamma$ would belong to $[y_{2k},y_{2k+1}] \cap [y_{2k+2},y_{2k+3}]= \emptyset$. $\square$

\medskip \noindent
\textbf{Proof of Theorem \ref{géodésiques contractantes 2}.} Suppose $\gamma$ contracting. In particular, $\gamma$ is quasiconvex, so, by applying Proposition \ref{géodésiques quasi-convexes 2}, we find a constant $L \geq 1$, a collection of pairwise disjoint hyperplanes $\{ J_i, i \in \mathbb{Z} \}$, and a collection of vertices $\{ y_i \in  \gamma \cap N(J_i), i \in \mathbb{Z} \}$, such that $d(y_i,y_{i+1}) \leq L$ for all $i \in \mathbb{Z}$. Let $C$ be the constant given by Corollary \ref{géodésiques contractantes} and set $x_i=y_{iC}$ for all $i \in \mathbb{Z}$. Notice that
\begin{center}
$d(x_i,x_{i+1})= d(y_{iC},y_{iC+C}) \leq d(y_{iC},y_{iC+1})+ \cdots + d(y_{iC+C-1},y_{iC+C}) \leq (C+1)L$.
\end{center}
Now, we want to prove that the hyperplanes $J_{nC}$ and $J_{(n+1)C}$ are $C$-well separated for every $n \in \mathbb{Z}$. So let $\mathcal{H}$ be a collection of hyperplanes transverse to both $J_{nC}$ and $J_{(n+1)C}$, which contains no facing triple. Because every hyperplane $H \in \mathcal{H}$ is transverse to any $J_{nC+k}$, for $0 \leq k \leq C$, we obtain a join of hyperplanes $(V_{C+1}, V_{\# \mathcal{H}})$ satisfying $V_{C+1} \subset \mathcal{H}(\gamma)$. By definition of $C$, we deduce $\# \mathcal{H} \leq C$.

\medskip \noindent
Conversely, suppose that there exist constants $\ell,L \geq 1$, hyperplanes $\{ J_i , i \in \mathbb{Z} \}$, and vertices $\{ x_i \in \gamma \cap N(J_i) , i \in \mathbb{Z} \}$ such that, for every $i \in \mathbb{Z}$:
\begin{itemize}
	\item[$\bullet$] $d(x_i,x_{i+1}) \leq \ell$,
	\item[$\bullet$] the hyperplanes $J_i$ and $J_{i+1}$ are $L$-well separated.
\end{itemize}
Let $(\mathcal{V}_p , \mathcal{V}_q)$ be a join of hyperplanes with $\mathcal{V}_p = \{ V_1, \ldots, V_p\} \subset \mathcal{H}(\gamma)$. For convenience, we may suppose without loss of generality that each $V_i$ intersects $\gamma$ at $y_i$ with the property that $y_i$ separates $y_{i-1}$ and $y_{i+1}$ along $\gamma$ for all $2 \leq i \leq p-1$. If $p > 3 \ell+2L$, there exist $L< r < s < p-L $ such that $y_r$ and $y_s$ are separated by $J_k$, $J_{k+1}$, $J_{k+2}$ and $J_{k+3}$ for some $k \in \mathbb{Z}$. Because $J_k$ and $J_{k+1}$ are $L$-well separated, the hyperplanes $\{ V_1, \ldots, V_r \}$ cannot be transverse to both $J_k$ and $J_{k+1}$ since $r >L$, so there exists $1 \leq \alpha \leq r$ such that $V_{\alpha}$ and $J_{k+1}$ are disjoint. Similarly, the hyperplanes $\{ V_s, \ldots, V_p \}$ cannot be transverse to both $J_{k+2}$ and $J_{k+3}$ since $p-s >L$, so there exists $1 \leq \omega \leq r$ such that $V_{\omega}$ and $J_{k+2}$ are disjoint. Notice that the hyperplanes $J_{k+1}$ and $J_{k+2}$ separate $V_{\alpha}$ and $V_{\omega}$, so that the hyperplanes of $\mathcal{V}_q$, which are all transverse to both $V_{\alpha}$ and $V_{\omega}$, are transverse to both $J_{k+1}$ and $J_{k+2}$. But $J_{k+1}$ and $J_{k+2}$ are $L$-well separated, hence $q \leq L$. 

\medskip \noindent
Thus, we have proved that $\min(p,q) \leq 3 \ell +2L$. Finally, Corollary \ref{géodésiques contractantes} implies that $\gamma$ is contracting. $\square$

\subsection{Contracting isometries I}

\noindent
Finally, we apply our criteria of contractivity to the axis of some isometry, providing a characterization of contracting isometries. We first need the following definition:

\begin{definition}
Let $X$ be a CAT(0) cube complex and $g \in \mathrm{Isom}(X)$ an isometry. We say that $g$ \textit{skewers} a pair of hyperplanes $(J_1,J_2)$ if $g^nD_1 \subsetneq D_2 \subsetneq D_1$ for some $n \in \mathbb{Z} \backslash \{ 0\}$ and for some half-spaces $D_1,D_2$ delimited by $J_1,J_2$ respectively.
\end{definition}

\noindent
Our main criterion is:

\begin{thm}\label{equivalence - contracting isometry}
Let $X$ be a CAT(0) cube complex and $g \in \mathrm{Isom}(X)$ an isometry with a combinatorial axis $\gamma$. The following statements are equivalent:
\begin{itemize}
	\item[(i)] $g$ is a contracting isometry;
	\item[(ii)] there exists $C \geq 1$ such that any join of hyperplanes $(\mathcal{H}, \mathcal{V})$ with $\mathcal{H} \subset \mathcal{H}(\gamma)$ is $C$-thin;
	\item[(iii)] there exists $C \geq 1$ such that:
		\begin{itemize}
			\item $\mathcal{H}(\gamma)$ does not contain $C$ pairwise transverse hyperplanes,
			\item any grid of hyperplanes $(\mathcal{H}, \mathcal{V})$ with $\mathcal{H} \subset \mathcal{H}(\gamma)$ is $C$-thin;
		\end{itemize}
	\item[(iv)] $g$ skewers a pair of well-separated hyperplanes.
\end{itemize}
\end{thm}

\noindent
\textbf{Proof.} The equivalence $(i) \Leftrightarrow (ii)$ is a direct consequence of Corollary \ref{géodésiques contractantes}. Then, because a grid of hyperplanes or a collection of pairwise transverse hyperplanes gives rise to a join of hyperplanes, $(ii)$ clearly implies $(iii)$. 

\medskip \noindent
Now, we want to prove $(iii) \Rightarrow (ii)$. Let $C$ denote the constant given by $(iii)$ and let $(\mathcal{H}, \mathcal{V})$ be a join of hyperplanes satisfying $\mathcal{H} \subset \mathcal{H}(\gamma)$. If $\# \mathcal{H}, \# \mathcal{V} \geq \mathrm{Ram}(C)$, then $\mathcal{H}$ and $\mathcal{V}$ each contain a subfamily of at least $C$ pairwise disjoint hyperplanes, say $\mathcal{H}'$ and $\mathcal{V}'$ respectively. But now $(\mathcal{H}',\mathcal{V}')$ defines a grid of hyperplanes satisfying $\# \mathcal{H}', \# \mathcal{V}' \geq C$, contradicting the definition of $C$. Thus, the join $(\mathcal{H},\mathcal{V})$ is necessarily $\mathrm{Ram}(C)$-thin.

\medskip \noindent
Now, we want to prove $(iv) \Rightarrow (i)$. So suppose that there exist two half-spaces $D_1,D_2$ respectively delimited by two well separated hyperplanes $J_1,J_2$ such that $g^nD_1 \subsetneq D_2 \subsetneq D_1$ for some $n \in \mathbb{Z}$. Notice that
\begin{center}
$\cdots \subsetneq g^{2n}D_1 \subsetneq g^nD_2 \subsetneq g^n D_1 \subsetneq D_2 \subsetneq D_1 \subsetneq g^{-n}D_2 \subsetneq g^{-n} D_1 \subsetneq g^{-2n} D_2 \subsetneq \cdots$.
\end{center}
We claim that $J_1$ intersects $\gamma$. Suppose by contradiction that it is not the case. Then, because $\gamma$ is $\langle g \rangle$-invariant, $g^k J_1$ does not intersect $\gamma$ for every $k \in \mathbb{Z}$. As a consequence, there exists some $m \in \mathbb{Z}$ such that $\gamma \subset g^{mn}D_1 \backslash g^{(m+1)n}D_1$. A fortiori, $g^{m} \gamma \subset g^{(m+1)n}D_1 \backslash g^{(m+2)n}D_1$, and because $\gamma$ is $\langle g \rangle$-invariant, we deduce that
$$\gamma \subset \left( g^{mn}D_1 \backslash g^{(m+1)n}D_1 \right) \cap \left( g^{(m+1)n}D_1 \backslash g^{(m+2)n}D_1 \right) = \emptyset,$$
a contradiction. So $J_1$ intersects $\gamma$, and a fortiori, $g^{kn}J_1$ intersects $\gamma$ for every $k \in \mathbb{Z}$. For every $k \in \mathbb{Z}$, the intersection between $g^{kn}N(J_1)$ and $\gamma$ contains exacly two vertices; fix an orientation along $\gamma$ and denote by $y_k$ the first vertex along $\gamma$ which belongs to $g^{kn}N(J_1)$. Clearly,
\begin{center}
$d(y_k,y_{k+1}) = d(g^{kn}y_0,g^{(k+1)n}y_0)=d(y_0,g^ny_0)= \| g \|^n$
\end{center}
since $y_0$ belongs to the combinatorial axis of $g$, where $\|g \|$ denotes the combinatorial translation length of $g$ (see \cite{arXiv:0705.3386}). Furthermore, $J_1$ and $g^nJ_1$ are well separated: indeed, because $J_2$ separates $J_1$ and $g^n J_1$, we deduce that any collection of hyperplanes transverse to both $J_1$ and $g^n J_1$ must also be transverse to both $J_1$ and $J_2$, so that we conclude that $J_1$ and $g^nJ_1$ are well separated since $J_1$ and $J_2$ are themselves well separated. Therefore, $\{ g^{kn}J_1 \mid k \in \mathbb{N} \}$ defines a family of pairwise well separated hyperplanes intersecting $\gamma$ at uniformly separated points. Theorem \ref{géodésiques contractantes 2} implies that $\gamma$ is contracting.

\medskip \noindent
Conversely, suppose that $g$ is contracting. According to Theorem \ref{géodésiques contractantes 2}, there exist three pairwise well separated hyperplanes $J_1,J_2,J_3$ intersecting $\gamma$. Say that they respectively delimit three half-spaces $D_1,D_2,D_3$ satisfying $D_3 \subsetneq D_2 \subsetneq D_1$, where $D_3$ contains $\gamma(+ \infty)$. We claim that $g$ skewers the pair $(J_1,J_2)$. Fix two vertices $x \in N(J_1) \cap \gamma$ and $y \in D_3 \cap \gamma$. Notice that, if $I$ denotes the set of $n \in \mathbb{N}$ such that $g^n J_1$ and $J_1$ are transverse and $g^nJ_1$ does not separate $x$ and $y$, then $\{ g^n J_1 \mid n \in I \}$ defines a set of hyperplanes (without facing triples, since they all intersect $\gamma$) transverse to both $J_2$ and $J_3$. Because $J_2$ and $J_3$ are well separated, necessarily $I$ has to be finite. A fortiori, there must exist at most $\#I + d(x,y)$ integers $n \in \mathbb{N}$ such that $g^n J_1$ and $J_1$ are transverse. As a consequence, there exists an integer $n \in \mathbb{N}$ such that $g^n J_1 \subset D_2$. Finally, we have proved $(i) \Rightarrow (iv)$. $\square$

\section{Combinatorial boundary}

\subsection{Generalities}

\noindent
In this section, we will use the following notation: given two subsets of hyperplanes $\mathcal{H}_1, \mathcal{H}_2$, we define the \emph{almost inclusion} $\mathcal{H}_1 \underset{a}{\subset} \mathcal{H}_2$ if all but finitely many hyperplanes of $\mathcal{H}_1$ belong to $\mathcal{H}_2$. In particular, $\mathcal{H}_1$ and $\mathcal{H}_2$ are \emph{commensurable} provided that $\mathcal{H}_1 \underset{a}{\subset} \mathcal{H}_2$ and $\mathcal{H}_2 \underset{a}{\subset} \mathcal{H}_1$. Notice that commensurability defines an equivalence relation on the set of all the collections of hyperplanes, and the almost-inclusion induces a partial order on the associated quotient set. This allows the following definition.

\begin{definition}
Let $X$ be a CAT(0) cube complex. Its \emph{combinatorial boundary} is the poset $(\partial^c X,\prec)$, where $\partial^c X$ is the set of the combinatorial rays modulo the relation: $r_1 \sim r_2$ if $\mathcal{H}(r_1)$ and $\mathcal{H}(r_2)$ are commensurable; and where the partial order $\prec$ is defined by: $r_1 \prec r_2$ whenever $\mathcal{H}(r_1) \underset{a}{\subset} \mathcal{H}(r_2)$.
\end{definition}

\noindent
Notice that this construction is equivariant, ie., if a group $G$ acts by isometries on $X$, then $G$ acts on $\partial^c X$ by $\prec$-isomorphisms.

\medskip \noindent
The following two lemmas essentially state that we will be able to choose a ``nice" combinatorial ray representing a given point in the combinatorial boundary. They will be useful in the next sections.

\begin{lemma}\label{basepoint}
Let $X$ be a CAT(0) cube complex, $x_0 \in X$ a base vertex and $\xi \in \partial^cX$. There exists a combinatorial ray $r$ with $r(0)=x_0$ and $r= \xi$ in $\partial^cX$.
\end{lemma}

\noindent
\textbf{Proof.} Let $r_n$ be a combinatorial path which is the concatenation of a combinatorial geodesic $[x_0,\xi(n)]$ between $x_0$ and $\xi(n)$ with the subray $\xi_n$ of $\xi$ starting at $\xi(n)$. If $r_n$ is not geodesic, then there exists a hyperplane $J$ intersecting both $[x_0,\xi(n)]$ and $\xi_n$. Thus, $J$ separates $x_0$ and $\xi(n)$ but cannot separate $\xi(0)$ and $\xi(n)$ since otherwise $J$ would intersect the ray $\xi$ twice. We deduce that $J$ necessarily separates $x_0$ and $\xi(0)$. Therefore, if we choose $n$ large enough so that the hyperplanes separating $x_0$ and $\xi(0)$ do not intersect the subray $\xi_n$, then $r_n$ is a combinatorial ray. By construction, we have $r_n(0)=x_0$, and $\mathcal{H}(r_n)$ and $\mathcal{H}(\xi)$ are clearly commensurable so $r_n=\xi$ in $\partial^c X$. $\square$

\begin{lemma}\label{basepointbis}
Let $r, \rho$ be two combinatorial rays satisfying $r \prec \rho$. There exists a combinatorial ray $p$ equivalent to $r$ satisfying $p(0)=\rho(0)$ and $\mathcal{H}(p) \subset \mathcal{H}(\rho)$. 
\end{lemma}

\noindent
\textbf{Proof.} According to the previous lemma, we may suppose that $r(0)= \rho(0)$. By assumption, $\mathcal{H}(r) \underset{a}{\subset} \mathcal{H}(\rho)$. Let $J$ be the last hyperplane of $\mathcal{H}(r) \backslash \mathcal{H}(\rho)$ intersected by $r$, ie., there exists some $k \geq 1$ such that $J$ is dual to the edge $[r(k),r(k+1)]$ and the hyperplanes intersecting the subray $r_k \subset r$ starting at $r(k+1)$ all belong to $\mathcal{H}(\rho)$. We claim that $r_k \subset \partial N(J)$. Indeed, if there exists some $j \geq k+1$ such that $d(r_k(j),N(J)) \geq 1$, then some hyperplane $H$ would separate $r_k(j)$ and $N(J)$ according to Lemma \ref{hyperplan séparant}; a fortiori, $H$ would intersect $r$ but not $\rho$, contradicting the choice of $J$. Let $r'$ denote the combinatorial ray obtained from $r$ by replacing the subray $r_{k-1}$ by the symmetric of $r_k$ with respect to $J$ in $N(J) \simeq J \times [0,1]$. The set of hyperplanes intersecting $r'$ is precisely $\mathcal{H}(r) \backslash \{ J \}$. Because $|\mathcal{H}(r') \backslash \mathcal{H}(\rho)| < |\mathcal{H}(r) \backslash \mathcal{H}(\rho)|$, iterating the process will produce after finitely many steps a combinatorial ray $p$ satisfying $p(0)= r(0)=\rho(0)$ and $\mathcal{H}(p) \subset \mathcal{H}(\rho)$. $\square$

\begin{definition}
Let $X$ be a CAT(0) cube complex and $Y \subset X$ a subcomplex. We define the \emph{relative combinatorial boundary} $\partial^c Y$ of $Y$ as the subset of $\partial^cX$ corresponding the set of the combinatorial rays included into $Y$.
\end{definition}

\begin{lemma}
Let $r$ be a combinatorial ray. Then $r \in \partial^c Y$ if and only if $\mathcal{H}(r) \underset{a}{\subset} \mathcal{H}(Y)$.
\end{lemma}

\noindent
\textbf{Proof.} If $r \in \partial^c Y$, then by definition there exists a combinatorial ray $\rho \subset Y$ equivalent to $r$. Then
\begin{center}
$\mathcal{H}(r) \underset{a}{\subset} \mathcal{H}(\rho) \subset \mathcal{H}(Y)$.
\end{center}
Conversely, suppose that $\mathcal{H}(r) \underset{a}{\subset} \mathcal{H}(Y)$. According to Lemma \ref{basepoint}, we may suppose that $r(0) \in Y$. Using the same argument as in the previous proof, we find an equivalent combinatorial ray $\rho$ with $\mathcal{H}(\rho) \subset \mathcal{H}(Y)$ and $\rho(0)=r(0) \in Y$. Because $Y$ is combinatorially convex, it follows that $\rho \subset Y$, hence $r \in \partial^c Y$. $\square$

\subsection{Quasiconvex geodesics II}

\noindent
The goal of this section is to prove the following criterion, determining when a combinatorial axis is quasiconvex just from its endpoints in the combinatorial boundary.

\begin{prop}\label{prop-qc2}
Let $X$ be a locally finite CAT(0) cube complex and $g\in \mathrm{Isom}(X)$ an isometry with a combinatorial axis $\gamma$. A subray $r \subset \gamma$ is quasiconvex if and only if $r$ is minimal in $\partial^cX$.
\end{prop}

\noindent
\textbf{Proof.} Suppose that $r$ is not quasiconvex. According to Proposition \ref{prop-qc}, for every $n \geq 1$, there exist a join of hyperplanes $(\mathcal{H}_n , \mathcal{V}_n)$ with $\mathcal{H}_n, \mathcal{V}_n \subset \mathcal{H}(r)$ and $\# \mathcal{H}_n, \# \mathcal{V}_n= n$. The hyperplanes of $\mathcal{H}_n \cup \mathcal{V}_n$ are naturally ordered depending on the order of their intersections along $r$. Without loss of generality, we may suppose that the hyperplane of $\mathcal{H}_n \cup \mathcal{V}_n$ which is closest to $r(0)$ belongs to $\mathcal{H}_n$; because $\langle g \rangle$ acts on $\mathcal{H}(\gamma)$ with finitely many orbits, up to translating by powers of $g$ and taking a subsequence, we may suppose that this hyperplane $J \in \mathcal{H}_n$ does not depend on $n$. Finally, let $V_n$ denote the hyperplane of $\mathcal{V}_n$ which is farthest from $r(0)$. 

\medskip \noindent
Thus, if $J^-$ (resp. $V_n^+)$ denotes the halfspace delimited by $J$ (resp. $V_n$) which does not contain $r(+ \infty)$ (resp. which contains $r(+ \infty)$), then $\mathcal{C}_n= (J^-,V_n^+,r)$ defines a cycle of three subcomplexes. Let $D_n \to X$ be a disc diagram of minimal complexity bounded by $\mathcal{C}_n$. Because no hyperplane can intersect $J$, $V_n$ or $r$ twice, it follows from Proposition \ref{disc embedding} that $D_n \to X$ is an isometric embedding, so we will identify $D_n$ with its image in $X$. Let us write the boundary $\partial D_n$ as a concatenation of three combinatorial geodesics
\begin{center}
$\partial D_n = \mu_n \cup \nu_n \cup \rho_n$,
\end{center}
where $\mu_n \subset J^-$, $\nu_n \subset V_n^+$ and $\rho_n \subset r$. Because $X$ is locally finite, up to taking a subsequence, we may suppose that the sequence of combinatorial geodesics $(\mu_n)$ converges to a combinatorial ray $\mu_{\infty} \subset J^-$. Notice that, if $H \in \mathcal{H}(\mu_{\infty})$, then $H \in \mathcal{H}(\mu_k)$ for some $k \geq 1$, it follows from Theorem \ref{disc diagram} that $H$ does not intersect $\nu_k$ in $D_k$, hence $H \in \mathcal{H}(\rho_k)$. Therefore, $\mathcal{H}(\mu_{\infty}) \subset \mathcal{H}(r)$. 

\medskip \noindent
According to Lemma \ref{infinite chain} below, there exists an infinite collection of pairwise disjoint hyperplanes $J_1,J_2, \ldots \in \mathcal{H}(r)$. Because $\langle g \rangle$ acts on $\mathcal{H}(\gamma)$ with finitely many orbits, necessarily there exist some $r,s \geq 1$ and some $m \geq 1$ such that $g^mJ_r=J_s$. Thus, $\{J_r,g^mJ_r,g^{2m}J_r,\ldots \}$ defines a collection of infinitely many pairwise disjoint hyperplanes in $\mathcal{H}(r)$ stable by the semigroup generated by $g^m$. For convenience, let $\{ H_0,H_1, \ldots \}$ denote this collection.

\medskip \noindent
If $J$ is disjoint from some $H_k$, then $H_k,H_{k+1}, \ldots \notin \mathcal{H}(\mu_{\infty})$ since $\mu_{\infty} \subset J^-$. On the other hand, $H_k,H_{k+1},\ldots \in \mathcal{H}(r)$, so we deduce that $\mu_{\infty} \prec r$ with $\mu_{\infty} \neq r$ in $\partial^c X$: it precisely means that $r$ is not minimal in $\partial^c X$.

\medskip \noindent
From now on, suppose that $J$ is transverse to all the $H_k$'s. As a consequence, $g^{km}J$ is transverse to $H_n$ for every $n \geq k$. For every $k \geq 1$, let $H_k^+$ denote the halfspace delimited by $H_k$ which contains $r(+ \infty)$ and let $p_k : X \to H_k^+$ be the combinatorial projection onto $H_k^+$. We define a sequence of vertices $(x_n)$ by:
\begin{center}
$x_0=r(0)$ and $x_{n+1}=p_{n+1}(x_n)$ for every $n \geq 0$.
\end{center}
Because of Lemma \ref{inclusion}, we have $x_n=p_n(x_0)$ for every $n \geq 1$. Therefore, it follows from Proposition \ref{projection} that, for every $n \geq 1$, there exists a combinatorial geodesic between $x_0$ and $x_n$ passing through $x_k$ for $k \leq n$. Thus, there exists a combinatorial ray $\rho$ passing through all the $x_k$'s. 

\medskip \noindent
Let $H \in \mathcal{H}(\rho)$. Then $H$ separates $x_k$ and $x_{k+1}$ for some $k \geq 0$, and it follows from Lemma \ref{hyperplan séparant} that $H$ separates $x_k$ and $H_k$. As a first consequence, we deduce that, for every $k \geq 1$, $g^{km}J$ cannot belong to $\mathcal{H}(\rho)$ since $g^{km}J$ is transverse to $H_n$ for every $n \geq k$, hence $r \neq \rho$ in $\partial^cX$. On the other hand, $H$ separates $x_0=r(0)$ and $H_k$, and $r(n) \in H_k^+$ for $n$ large enough, so $H \in \mathcal{H}(r)$. We have proved that $\mathcal{H}(\rho) \subset \mathcal{H}(r)$, ie., $\rho \prec r$. Thus, $r$ is not minimal in $\partial^c X$. 

\medskip \noindent
Conversely, suppose that $r$ is not minimal in $\partial^cX$. Thus, there exists a combinatorial ray $\rho$ with $\rho(0)=r(0)$ such that $\mathcal{H}(\rho) \subset \mathcal{H}(r)$ and $|\mathcal{H}(r) \backslash \mathcal{H}(\rho)| =+ \infty$. Let $J_1,J_2, \ldots \in \mathcal{H}(r) \backslash \mathcal{H}(\rho)$ be an infinite collection. Without loss of generality, suppose that, for every $i>j \geq 1$, the edge $J_j \cap r$ is closer to $r(0)$ than the edge $J_i \cap r$. Let $N \geq 1$, and let $d(N)$ denote the distance $d(r(0),\omega)$ where $\omega$ is the endpoint of the edge $J_N \cap r$ which is farthest from $r(0)$. Choose any collection of hyperplanes $\mathcal{V} \subset \mathcal{H}(\rho)$ with $\# \mathcal{V} \geq d(N)+N$. Then $\mathcal{V}$ contains a subcollection $\{ V_1, \ldots, V_N \}$ such that the edge $V_i \cap r$ is farther from $r(0)$ than the edge $J_N \cap r$. We know that each $V_j$ separates $\{r(0),\omega\}$ and $\{ r(k),\rho(k) \}$ for $k$ large enough, and each $J_i$ separates $\{ r(k), \omega \}$ and $\{ r(0) , \rho(k) \}$ for $k$ large enough (since $J_i$ and $\rho$ are disjoint); we deduce that $J_i$ and $V_j$ are transverse for any $1 \leq i,j \leq N$. Therefore, $(\{ J_1, \ldots, J_N\}, \{ V_1, \ldots, V_N \})$ defines a join of hyperplanes in $\mathcal{H}(r)$. In fact, we have proved

\begin{fact}\label{fait-rayon}
Let $r \prec \rho$ be two combinatorial rays with the same origin. If $\mathcal{J} \subset \mathcal{H}(r) \backslash \mathcal{H}(\rho)$ is an infinite collection, then for every $N \geq 1$ there exists a join of hyperplanes $(\mathcal{H}, \mathcal{V})$ with $\mathcal{H} \subset \mathcal{J}$, $\mathcal{V} \subset \mathcal{H}(\rho)$ and $\# \mathcal{H}, \# \mathcal{V} \geq n$. 
\end{fact}

\noindent
Since $N$ can be chosen arbitrarily large, it follows from Proposition \ref{prop-qc} that $r$ is not quasiconvex. $\square$

\begin{lemma}\label{infinite chain}
Let $X$ be a complete CAT(0) cube complex and $r$ a combinatorial ray. Then $\mathcal{H}(r)$ contains an infinite sequence of pairwise disjoint hyperplanes.
\end{lemma}

\noindent
\textbf{Proof.} We begin our proof with a completely general argument, without assuming that $X$ is complete. First, we decompose $\mathcal{H}(r)$ as the disjoint union $\mathcal{H}_0 \sqcup \mathcal{H}_1 \sqcup \cdots$ where we define $\mathcal{H}_i = \{ J \in \mathcal{H}(r) \mid d_{\infty}(r(0),J)=i \}$ for all $i \geq 0$. A priori, some $\mathcal{H}_i$ might be empty or infinite; moreover, $\mathcal{H}_i$ could be infinite and $\mathcal{H}_{i+1}$ non-empty. We want to understand the behaviour of this sequence of collections of hyperplanes. Our first two claims state that each $\mathcal{H}_i$ is a collection of pairwise transverse hyperplanes and that $\mathcal{H}_i$ non-empty implies $\mathcal{H}_j$ non-empty for every $j<i$. 

\begin{claim}\label{claim-transverse}
For all $i \geq 0$ and $J_1,J_2 \in \mathcal{H}_i$, $J_1$ and $J_2$ are transverse.
\end{claim}

\noindent
Suppose by contradiction that $J_1$ and $J_2$ are disjoint, and, for convenience, say that $J_1$ separates $r(0)$ and $J_2$. Because $d_{\infty}(r(0),J_1)=i$, there exists $i$ pairwise disjoint hyperplanes $V_1, \ldots, V_i$ separating $r(0)$ and $J_1$. Then $V_1, \ldots, V_i,J_1$ are $i+1$ pairwise transverse hyperplanes separating $r(0)$ and $J_2$, hence $d_{\infty}(r(0),J_2) \geq i+1$, a contradiction.

\begin{claim}\label{claim-chain}
Let $J \in \mathcal{H}_i$, ie., there exists a collection of $i$ pairwise disjoint hyperplanes $V_0, \ldots, V_{i-1}$ separating $r(0)$ and $J$. Then $V_j \in \mathcal{H}_j$ for every $0 \leq j \leq i-1$.
\end{claim}

\noindent
First, because $V_0, \ldots, V_{j-1}$ separate $r(0)$ and $V_j$, necessarily we have $d_{\infty}(r(0),V_j) \geq j$. Let $H_1 \ldots, H_k$ be $k$ pairwise disjoint hyperplanes separating $r(0)$ and $V_j$. Then $H_1, \ldots, H_k,V_j, \ldots, V_{i-1}$ define $k+i-j$ pairwise disjoint hyperplanes separating $r(0)$ and $J$, hence
\begin{center}
$k+i-j \leq d_{\infty}(r(0),J) = i$, 
\end{center}
ie., $k \leq j$. We deduce that $d_{\infty}(r(0),V_j) \leq j$. Finally, we conclude that $d_{\infty}(r(0),V_j) = j$, that is $V_j \in \mathcal{H}_j$.

\medskip \noindent
Our third and last claim states that, if the first collections $\mathcal{H}_0, \ldots, \mathcal{H}_p$ are finite, then there exists a sequence of cubes $Q_0, \ldots, Q_p$ such that the intersection between two successive cubes is a single vertex and the hyperplanes dual to the cube $Q_i$ are precisely the hyperplanes of $\mathcal{H}_i$. For this purpose, we define by induction the sequence of vertices $x_0, x_1, \ldots$ by:
\begin{itemize}
	\item $x_0=r(0)$;
	\item if $x_j$ is defined and $\mathcal{H}_j$ is finite, $x_{j+1}$ is the projection of $x_j$ onto $C_j:= \bigcap\limits_{J \in \mathcal{H}_j} J^+$, where $J^+$ denotes the halfspace delimited by $J$ which does not contain $r(0)$.
\end{itemize}
Notice that the intersection $C_j$ is non-empty precisely because of Claim \ref{claim-transverse} and because we assumed that $\mathcal{H}_j$ is finite.

\begin{claim}
For every $i \geq 0$, if $\mathcal{H}_{i-1}$ is finite then any hyperplane of $\mathcal{H}_i$ is adjacent to $x_{i}$; and, if $\mathcal{H}_i$ is finite, there exists a cube $Q_i$ containing $x_i$ and $x_{i+1}$ as opposite vertices such that $\mathcal{H}_i$ is the set of hyperplanes intersecting $Q_i$.
\end{claim}

\noindent
We argue by induction. 

\medskip \noindent
First, we notice that any hyperplane $J \in \mathcal{H}_0$ is adjacent to $x_0$. Suppose by contradiction that this not the case, and let $p$ denote the combinatorial projection of $x_0$ onto $J^+$. By assumption, $d(x_0,p) \geq 2$ so there exists a hyperplane $H$ separating $x_0$ and $p$ which is different from $J$. According to Lemma \ref{hyperplan séparant}, $H$ separates $x_0$ and $J$. This contradicts $J \in \mathcal{H}_0$. Therefore, any hyperplane $J \in \mathcal{H}_0$ is dual to an edge $e(J)$ with $x_0$ as an endpoint. If $\mathcal{H}_0$ is infinite, there is nothing to prove. Otherwise, $\{ e(J) \mid J \in \mathcal{H}_0 \}$ is a finite set of edges with a common endpoint and whose associated hyperplanes are pairwise transverse. Thus, because a CAT(0) cube complex does not contain inter-osculating hyperplanes, these edges span a cube $Q_0$. By construction, the set of hyperplanes intersecting $Q_0$ is $\mathcal{H}_0$. Furthermore, because the vertex of $Q_0$ opposite to $x_0$ belongs to $C_0$ and is at distance $\# \mathcal{H}_0$ from $x_0$, we deduce that it is precisely $x_1$. 

\medskip \noindent
Now suppose that $x_0, \ldots, x_i$ are well-defined that there exist some cubes $Q_0, \ldots, Q_{i-1}$ satisfying our claim. We first want to prove that any hyperplane $J \in \mathcal{H}_{i}$ is adjacent to $x_i$. Suppose by contradiction that this is not the case, ie., there exists some $J \in \mathcal{H}_i$ which is not adjacent to $x_i$. As a consequence, there must exist a hyperplane $H$ separating $x_i$ and $J$. 

\medskip \noindent
\underline{Case 1:} $H$ does not separate $x_0$ and $x_{i}$. So $H$ separates $x_0$ and $J$, and $H$ does not belong to $\mathcal{H}_k$ for $k \leq i-1$. Noticing that $r(k) \in J^+$ for sufficiently large $k \geq 0$, we deduce that $H$ necessarily intersects $r$, hence $H \in \mathcal{H}_j$ for some $j \geq i$. Therefore, there exist $i$ pairwise disjoint hyperplanes $V_1, \ldots, V_i$ separating $r(0)$ and $H$. A fortiori, $V_1, \ldots, V_i,H$ define $i+1$ pairwise disjoint hyperplanes separating $r(0)$ and $J$, contradicting $J \in \mathcal{H}_i$. 

\medskip \noindent
\underline{Case 2:} $H$ separates $x_0$ and $x_i$. In particular, $H$ intersects a cube $Q_j$ for some $0 \leq j \leq i-1$, ie., $H \in \mathcal{H}_j$. Let $p_0$ (resp. $p_i$) denote the combinatorial projection of $x_0$ (resp. $x_i$) onto $J^+$. Notice that, because $H$ is disjoint from $J$, $H$ does not separate $p_0$ and $p_i$: these two vertices belong to the same half-space delimited by $H$, say $H^+$. Because $H$ separates $x_i$ and $J$, it separates $x_i$ and $p_i$, so $x_i$ belongs to the second half-space delimited by $H$, say $H^-$. Then, since $H$ separates $x_0$ and $x_i$ by hypothesis, we deduce that $x_0$ belongs to $H^+$; in particular, $H$ does not separate $x_0$ and $p_0$. On the other hand, if $k \geq 0$ is sufficiently large so that $r(k) \in J^+$, $H$ must separate $r(0)$ and $r(k)$ since $H \in \mathcal{H}(r)$. According to Proposition \ref{projection}, there exists a combinatorial geodesic between $x_0$ and $r(k)$ passing through $p_0$. Because $H$ is disjoint from $J^+$, we conclude that $H$ must separate $x_0$ and $p_0$, a contradiction.

\medskip \noindent
We have proved that any hyperplane of $\mathcal{H}_i$ is adjacent to $x_i$. Thus, if $\mathcal{H}_i$ is finite, we can construct our cube $Q_i$ as above. This concludes the proof of our claim.

\medskip \noindent
From now on, we assume that $X$ is a complete CAT(0) cube complex. First, as a consequence of the previous claim, we deduce that each $\mathcal{H}_i$ is finite. Indeed, suppose by contradiction that some $\mathcal{H}_i$ is infinite. Without loss of generality, we may suppose that $\mathcal{H}_j$ is finite for every $j<i$, so that the points $x_0, \ldots, x_i$ and the cubes $Q_0, \ldots, Q_{i-1}$ are well-defined. According to our previous claim, any hyperplane of $\mathcal{H}_i$ is adjacent to $x_i$: thus, the infinite set of edges adjacent to $x_i$ and dual to the hyperplanes of $\mathcal{H}_i$ define an infinite cube, contradicting the completeness of $X$.

\medskip \noindent
Therefore, we have defined an infinite sequence of vertices $x_0, x_1, \ldots$ and an infinite sequence of cubes $Q_0,Q_1 , \ldots$. Thanks to Claim \ref{claim-chain}, for every $i \geq 0$, there exists a sequence of pairwise disjoint hyperplanes $V_0^i, \ldots, V_i^i$ with $V_j^i \in \mathcal{H}_j$. Because each $\mathcal{H}_k$ is finite, up to taking a subsequence, we may suppose that, for every $k \geq 0$, the sequence $(V_k^i)$ is eventually constant to some hyperplane $V_k \in \mathcal{H}_k$. By construction, our sequence $V_0, V_1, \ldots$ defines a collection of pairwise disjoint hyperplanes in $\mathcal{H}(r)$, concluding the proof of our lemma. $\square$

\subsection{Contracting convex subcomplexes II}

\noindent
Once again, according to Proposition \ref{prop-qc}, the previous section reduces the problem of determining when a combinatorial axis is contracting to the problem of determining when a (hyperbolic, cocompact) combinatorially convex subcomplex is contracting. To do so, we need the following definition:

\begin{definition}
A subset $S \subset \partial^c X$ is \emph{full} if any element of $\partial^c X$ which is comparable with an element of $S$ necessarily belongs to $S$.
\end{definition}

\noindent
Our main criterion is:

\begin{thm}\label{contracting subcomplex2}
Let $X$ be a locally finite CAT(0) cube complex and $Y \subset X$ a hyperbolic combinatorially convex $\mathrm{Aut}(X)$-cocompact subcomplex. Then $Y$ is contracting if and only if $\partial^c Y$ is full in $\partial^c X$. 
\end{thm}

\noindent
\textbf{Proof.} Suppose that $Y$ is not contracting. Notice that $Y$ is a cocompact subcomplex, so its dimension, say $d$, must be finite. According to Theorem \ref{contracting subcomplexes}, for every $n \geq d$, there exists a join of hyperplanes $(\mathcal{H}_n,\mathcal{V}_n)$ with $\mathcal{H}_n \subset \mathcal{H}(Y)$, $\mathcal{V}_n \cap \mathcal{H}(Y) = \emptyset$ and $\# \mathcal{H}_n= \# \mathcal{V}_n \geq \mathrm{Ram}(n)$. Next, up to taking subcollections of $\mathcal{H}_n$ and $\mathcal{V}_n$, we may suppose thanks to Ramsey's theorem that $\mathcal{H}_n$ and $\mathcal{V}_n$ are collections of pairwise disjoint hyperplanes of cardinalities exactly $n$, ie., $(\mathcal{H}_n,\mathcal{V}_n)$ is a grid of hyperplanes. For convenience, write $\mathcal{H}_n= (H_1^n, \ldots, H_n^n)$ (resp. $\mathcal{V}_n= (V_1^n, \ldots, V_n^n)$ where $H_i^n$ (resp. $V_i^n$) separates $H_{i-1}^n$ and $H_{i+1}^n$ (resp. $V_{i-1}^n$ and $V_{i+1}^n$) for every $2 \leq i \leq n-1$; we also suppose that $V_1^n$ separates $Y$ and $V_n^n$. 

\medskip \noindent
Let $\mathcal{C}_n$ be the cycle of subcomplexes $(N(H_1^n),N(V_n^n), N(H_n^n),Y)$. According to Corollary \ref{flat rectangle}, a disc diagram of minimal complexity bounded by $\mathcal{C}_n$ defines a flat rectangle $D_n$. Say that a hyperplane intersecting the $H_i^n$'s in $D_n$ is \emph{vertical}, and a hyperplane intersecting the $V_i^n$'s in $D_n$ is \emph{horizontal}. 

\begin{claim}
If the grids of hyperplanes of $Y$ are all $C$-thin and $n>C$, then at most $\mathrm{Ram}(C)$ horizontal hyperplanes intersect $Y$. 
\end{claim}

\noindent
Let $\mathcal{V}$ be a collection of horizontal hyperplanes which intersect $Y$; notice that $\mathcal{V}$ does not contain any facing triple since there exists a combinatorial geodesic in $D_n$ (which is a combinatorial geodesic in $X$) intersecting all the hyperplanes of $\mathcal{V}$. If $\# \mathcal{V} \geq \mathrm{Ram}(s)$ for some $s \geq \dim Y$, then $\mathcal{V}$ contains a subcollection $\mathcal{V}'$ with $s$ pairwise disjoint hyperplanes, so that $(\mathcal{H}_n, \mathcal{V}')$ defines a $(n,s)$-grid of hyperplanes. If $n >C$, this implies $s \leq C$. Therefore, $\# \mathcal{V} \leq \mathrm{Ram}(C)$. This proves the claim.

\medskip \noindent
Now, because $Y$ is a cocompact subcomplex, up to translating the $D_n$'s, we may suppose that a corner of each $D_n$ belongs to a fixed compact fundamental domain $C$; and because $X$ is locally finite, up to taking a subsequence, we may suppose that, for every ball $B$ centered in $C$, $(D_n \cap B)$ is eventually constant, so that $(D_n)$ converges naturally to a subcomplex $D_{\infty}$. Notice that $D_{\infty}$ is isomorphic to the square complex $[0,+ \infty) \times [0,+ \infty)$ with $[0,+ \infty) \times \{ 0 \} \subset Y$; let $\rho$ denote this combinatorial ray. Clearly, if $r \subset D_{\infty}$ is a diagonal combinatorial ray starting from $(0,0)$, then $\mathcal{H}(\rho) \subset \mathcal{H}(r)$, ie. $r$ and $\rho$ are comparable in $\partial^c X$. Furthermore, the hyperplanes intersecting $[0,+ \infty) \times \{ 0 \} \subset D_{\infty}$ are horizontal hyperplanes in some $D_n$, and so infinitely many of them are disjoint from $Y$ according to our previous claim: this implies that $r \notin \partial^c Y$. We conclude that $\partial^c Y$ is not full in $\partial^c X$.

\medskip \noindent
Conversely, suppose that $\partial^c Y$ is not full in $\partial^c X$, ie., there exist two $\prec$-comparable combinatorial rays $r,\rho$ satisfying $\rho \subset Y$ and $r \notin \partial^c Y$.

\medskip \noindent 
Suppose that $r \prec \rho$. According to Lemma \ref{basepointbis}, we may suppose that $r(0)= \rho(0)$ and $\mathcal{H}(r) \subset \mathcal{H}(\rho)$. But $r(0) \in Y$ and $\mathcal{H}(r) \subset \mathcal{H}(Y)$ imply $r \subset Y$, contradicting the assumption $r \notin \partial^cY$.

\medskip \noindent
Suppose that $\rho \prec r$. According to Lemma \ref{basepoint}, we may suppose that $r(0)= \rho(0)$. Because $r \notin \partial^cY$, there exists an infinite collection $\mathcal{J} \subset \mathcal{H}(r) \backslash \mathcal{H}(Y)$; a fortiori, $\mathcal{J} \subset \mathcal{H}(r) \backslash \mathcal{H}(\rho)$ since $\rho \subset Y$. It follows from Fact \ref{fait-rayon} that, for every $N \geq 1$, there exists a join of hypergraphs $(\mathcal{H}, \mathcal{V})$ with $\mathcal{H} \subset \mathcal{J}$, $\mathcal{V} \subset \mathcal{H}(\rho)$ and $\# \mathcal{H}, \# \mathcal{V} \geq N$. Therefore, Theorem \ref{contracting subcomplexes} implies that $Y$ is not contracting. $\square$

\begin{remark}
Notice that neither the hyperbolicity of the subcomplex $Y$ nor its cocompactness was necessary to prove the converse of the previous theorem. Thus, the relative combinatorial boundary of a combinatorially convex subcomplex is always full.
\end{remark}

\subsection{Contracting isometries II}

\noindent
Finally, we want to apply the results we have found in the previous sections to characterize contracting isometries from the combinatorial boundary. We begin by introducing the following definition:

\begin{definition}
A point $\xi \in \partial^cX$ is \emph{isolated} if it is not comparable with any other element of $\partial^cX$.
\end{definition}

\noindent
Our main criterion is the following:

\begin{thm}\label{contracting isometry and boundary}
Let $X$ be a locally finite CAT(0) cube complex and $g\in \mathrm{Isom}(X)$ an isometry with a combinatorial axis $\gamma$. Then $g$ is a contracting isometry if and only if $\gamma(+ \infty)$ is isolated in $\partial^cX$. 
\end{thm}

\noindent
\textbf{Proof.} If $g$ is a contracting isometry then a subray $r \subset \gamma$ containing $\gamma(+ \infty)$ is contracting. Because $N(r)$ is a combinatorially convex subcomplex quasi-isometric to a line, it follows that $\partial^cN(r)= \{ r \}$. On the other hand, since $N(r)$ is contracting, it follows from Theorem \ref{contracting subcomplex2} that $\partial^c N(r)$ is full in $\partial^c X$. This precisely means that $r$ is isolated.

\medskip \noindent
Conversely, suppose that $r$ is isolated in $\partial^c X$. It follows from Proposition \ref{prop-qc2} that $r$ is quasiconvex. Thus, $r$ is contracting if and only if $N(r)$ is contracting. The quasiconvexity of $r$ also implies $\partial^cN(r)= \{ r \}$. Because $r$ is isolated in $\partial^c X$, $\partial^c N(r)$ is full in $\partial^c X$, and so $N(r)$ is contracting according to Theorem \ref{contracting subcomplex2}. $\square$

\medskip \noindent
Thus, containing an isolated point in the combinatorial boundary is a necessary condition in order to admit a contracting isometry. The converse is also true if we strengthen our assumptions on the action of our group. The first result in this direction is the following theorem, where we assume that there exist only finitely many orbits of hyperplanes:

\begin{thm}\label{isolatedvertex2}
Let $G$ be a group acting on a locally finite CAT(0) cube complex $X$ with finitely many orbits of hyperplanes. Then $G \curvearrowright X$ contains a contracting isometry if and only if $\partial^c X$ has an isolated vertex.
\end{thm}

\noindent
This result will be essentially a direct consequence of our main technical lemma:

\begin{lemma}\label{wsh-existence}
Let $X$ be a locally finite CAT(0) cube complex. If $\partial^c X$ has an isolated point then, for every $N \geq 2$, there exists a collection of $N$ pairwise well-separated hyperplanes of $X$ which does not contain any facing triple. 
\end{lemma}

\noindent
\textbf{Proof.} Let $r \in \partial^c X$ be an isolated point. According to Lemma \ref{infinite chain}, there exists an infinite collection $V_0,V_1, \ldots \in \mathcal{H}(r)$ of pairwise disjoint hyperplanes.

\begin{claim}
For every $k \geq 0$, $\mathcal{H}(V_k) \cap \mathcal{H}(r)$ is finite.
\end{claim}

\noindent
Suppose by contradiction that $\mathcal{H}(V_k) \cap \mathcal{H}(r)$ is infinite for some $k \geq 0$. Fix a vertex $x_0 \in r$ adjacent to $V_k$ and $H_1, H_2, \ldots \in \mathcal{H}(V_k) \cap \mathcal{H}(r)$. Now, set $x_i$ the combinatorial projection of $x_0$ onto $N(H_i)$ for every $i \geq 1$; notice that $d(x_0,x_i) \underset{i \to + \infty}{\longrightarrow} + \infty$ since $X$ is locally finite so that only finitely many hyperplanes intersect a given ball centered at $x_0$. Then, once again because $X$ is locally finite, if we choose a combinatorial geodesic $[x_0,x_i]$ between $x_0$ and $x_i$ for every $i \geq 1$, up to taking a subsequence, we may suppose that the sequence $([x_0,x_i])$ converges to some combinatorial ray $\rho$. 

\medskip \noindent
Let $J \in \mathcal{H}(\rho)$. Then $J$ separates $x_0$ and $x_i$ for some $i \geq 1$. From Lemma \ref{hyperplan séparant}, we deduce that $J$ separates $x_0$ and $N(H_i)$. On the other hand, we know that $H_i$ intersects the combinatorial ray $r$ which starts from $x_0$, so necessarily $J \in \mathcal{H}(r)$. We have proved that $\mathcal{H}(\rho) \subset \mathcal{H}(r)$holds. 

\medskip \noindent
Next, notice that $x_i \in N(V_k)$ for every $i \geq 1$. Indeed, if $p_i$ denotes the combinatorial projection of $x_0$ onto $N(V_k) \cap N(H_i)$ and $m$ the median vertex of $\{x_0,x_i,p_i \}$, then $x_0,p_i \in N(V_k)$ implies $m \in N(V_k)$ and $x_i,p_i \in N(H_i)$ implies $m \in N(H_i)$, hence $m \in N(H_i) \cap N(V_k)$. Since $x_i$ minimizes the distance to $x_0$ in $N(H_i)$ and that $m$ belongs to to a geodesic between $x_0$ and $x_i$, we conclude that $x_i=m \in N(V_k)$. A fortiori, $\rho \subset N(V_k)$. As a consequence, $V_{k+1},V_{k+2}, \ldots \notin \mathcal{H}(\rho)$, so $\rho$ defines a point of $\partial^cX$, different from $r$, which is $\prec$-comparable with $r$: this contradicts the fact that $r$ is isolated in $\partial^c X$. Our claim is proved.

\begin{claim}\label{mainclaim}
For every $k \geq 0$, there exists some $i \geq 1$ such that $V_k$ and $V_{k+i}$ are $(2i)$-well-separated.
\end{claim}

\noindent
Suppose by contradiction that $V_k$ and $V_{k+i}$ are not well-separated for every $i \geq 1$, so that there exists a collection of hyperplanes $\mathcal{H}^i \subset \mathcal{H}(V_k) \cap \mathcal{H}(V_{k+i})$ which does not contain any facing triple and satisfying $\# \mathcal{H}^i \geq 2i$. Let $x_0 \in r$ be a vertex adjacent to $V_k$. Because $\mathcal{H}^i$ does not contain any facing triple and that $\mathcal{H}^i \subset \mathcal{H}(V_k)$, there exists a vertex $x_i \in N(V_k)$ such that $x_0$ and $x_i$ are separated by at least $i$ hyperplanes of $\mathcal{H}^i$; let $\mathcal{K}^i$ denote this set of hyperplanes. For every $i \geq 1$, let
\begin{center}
$C_i= \bigcap\limits_{J \in \mathcal{K}^i} \{ \text{halfspace delimited by} \ J \ \text{containing} \ x_i\}$.
\end{center}
Let $\mathcal{C}_i$ be the cycle of subcomplexes $(N(V_k),C_i,N(V_{k+i}),r)$, and let $D_i \to X$ be a disc diagram of minimal complexity bounded by $\mathcal{C}_i$. According to Corollary \ref{embedding}, $D_i \to X$ is an isometric embedding so that we will identify $D_i$ with its image in $X$. Because $\mathcal{H}(V_k) \cap \mathcal{H}(r)$ is finite according to our previous claim, we know that, for sufficently large $i \geq 1$, $\mathcal{K}^i$ will contain a hyperplane disjoint from $r$, so that $C_i \cap r= \emptyset$. Up to taking a subsequence, we will suppose that this is always the case. In particular, $\partial D_i$ can be decomposed as the concatenation of four non trivial combinatorial geodesics
\begin{center}
$\partial D_i = \alpha_i \cup \beta_i \cup \gamma_i \cup r_i$,
\end{center}
where $\alpha_i \subset N(V_k)$, $\beta_i \subset C_i$, $\gamma_i \subset N(V_{k+i})$ and $r_i \subset r$. Notice that the intersection between $\alpha_i$ and $r_i$ is a vertex adjacent to $x_0$, so, if we use the local finiteness of $X$ to define, up to taking a subsequence, a subcomplex $D_{\infty}$ as the limit of $(D_i)$, then $D_{\infty}$ is naturally bounded by two combinatorial rays $\alpha_{\infty}$ and $r_{\infty}$ which are the limits of $(\alpha_i)$ and $(r_i)$ respectively; in particular, $\alpha_{\infty} \subset N(V_k)$ and $r_{\infty}$ is a subray of $r$. 

\medskip \noindent
For each $i \geq 1$, we will say that a hyperplane of $D_i$ intersecting $\beta_i$ is \emph{vertical}, and a hyperplane of $D_i$ intersecting $\alpha_i$ is \emph{horizontal}; a horizontal hyperplane $J \in \mathcal{H}(D_i)$ is \emph{high} if $d(x_0,J \cap \alpha_i)> \# \mathcal{H}(V_k) \cap \mathcal{H}(r)$. 

\begin{fact}\label{hv-hyperplanes}
Let $i \geq 1$. The vertical hyperplanes of $D_i$ are pairwise disjoint and intersect $r_i$. The horizontal hyperplanes of $D_i$ are pairwise disjoint. The high horizontal hyperplanes of $D_i$ intersect $\gamma_i$. 
\end{fact} 

\noindent
It follows directly from Theorem \ref{disc diagram} that the vertical hyperplanes of $D_i$ are pairwise disjoint and intersect $r_i$; and the horizontal hyperplanes of $D_i$ are pairwise disjoint. Suppose by contradiction that a high horizontal hyperplane $H$ does not intersect $\gamma_i$. According to Theorem \ref{disc diagram}, necessarily $H$ intersects $r$. On the other hand, because $d(x_0,H \cap \alpha_i)> \# \mathcal{H}(V_k) \cap \mathcal{H}(r)$, there must exist a horizontal hyperplane $J$ intersecting $\alpha_i$ between $H$ and $x_0$ which does not intersect $r_i$; once again, it follows from Theorem \ref{disc diagram} that $J$ has to intersect $\gamma_i$. A fortiori, $J$ and $H$ are necessarily transverse, contradicting the fact that two horizontal hyperplanes are disjoint. This proves the fact.

\medskip \noindent
Now, we want to define a particular subdiagram $D_i' \subset D_i$. Let $H_i$ be the nearest high horizontal hyperplane of $D_i$ from $\alpha_i \cap r_i$. According to the previous fact, $H_i$ intersects $\gamma_i$. Let $D_i'$ be the subdiagram delimited by $H_i$ which contains $\beta_i$. Notice that the hyperplanes intersecting $D_i'$ are either vertical or high horizontal. Thus, from the decomposition of $\mathcal{H}(D_i')$ as horizontal and vertical hyperplanes given by Fact \ref{hv-hyperplanes}, we deduce that $D_i'$ is isometric to some square complex $[0,a_i] \times [0,b_i]$. Notice that $b_i \underset{i \to + \infty}{\longrightarrow} + \infty$ since $V_{k+1}, \ldots, V_{k+i-1}$ separate $\alpha_i$ and $\gamma_i$, and similarly $a_i \underset{i \to + \infty}{\longrightarrow} + \infty$ since $a_i= \mathrm{length}(\alpha_i)- | \mathcal{H}(V_k) \cap \mathcal{H}(r)|$ and $\# \mathcal{K}^i \geq i$. Therefore, if $D_{\infty}' \subset D_{\infty}$ denotes the limit of $(D_i')$, $D_{\infty}'$ is naturally isomorphic to the square complex $[0,+ \infty) \times [0,+ \infty)$. Let $\rho_{\infty} \subset D_{\infty}'$ be the combinatorial ray associated to $[0,+ \infty) \times \{ 0 \}$, and $\mu \subset D_{\infty}'$ be a \emph{diagonal} combinatorial ray, ie., a ray passing through the vertices $\{(i,i) \mid i \geq 0 \}$.

\medskip \noindent
Notice that $\rho_{\infty}$ is naturally the limit of $(\rho_i):=(\partial D_i' \cap N(H_i))$. Thus, any hyperplane intersecting $\rho_{\infty}$ is a vertical hyperplane in some $D_i$. It follows from Fact \ref{hv-hyperplanes} that $\mathcal{H}(\rho_{\infty}) \subset \mathcal{H}(r)$. Because $r$ is isolated in $\partial^c X$, we deduce that $r=\rho_{\infty}$ in $\partial^c X$. 

\medskip \noindent
Now, we clearly have $\mathcal{H}(\rho_{\infty}) \subset \mathcal{H}(\mu)$, hence $\mathcal{H}(r) \underset{a}{\subset} \mathcal{H}(\mu)$. On the other hand, $\mu$ intersects infinitely many high horizontal hyperplanes, which are disjoint from $r$ according to Fact \ref{hv-hyperplanes}. Therefore, $\mu$ is a point of $\partial^cX$ which is $\prec$-comparable with $r$, but different from it, contradicting the assumption that $r$ is isolated in $\partial^cX$. This concludes the proof of our claim.

\medskip \noindent
Now, we are able to conclude the proof of our lemma. Applying Claim \ref{mainclaim} $N-1$ times, we find a sequence of hyperplanes $V_{i(1)},\ldots, V_{i(N)}$ such that $V_{i(k)}$ and $V_{i(k+1)}$ are $L(k)$-well-separated for some $L(k) \geq 1$. Thus, $\{V_{i(1)}, \ldots, V_{i(N)} \}$ defines a collection of $N$ pairwise $L$-well-separated hyperplanes, with $L= \max(L(1), \ldots, L(N-1))$, which does not contain any facing triple. $\square$

\medskip \noindent
\textbf{Proof of Theorem \ref{isolatedvertex2}.} If $G$ contains a contracting isometry, then $\partial^cX$ contains an isolated point according to Theorem \ref{contracting isometry and boundary}. Conversely, suppose that $\partial^cX$ contains an isolated point. Let $N$ denote the number of orbits of hyperplanes for the action $G \curvearrowright X$. According to Lemma \ref{wsh-existence}, $X$ contains a family of $3N$ pairwise well-separated hyperplanes which does not contain any facing triple. By our choice of $N$, we deduce that there exist a hyperplane $J$ and $g,h \in G$ such that $J,gJ,hgJ$ are pairwise well-separated. Without loss of generality, suppose that $hgJ$ separates $J$ and $gJ$, and let $J^+$ be the halfspace delimited by $J$ which contains $gJ$ and $hgJ$. If $gJ^+ \subset J^+$ or $hgJ^+ \subset J^+$, we deduce that $g$ or $hg$ skewers a pair of well-separated hyperplanes, and we deduce from Theorem \ref{equivalence - contracting isometry} that $G$ contains a contracting isometry. Otherwise, $hgJ^+ \subset gJ^+$ since $hgJ$ separates $J$ and $gJ$. Therefore, $h$ skewers a pair of well-separated hyperplanes, and we deduce from Theorem \ref{equivalence - contracting isometry} that $G$ contains a contracting isometry. $\square$

\medskip \noindent
Our second result stating that the existence of an isolated point in the combinatorial boundary assures the existence of a contracting isometry is the following:

\begin{thm}\label{isolatedvertex3}
Let $G$ be a countable group acting on a locally finite CAT(0) cube complex $X$. Suppose that the action $G \curvearrowright X$ is minimal (ie., $X$ does not contain a proper $G$-invariant combinatorially convex subcomplex) and $G$ does not fix a point of $\partial^cX$. Then $G$ contains a contracting isometry if and only if $\partial^c X$ contains an isolated point.
\end{thm}

\noindent
Given a hyperplane $J$, delimiting two haflspaces $J^-$ and $J^+$, an \emph{orientation} of $J$ is the choice of an ordered couple $(J^-,J^+)$. The following terminology was introduced in \cite{MR2827012}:

\begin{definition}
Let $G$ be a group acting on a CAT(0) cube complex. A hyperplane $J$ is \emph{$G$-flippable} if, for every orientation $(J^-,J^+)$ of $J$, there exists some $g \in G$ such that $gJ^- \subsetneq J^+$. 
\end{definition}

\noindent
Our main tool to prove Theorem \ref{isolatedvertex3} is:

\begin{lemma}\label{flipping}
Let $G$ be a countable group acting minimally on a CAT(0) cube complex $X$ without fixing any point of $\partial^cX$. Then any hyperplane of $X$ is $G$-flippable.
\end{lemma}

\noindent
\textbf{Proof.} Suppose that there exists a hyperplane $J$ which is not $G$-flippable, ie., $J$ admits an orientation $(J^-,J^+)$ such that $gJ^+ \cap J^+ \neq \emptyset$ for every $g \in G$. Thus, $\{gJ^+ \mid g \in G\}$ defines a family of pairwise intersecting halfspaces. If the intersection $\bigcap\limits_{g \in G} gJ^+$ is non-empty, it defines a proper $G$-invariant combinatorially convex subcomplex, so that the action $G \curvearrowright X$ is not minimal.

\medskip \noindent
From now on, suppose that $\bigcap\limits_{g \in G} gJ^+ = \emptyset$. Fix an enumeration $G= \{ g_1,g_2, \ldots\}$ and let $I_n$ denote the non-empty intersection $\bigcap\limits_{i=1}^n g_iJ^+$. In particular, we have $I_1 \supset I_2 \supset \cdots$ and $\bigcap\limits_{n \geq 1} I_n = \emptyset$. Define a sequence of vertices $(x_n)$ by:
\begin{center}
$x_0 \notin I_1$ and $x_{n+1}=p_{n+1}(x_n)$ for every $n \geq 0$,
\end{center}
where $p_k : X \to I_k$ denotes the combinatorial projection onto $I_k$. According to Proposition \ref{projection}, for every $n \geq 1$, there exists a combinatorial geodesic between $x_0$ and $x_n$ passing through $x_k$ for any $1 \leq k \leq n$, so there exists a combinatorial ray $r$ passing through all the $x_n$'s. Notice that, because we also have $x_n=p_n(x_0)$ for every $n \geq 1$ according to Lemma \ref{inclusion}, $\mathcal{H}(r)$ is precisely the set of the hyperplanes separating $x_0$ from some $I_n$. In particular, if $\rho$ is any combinatorial ray starting from $x_0$ and intersecting all the $I_n$'s, then $\mathcal{H}(r) \subset \mathcal{H}(\rho)$. 

\medskip \noindent
Therefore, $\bigcap\limits_{k \geq 1} \partial^c I_k$ is a non-empty (since it contains $r$) subset of $\partial^c X$ with $r$ as a minimal element, ie., for any $\rho \in \bigcap\limits_{k \geq 1} \partial^c I_k$ necessarily $r \prec \rho$. 

\medskip \noindent
Let us prove that $\bigcap\limits_{k \geq 1} \partial^c I_k$ is $G$-invariant. Fix some $g \in G$, and, for every $k \geq 1$, let $L(k,g)$ be a sufficiently large integer so that $\{gg_1, \ldots, gg_k\} \subset \{ g_1, \ldots, g_{L(k,g)}\}$. Of course, $L(k,g) \geq k$ and $g I_k = \bigcap\limits_{i=1}^k gg_iJ^+ \supset I_{L(k,g)}$. Consequently, 
$$g \bigcap\limits_{k \geq 1} \partial^c I_k = \bigcap\limits_{k \geq 1} \partial^c (gI_k) \supset \bigcap\limits_{k \geq 1} \partial^c I_{L(k,g)} = \bigcap\limits_{k \geq 1} \partial^c I_k,$$
because $L(k,g) \underset{k \to \infty}{\longrightarrow} + \infty$. Thus, we have proved that $g \bigcap\limits_{k \geq 1} \partial^c I_k \supset \bigcap\limits_{k \geq 1} \partial^c I_k$ for every $g \in G$. It follows that our intersection is $G$-invariant. 

\medskip \noindent
Finally, because $G$ acts on $\partial^c X$ by order-automorphisms, we conclude that $G$ fixes $r \in \partial^c X$. $\square$

\begin{cor}\label{double flipping}
Let $G$ be a countable group acting minimally on a CAT(0) cube complex $X$ without fixing any point of $\partial^cX$. Then any pair of disjoint hyperplanes is skewered by an element of $G$.
\end{cor}

\noindent
\textbf{Proof.} Let $(J,H)$ be a pair of disjoint hyperplanes. Fix two orientations $J=(J^-,J^+)$ and $H=(H^-,H^+)$ so that $J^+ \subsetneq H^+$. Now, by applying the previous lemma twice, we find two elements $g,h \in G$ such that $gJ^- \subsetneq J^+$ and $hgH^+ \subsetneq gH^-$. Thus,
\begin{center}
$hgH^+ \subsetneq gH^-= X \backslash gH^+ \subsetneq X \backslash gJ^+=gJ^- \subsetneq J^+$.
\end{center}
We deduce that $hg$ skewers the pair $(J,H)$. $\square$

\begin{remark}
Lemma \ref{flipping} may be thought of as a generalisation of the Flipping Lemma proved in \cite{MR2827012}, where the Tits boundary is replaced with the combinatorial boundary in possibly infinite dimension. Corollary \ref{double flipping} corresponds to the Double Skewering Lemma \cite{MR2827012}. 
\end{remark}

\noindent
\textbf{Proof of Theorem \ref{isolatedvertex3}.} If $G$ contains a contracting isometry, then $\partial^cX$ contains an isolated point according to Theorem \ref{contracting isometry and boundary}. Conversely, suppose that $\partial^cX$ contains an isolated point. It follows from Lemma \ref{wsh-existence} that $X$ contains two well-separated hyperplanes, and it suffices to invoke Corollary \ref{double flipping} to find an element $g \in G$ which skewers this pair of hyperplanes: according to Theorem \ref{equivalence - contracting isometry}, this is a contracting isometry. $\square$

\section{Acylindrical hyperbolicity of diagram groups}

\subsection{Diagram groups and their cube complexes}\label{section:preliminarydiag}

\noindent
We refer to \cite[§3 and §5]{MR1396957} for a detailed introduction to \textit{semigroup diagrams} and \textit{diagram groups}. 

For an alphabet $\Sigma$, let $\Sigma^+$ denote the free semigroup over $\Sigma$. If $\mathcal{P}= \langle \Sigma \mid \mathcal{R} \rangle$ is a semigroup presentation, where $\mathcal{R}$ is a set of pairs of words in $\Sigma^+$, the semigroup associated to $\mathcal{P}$ is the one given by the factor-semigroup $\Sigma^+ / \sim$ where $\sim$ is the smallest equivalent relation on $\Sigma^+$ containing $\mathcal{R}$. We will always assume that if $u=v \in \mathcal{R}$ then $v=u \notin \mathcal{R}$; in particular, $u=u \notin \mathcal{R}$.

A \textit{semigroup diagram over $\mathcal{P}$} is the analogue for semigroups to van Kampen diagrams for group presentations. Formally, it is a finite connected planar graph $\Delta$ whose edges are oriented and labelled by the alphabet $\Sigma$, satisfying the following properties:
\begin{itemize}
	\item $\Delta$ has exactly one vertex-source $\iota$ (which has no incoming edges) and exactly one vertex-sink $\tau$ (which has no outgoing edges);
	\item the boundary of each cell has the form $pq^{-1}$ where $p=q$ or $q=p \in \mathcal{R}$;
	\item every vertex belongs to a positive path connecting $\iota$ and $\tau$;
	\item every positive path in $\Delta$ is simple.
\end{itemize}
In particular, $\Delta$ is bounded by two positive paths: the \textit{top path}, denoted $\text{top}(\Delta)$, and the \textit{bottom path}, denoted $\text{bot}(\Delta)$. By extension, we also define $\text{top}(\Gamma)$ and $\text{bot}(\Gamma)$ for every \textit{subdiagram} $\Gamma$. In the following, the notations $\text{top}(\cdot)$ and $\text{bot}(\cdot)$ will refer to the paths and to their labels. Also, a \textit{$(u,v)$-cell} (resp. a \textit{$(u,v)$-diagram}) will refer to a cell (resp. a semigroup diagram) whose top path is labelled by $u$ and whose bottom path is labelled by $v$.

Two words $w_1,w_2$ in $\Sigma^+$ are \textit{equal modulo $\mathcal{P}$} if their images in the semigroup associated to $\mathcal{P}$ coincide. In particular, there exists a \textit{derivation} from $w_1$ to $w_2$, i.e., a sequence of relations of $\mathcal{R}$ allowing us to transform $w_1$ into $w_2$. To any such derivation is associated a semigroup diagram, or more precisely a $(w_1,w_2)$-diagram, whose construction is clear from the example below. As in the case for groups, the words $w_1,w_2$ are equal modulo $\mathcal{P}$ if and only if there exists a $(w_1,w_2)$-diagram.

\begin{ex}
Let $\mathcal{P}= \langle a,b,c \mid ab=ba, ac=ca, bc=cb \rangle$ be a presentation of the free Abelian semigroup of rank three. In particular, the words $a^2bc$ and $caba$ are equal modulo $\mathcal{P}$, with the following possible derivation:
\begin{center}
$aabc \overset{(a,ab \to ba,c)}{\longrightarrow} abac \overset{(ab,ac \to ca, \emptyset)}{\longrightarrow} abca \overset{(a,bc \to cb,a)}{\longrightarrow} acba \overset{(\emptyset,ac \to ca,ba)}{\longrightarrow} caba$.
\end{center}
Then, the associated $(a^2bc,caba)$-diagram $\Delta$ is:
\begin{center}
\includegraphics[scale=0.6]{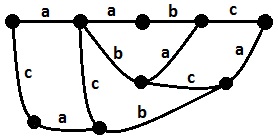}
\end{center}
On such a graph, the edges are supposed oriented from left to right. Here, the diagram $\Delta$ has nine vertices, twelve edges and four cells; notice that the number of cells of a diagram corresponds to the length of the associated derivation. The paths $\text{top}(\Delta)$ and $\text{bot}(\Delta)$ are labelled respectively by $a^2bc$ and $caba$. 
\end{ex}

Since we are only interested in the combinatorics of semigroup diagrams, we will not distinguish isotopic diagrams. For example, the two diagrams below will be considered as equal.
\begin{center}
\includegraphics[scale=0.6]{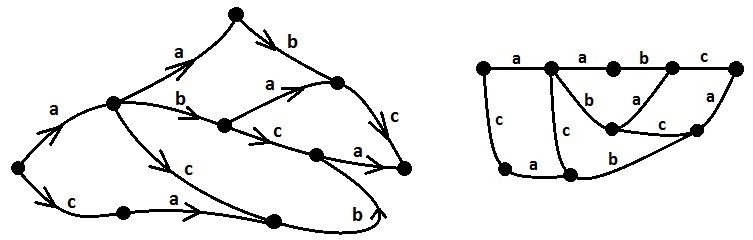}
\end{center}
If $w \in \Sigma^+$, we define the \textit{trivial diagram} $\epsilon(w)$ as the semigroup diagram without cells whose top and bottom paths, labelled by $w$, coincide. Any diagram without cells is trivial. A diagram with exactly one cell is \textit{atomic}.

If $\Delta_1$ is a $(w_1,w_2)$-diagram and $\Delta_2$ a $(w_2,w_3)$-diagram, we define the \textit{concatenation} $\Delta_1 \circ \Delta_2$ as the semigroup diagram obtained by identifying the bottom path of $\Delta_1$ with the top path of $\Delta_2$. In particular, $\Delta_1 \circ \Delta_2$ is a $(w_1,w_3)$-diagram. Thus, $\circ$ defines a partial operation on the set of semigroup diagrams over $\mathcal{P}$. However, restricted to the subset of $(w,w)$-diagrams for some $w \in \Sigma^+$, it defines a semigroup operation; such diagrams are called \textit{spherical with base $w$}. We also define the \textit{sum} $\Delta_1+ \Delta_2$ of two diagrams $\Delta_1,\Delta_2$ as the diagram obtained by identifying the rightmost vertex of $\Delta_1$ with the leftmost vertex of $\Delta_2$.
\begin{center}
\includegraphics[scale=0.6]{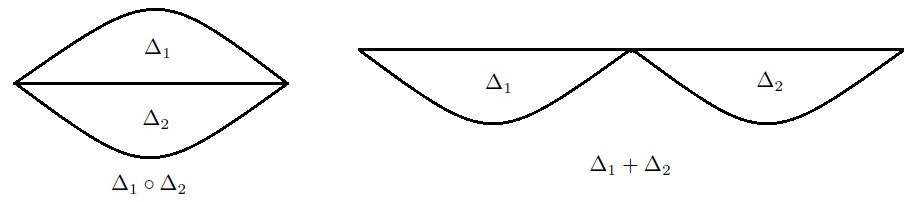}
\end{center}
Notice that any semigroup diagram can be viewed as a concatenation of atomic diagrams. In the following, if $\Delta_1,\Delta_2$ are two diagrams, we will say that $\Delta_1$ is a \textit{prefix} (resp. a \textit{suffix}) of $\Delta_2$ if there exists a diagram $\Delta_3$ satisfying $\Delta_2= \Delta_1 \circ \Delta_3$ (resp. $\Delta_2= \Delta_3 \circ \Delta_1$). Throughout this paper, the fact that $\Delta$ is a prefix of $\Gamma$ will be denoted by $\Delta \leq \Gamma$.

Suppose that a diagram $\Delta$ contains a $(u,v)$-cell and a $(v,u)$-cell such that the top path of the first cell is the bottom path of the second cell. Then, we say that these two cells form a \textit{dipole}. In this case, we can remove these two cells by first removing their common path, and then identifying the bottom path of the first cell with the top path of the second cell; thus, we \textit{reduce the dipole}. A diagram is called \textit{reduced} if it does not contain dipoles. By reducing dipoles, a diagram can be transformed into a reduced diagram, and a result of Kilibarda [Kil94] proves that this reduced form is unique. If $\Delta_1,\Delta_2$ are two diagrams for which $\Delta_1 \circ \Delta_2$ is well defined, let us denote by $\Delta_1 \cdot \Delta_2$ the reduced form of $\Delta_1 \circ \Delta_2$.

Thus, the set of reduced semigroup diagrams, endowed with the product $\cdot$ we defined above, naturally defines a groupoid $\mathcal{G}(\mathcal{P})$, ie., loosely speaking, a ``group'' where the product is only partially defined. The neutral elements correspond to the trivial diagrams, and the inverse of a diagram is constructed in the following way: if $\Delta$ is a $(w_1,w_2)$-diagram, its inverse $\Delta^{-1}$ is the $(w_2,w_1)$-diagram obtained from $\Delta$ by a mirror reflection with respect to $\mathrm{top}(\Delta)$. A natural way to deduce a group from this groupoid is to fix a base word $w \in \Sigma^+$, and to define the \textit{diagram group} $D(\mathcal{P},w)$ as the set of reduced $(w,w)$-diagrams endowed with the product $\cdot$ we defined above.

\paragraph{Farley cube complexes.} The groupoid $\mathcal{G}(\mathcal{P})$ has a natural generating set, namely the set of atomic diagrams, and the group $D(\mathcal{P},w)$ acts by left-multiplication on the connected component of the associated Cayley graph which contains $\epsilon(w)$. Surprisingly, this Cayley graph turns out to be naturally the 1-skeletton of a CAT(0) cube complex. Below, we detail the construction of this cube complex as explained in \cite{MR1978047}.

A semigroup diagram is \textit{thin} whenever it can be written as a sum of atomic diagrams. We define the \textit{Farley complex} $X(\mathcal{P},w)$ as the cube complex whose vertices are the reduced semigroup diagrams $\Delta$ over $\mathcal{P}$ satisfying $\mathrm{top}(\Delta)=w$, and whose $n$-cubes are spanned by the vertices $\{ \Delta \cdot P \mid P \leq \Gamma \}$ for some vertex $\Delta$ and some thin diagram $\Gamma$ with $n$ cells. In particular, two diagrams $\Delta_1$ and $\Delta_2$ are linked by an edge if and only if there exists an atomic diagram $A$ such that $\Delta_1= \Delta_2 \cdot A$.

\begin{thm}\emph{\cite[Theorem 3.13]{MR1978047}}
$X(\mathcal{P},w)$ is a CAT(0) cube complex. Moreover it is complete, i.e., there are no infinite increasing sequences of cubes in $X(\mathcal{P},w)$.
\end{thm}

\noindent
There is a natural action of $D(\mathcal{P},w)$ on $X(\mathcal{P},w)$, namely $(g, \Delta) \mapsto g \cdot \Delta$. Then

\begin{prop}\emph{\cite[Theorem 3.13]{MR1978047}} 
The action $D(\mathcal{P},w) \curvearrowright X(\mathcal{P},w)$ is free. Moreover, it is cocompact if and only if the class $[w]_{\mathcal{P}}$ of words equal to $w$ modulo $\mathcal{P}$ is finite.
\end{prop}

\noindent
In this paper, we will always suppose that the cube complex $X(\mathcal{P},w)$ is locally finite; in particular, this implies that the action $D(\mathcal{P},w) \curvearrowright X(\mathcal{P},w)$ is properly discontinuous. For instance, this is case if $\mathcal{P}$ is a finite presentation, and because we may always suppose that $\mathcal{P}$ is a finite presentation if our diagram group is finitely generated, this assumption is not really restrictive. However, notice that $X(\mathcal{P},w)$ may be locally finite even if $\mathcal{P}$ is an infinite presentation. For example, the derived subgroup $[F,F]$ of Thompson's group $F$ is isomorphic to the diagram group $D(\mathcal{P},a_0b_0)$, where
\begin{center}
$\mathcal{P} = \langle x,a_i,b_i, i \geq 0 \mid x=x^2, a_i=a_{i+1}x, b_i=xb_{i+1}, i \geq 0 \rangle$.
\end{center}
See \cite[Theorem 26]{MR1725439} for more details. In this situation, $X(\mathcal{P},a_0b_0)$ is nevertheless locally finite. 

\medskip \noindent
Because $X(\mathcal{P},w)$ is essentially a Cayley graph, the combinatorial geodesics behave as expected:

\begin{lemma}\label{geodesic}\emph{\cite[Lemma 2.3]{arXiv:1505.02053}}
Let $\Delta$ be a reduced diagram. If $\Delta = A_1 \circ \cdots \circ A_n$ where each $A_i$ is an atomic diagram, then the path
\begin{center}
$\epsilon(w), \ A_1, \ A_1 \circ A_2, \ldots, \ A_1 \circ \cdots \circ A_n= \Delta$
\end{center}
defines a combinatorial geodesic from $\epsilon(w)$ to $\Delta$. Furthermore, every combinatorial geodesic from $\epsilon(w)$ to $\Delta$ has this form.
\end{lemma}

\begin{cor}\emph{\cite[Corollary 2.4]{arXiv:1505.02053}}
Let $A,B \in X(\mathcal{P},w)$ be two reduced diagrams. Then we have 
\begin{center}
$d(A,B)= \# (A^{-1} \cdot B)$,
\end{center}
where $\# ( \cdot)$ denotes the number of cells of a diagram. In particular, $d(\epsilon(w),A)= \# (A)$.
\end{cor}

\noindent
In \cite{arXiv:1505.02053}, we describe the hyperplanes of $X(\mathcal{P},w)$. We recal this description below. For any hyperplane $J$ of $X( \mathcal{P},w)$, we introduce the following notation:
\begin{itemize}
	\item[$\bullet$] $J^+$ is the halfspace associated to $J$ not containing $\epsilon(w)$,
	\item[$\bullet$] $J^-$ is the halfspace associated to $J$ containing $\epsilon(w)$,
	\item[$\bullet$] $\partial_{\pm} J$ is the intersection $\partial N(J) \cap J^{\pm}$.
\end{itemize}

\begin{definition}
A diagram $\Delta$ is \textit{minimal} if its maximal thin suffix $F$ has exactly one cell. (The existence and uniqueness of the maximal thin suffix is given by \cite[Lemma 2.3]{MR1978047}.)
\end{definition}

\noindent
In the following, $\overline{\Delta}$ will denote the diagram $\Delta \cdot F^{-1}$, obtained from $\Delta$ by removing the suffix $F$. The following result uses minimal diagrams to describe hyperplanes in Farley complexes.

\begin{prop}\label{prop:hypFarley}\emph{\cite[Proposition 2.2]{arXiv:1505.02053}} 
Let $J$ be a hyperplane of $X(\mathcal{P},w)$. Then there exists a unique minimal diagram $\Delta$ such that $J^+= \{ D \ \text{diagram} \mid \Delta \leq D \}$. Conversely, if $\Delta$ is a minimal diagram and $J$ the hyperplane dual to the edge $[\Delta, \overline{\Delta}]$, then $J^+ = \{ D \ \text{diagram} \mid \Delta \leq D \}$.
\end{prop}

\noindent
Thanks to the geometry of CAT(0) cube complexes, we will use this description of the hyperplanes of $X(\mathcal{P},w)$ in order to define the \emph{supremum}, when it exists, of a collection of minimal diagrams. If $\Delta$ is a diagram, let $\mathcal{M}(\Delta)$ denote the set of its minimal prefixes. A set of minimal diagrams $\mathcal{D}$ is \emph{consistent} if, for every $D_1,D_2 \in \mathcal{D}$, there exists some $D_3 \in \mathcal{D}$ satisfying $D_1,D_2 \leq D_3$. 

\begin{prop}\label{finite-consistent}
Let $\mathcal{D}$ be a finite consistent collection of reduced minimal $(w,\ast)$-diagrams. Then there exists a unique $\leq$-minimal reduced $(w,\ast)$-diagram admitting any element of $\mathcal{D}$ as a prefix; let $\mathrm{sup}( \mathcal{D})$ denote this diagram. Furthermore, $\mathcal{M}( \mathrm{sup}(\mathcal{D})) = \bigcup\limits_{D \in \mathcal{D}} \mathcal{M}(D)$.
\end{prop}

\noindent
\textbf{Proof.} Let $\mathcal{D}=\{ D_1, \ldots, D_r \}$. For every $1 \leq i \leq r$, let $J_i$ denote the hyperplane associated to $D_i$. Because $\mathcal{D}$ is consistent, for every $1 \leq i,j \leq r$ the intersection $J_i^+ \cap J_j^+$ is non-empty, hence $C=\bigcap\limits_{i=1}^r J_i^+ \neq \emptyset$. Let $\Delta$ denote the combinatorial projection of $\epsilon(w)$ onto $C$. For every $1 \leq i \leq r$, $\Delta$ belongs to $J_i^+$, hence $D_i \leq \Delta$. Then, if $D$ is another diagram admitting any element of $\mathcal{D}$ as a prefix, necessarily $D \in C$, and it follows from Proposition \ref{projection} that there exists a combinatorial geodesic between $\epsilon(w)$ and $D$ passing through $\Delta$, hence $\Delta \leq D$. This proves the $\leq$-minimality and the uniqueness of $\Delta$.

\medskip \noindent
Now, let $D \in \mathcal{M}(\Delta)$. If $J$ denotes the hyperplane associated to $D$, then $J$ separates $\epsilon(w)$ and $\Delta$, and according to Lemma \ref{hyperplan séparant} it must separate $\epsilon(w)$ and $C$. Therefore, either $J$ belongs to $\{J_1, \ldots, J_r \}$ or it separates $\epsilon(w)$ and some $J_k$; equivalently, either $D$ belongs to $\{D_1, \ldots, D_r\}$ or it is a prefix of some $D_k$. Hence $D \in \bigcup\limits_{i=1}^r \mathcal{M}(D_i)$. Thus, we have proved that $\mathcal{M}( \Delta) \subset \bigcup\limits_{D \in \mathcal{D}} \mathcal{M}(D)$; the inclusion $\bigcup\limits_{D \in \mathcal{D}} \mathcal{M}(D) \subset \mathcal{M}( \Delta))$ is clear since any element of $\mathcal{D}$ is a prefix of $\Delta$. $\square$

\medskip \noindent
As a consequence, we deduce that a diagram is determined by the set of its minimal prefixes.

\begin{lemma}\label{prefixesminimaux}
Let $\Delta_1$ and $\Delta_2$ be two diagrams. If $\mathcal{M}(\Delta_1)= \mathcal{M}(\Delta_2)$ then $\Delta_1= \Delta_2$.
\end{lemma}

\noindent
\textbf{Proof.} We begin with some general results. Let $\Delta$ be a diagram.

\begin{fact}\label{cardinal}
$\# \mathcal{M}(\Delta)= \# \Delta$.
\end{fact}

\noindent
Indeed, there exists a bijection between $\mathcal{M}(\Delta)$ and the hyperplanes separating $\epsilon(w)$ and $\Delta$, hence $\# \mathcal{M}(\Delta)=d(\epsilon(w),\Delta)= \# \Delta$. We deduce that

\begin{fact}
$\mathrm{sup}(\mathcal{M}(\Delta))=\Delta$.
\end{fact}

\noindent
Indeed, by definition any element of $\mathcal{M}(\Delta)$ is a prefix of $\Delta$, so that $\mathrm{sup}(\mathcal{M}(\Delta))$ exists (because $\mathcal{M}(\Delta)$ is consistent) and $\mathrm{sup}(\mathcal{M}(\Delta)) \leq \Delta$ by the $\leq$-minimality of $\mathrm{sup}(\mathcal{M}(\Delta))$. On the other hand, we deduce from the previous fact that
\begin{center}
$\# \mathrm{sup}(\mathcal{M}(\Delta)) = | \mathcal{M}( \mathrm{sup}(\mathcal{M}(\Delta))) | = \left| \bigcup\limits_{D \in \mathcal{M}(\Delta)} \mathcal{M}(D) \right| = \# \mathcal{M}(\Delta)= \# \Delta$.
\end{center}
Therefore, we necessarily have $\mathrm{sup}(\mathcal{M}(\Delta))=\Delta$.

\medskip \noindent
Finally, our lemma follows easily, since $\Delta_1= \mathrm{sup}(\mathcal{M}(\Delta_1))= \mathrm{sup}(\mathcal{M}(\Delta_2)) = \Delta_2$. $\square$

\paragraph{Squier cube complexes.} Because the action $D(\mathcal{P},w) \curvearrowright X(\mathcal{P},w)$ is free and properly discontinuous, we know that $D(\mathcal{P},w)$ is isomorphic to the fundamental group of the quotient space $X(\mathcal{P},w)/D(\mathcal{P},w)$. We give now a description of this quotient. We define the \textit{Squier complex} $S(\mathcal{P})$ as the cube complex whose vertices are the words in $\Sigma^+$; whose (oriented) edges can be written as $(a,u \to v,b)$, where $u=v$ or $v=u$ belongs to $\mathcal{R}$, linking the vertices $aub$ and $avb$; and whose $n$-cubes similarly can be written as 
\begin{center}
$(a_1,u_1 \to v_1, \ldots, a_n, u_n \to v_n, a_{n+1})$, 
\end{center}
spanned by the set of vertices $\{ a_1w_1 \cdots a_nw_na_{n+1} \mid w_i=v_i \ \text{or} \ u_i \}$.

Then, there is a natural morphism from the fundamental group of $S(\mathcal{P})$ based at $w$ to the diagram group $D(\mathcal{P},w)$. Indeed, a loop in $S(\mathcal{P})$ based at $w$ is just a series of relations of $\mathcal{R}$ applied to the word $w$ so that the final word is again $w$, and such a derivation may be encoded into a semigroup diagram. The figure below shows an example, where the semigroup presentation is $\mathcal{P}= \langle a,b,c \mid ab=ba, bc=cb,ac=ca \rangle$: 
\begin{center}
\includegraphics[scale=0.5]{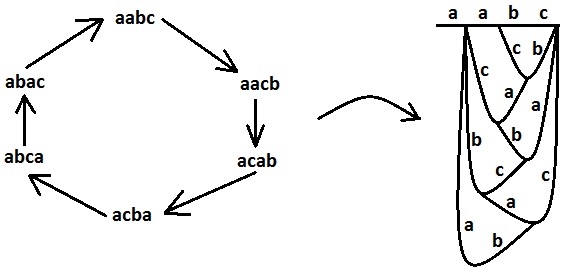}
\end{center}
Thus, this defines a map from the set of loops of $S(\mathcal{P})$ based at $w$ to the set of spherical semigroup diagrams. In fact, the map extends to a morphism which turns out to be an isomorphism:

\begin{thm}\label{iso}
\emph{\cite[Theorem 6.1]{MR1396957}} $D(\mathcal{P},w) \simeq \pi_1( S(\mathcal{P}),w)$.
\end{thm}

For convenience, $S(\mathcal{P},w)$ will denote the connected component of $S(\mathcal{P})$ containing $w$. Notice that two words $w_1,w_2 \in \Sigma^+$ are equal modulo $\mathcal{P}$ if and only if they belong to the same connected component of $S(\mathcal{P})$. Therefore, a consequence of Theorem \ref{iso} is:

\begin{cor}\label{morphism1}
If $w_1,w_2 \in \Sigma^+$ are equal modulo $\mathcal{P}$, then there exists a $(w_2,w_1)$-diagram $\Gamma$ and the map
\begin{center}
$\Delta \mapsto \Gamma \cdot \Delta \cdot \Gamma^{-1}$
\end{center}
induces an isomorphism from $D(\mathcal{P},w_1)$ to $D(\mathcal{P},w_2)$.
\end{cor}

\noindent
As claimed, we notice that the map $\Delta \mapsto \mathrm{bot}(\Delta)$ induces a universal covering $X(\mathcal{P},w) \to S(\mathcal{P},w)$ so that the action of $\pi_1(S(\mathcal{P},w))$ on $X(\mathcal{P},w)$ coincides with the natural action of $D(\mathcal{P},w)$. 

\begin{lemma}
\emph{\cite[Lemma 1.3.5]{arXiv:1505.02053}} The map $\Delta \mapsto \mathrm{bot}(\Delta)$ induces a cellular isomorphism from the quotient $X(\mathcal{P},w)/D(\mathcal{P},w)$ to $S(\mathcal{P},w)$.
\end{lemma}

\paragraph{Changing the base word.} 
When we work on the Cayley graph of a group, because the action of the group is vertex-transitive, we may always suppose that a fixed base point is the identity element up to a conjugation. In our situation, where $X(\mathcal{P},w)$ is the Cayley graph of a groupoid, the situation is slightly different but it is nevertheless possible to make something similar. Indeed, if $\Delta \in X(\mathcal{P},w)$ is a base point, then conjugating by $\Delta$ sends $\Delta$ to a trivial diagram, but which does not necessarily belong to $X(\mathcal{P},w)$: in general, it will belong to $X(\mathcal{P}, \mathrm{bot}(\Delta))$. In the same time, $D(\mathcal{P},w)$ becomes $D(\mathcal{P}, \mathrm{bot}(\Delta))$, which is naturally isomorphic to $D(\mathcal{P},w)$ according to Corollary \ref{morphism1}. Finally, we get a commutative diagram
\begin{displaymath}
\xymatrix{ D(\mathcal{P},w) \ar[rr]^{\text{isomorphism}} \ar[d] & & D(\mathcal{P}, \mathrm{bot}(\Delta)) \ar[d] \\ X(\mathcal{P},w) \ar[rr]^{\text{isometry}} & & X(\mathcal{P},\mathrm{bot}(\Delta)) }
\end{displaymath}
Thus, we may always suppose that a fixed base point is a trivial diagram up to changing the base word, which does not disturb either the group or the cube complex. In particular, a contracting isometry stays a contracting isometry after the process.

\paragraph{Infinite diagrams.} Basically, an infinite diagram is just a diagram with infinitely many cells. To be more precise, if finite diagrams are thought of as specific planar graphs, then an infinite diagram is the union of an increasing sequence finite diagrams with the same top path. Formally,

\begin{definition}
An \emph{infinite diagram} is a formal concatenation $\Delta_1 \circ \Delta_2 \circ \cdots$ of infinitely many diagrams $\Delta_1, \Delta_2, \ldots$ up to the following identification: two formal concatenations $A_1 \circ A_2 \circ \cdots$ and $B_1 \circ B_2 \circ \cdots$ define the same infinite diagram if, for every $i \geq 1$, there exists some $j \geq 1$ such that $A_1 \circ \cdots \circ A_i$ is a prefix of $B_1 \circ  \cdots \circ B_j$, and conversely, for every $i \geq 1$, there exists some $j \geq 1$ such that $B_1 \circ  \cdots \circ B_i$ is a prefix of $A_1 \circ \cdots \circ A_j$. \\
An infinite diagram $\Delta= \Delta_1 \circ \Delta_2 \circ \cdots$ is \emph{reduced} if $\Delta_1 \circ \cdots \circ \Delta_n$ is reduced for every $n \geq 1$. If $\mathrm{top}(\Delta_1)=w$, we say that $\Delta$ is an \emph{infinite $w$-diagram}.
\end{definition}

\noindent
Notice that, according to Lemma \ref{geodesic}, a combinatorial ray is naturally labelled by a reduced infinite diagram. This explains why infinite diagrams will play a central role in the description of the combinatorial boundaries of Farley cube complexes.

\begin{definition}
Let $\Delta=\Delta_1 \circ \Delta_2 \circ \cdots$ and $\Xi = \Xi_1 \circ \Xi_2 \circ \cdots$ be two infinite diagrams. We say that $\Delta$ is a prefix of $\Xi$ if, for every $i \geq 1$, there exists some $j \geq 1$ such that $\Delta_1 \circ \cdots \circ \Delta_i$ is a prefix of $\Xi_1 \circ \cdots \Xi_j$.
\end{definition}

\subsection{Combinatorial boundaries of Farley complexes}\label{section:Fboundary}

\noindent
The goal of this section is to describe the combinatorial boundaries of Farley cube complexes as a set of infinite reduced diagrams. First, we need to introduce a partial order on the set of infinite diagrams. Basically, we will say that a diagram $\Delta_1$ is \emph{almost a prefix} of $\Delta_2$ if there exists a third diagram $\Delta_0$ which is a prefix of $\Delta_1$ and $\Delta_2$ such that all but finitely many cells of $\Delta_1$ belong to $\Delta_0$. Notice that, if we fix a decomposition of a diagram $\Delta$ as a concatenation of atomic diagrams, then to any cell $\pi$ of $\Delta$ corresponds naturally a minimal prefix of $\Delta$, namely the smallest prefix containing $\pi$ (which is precisely the union of all the cells of $\Delta$ which are ``above'' $\pi$; more formally, using the relation introduced by Definition \ref{def:precprec}, this is the union of all the cells $\pi'$ of $\Delta$ satisfying $\pi' \prec \prec \pi$). Therefore, alternatively we can describe our order in terms of minimal prefixes.

\begin{definition}
Let $\Delta_1$ and $\Delta_2$ be two infinite reduced diagrams. Then,
\begin{itemize}
	\item $\Delta_1$ is \emph{almost a prefix} of $\Delta_2$ if $\mathcal{M}(\Delta_1) \underset{a}{\subset} \mathcal{M}(\Delta_2)$;
	\item $\Delta_1$ and $\Delta_2$ are \emph{commensurable} if $\mathcal{M}(\Delta_1) \underset{a}{=} \mathcal{M}(\Delta_2)$.
\end{itemize}
\end{definition}

\noindent
Notice that two infinite reduced diagrams are commensurable if and only if each one is almost a prefix of the other. Our model for the combinatorial boundaries of Farley cube complexes will be the following:

\begin{definition} 
If $\mathcal{P}= \langle \Sigma \mid \mathcal{R} \rangle$ is semigroup presentation and $w \in \Sigma^+$ a base word, let $(\partial (\mathcal{P},w), \underset{a}{<})$ denote the poset of the infinite reduced $(w,\ast)$-diagrams, up to commensurability, ordered by the relation $\underset{a}{<}$ of being almost a prefix.
\end{definition}

\noindent
Recall that to any combinatorial ray $r \subset X(\mathcal{P},w)$ starting from $\epsilon(w)$ is naturally associated an infinite reduced $(w,\ast)$-diagram $\Psi(r)$.

\begin{thm}\label{cb-Farley}
The map $\Psi$ induces a poset-isomorphism $\partial^c X(\mathcal{P},w) \to \partial(\mathcal{P},w)$.
\end{thm}

\noindent
The main tool to prove our theorem is provided by the following proposition, which is an infinite analogue of Proposition \ref{finite-consistent}. Recall that a collection of (finite) diagrams $\mathcal{D}$ is consistent if, for every $\Delta_1,\Delta_2 \in \mathcal{D}$, there exists a third diagram $\Delta$ satisfying $\Delta_1,\Delta_2 \leq \Delta$. 

\begin{prop}
Let $\mathcal{D}$ be a countable consistent collection of finite reduced minimal $(w,\ast)$-diagrams. Then there exists a unique $\leq$-minimal reduced $(w,\ast)$-diagram admitting any element of $\mathcal{D}$ as a prefix; let $\mathrm{sup}( \mathcal{D})$ denote this diagram. Furthermore, $\mathcal{M}( \mathrm{sup}(\mathcal{D})) = \bigcup\limits_{D \in \mathcal{D}} \mathcal{M}(D)$.
\end{prop}

\noindent
\textbf{Proof.} Let $\mathcal{D}= \{ D_1, D_2, \ldots \}$, with $D_i \neq D_j$ if $i \neq j$. For each $i \geq 1$, let $H_i$ denote the hyperplane of $X(\mathcal{P},w)$ associated to $D_i$. Because $\mathcal{D}$ is consistent, the halfspaces $H_i^+$ and $H_j^+$ intersect for every $i \geq 1$. Therefore, for every $k \geq 1$, the intersection $I_k = \bigcap\limits_{i=1}^k H_i^+$ is non-empty; let $\Delta_k$ denote the combinatorial projection of $\epsilon(w)$ onto $I_k$. As a consequence of Proposition \ref{projection}, there exists a combinatorial ray starting from $\epsilon(w)$ and passing through each $\Delta_k$. In particular, for every $k \geq 1$, $\Delta_k$ is a prefix of $\Delta_{k+1}$, so that it makes sense to define the infinite diagram $\Delta$ as the limit of the sequence $(\Delta_k)$. Notice that, for every $k \geq 1$, $\Delta_k \in H_k^+$ so $D_k \leq \Delta_k$ and a fortiori $D_k \leq \Delta$. 

\medskip \noindent
Now, let $\Xi$ be another reduced $w$-diagram admitting each $D_k$ as a prefix. Fixing a decomposition of $\Xi$ as an infinite concatenation of atomic diagrams, say $A_1 \circ A_2 \circ \cdots$, we associate to $\Xi$ a combinatorial ray $r$ starting from $\epsilon(w)$ defined as the path 
$$\epsilon(w), \ A_1, \ A_1 \circ A_2, \ A_1 \circ A_2 \circ A_3 \ldots$$
Let $n \geq 1$. Because $D_1, \ldots, D_n$ are prefixes of $\Xi$, there exists some $k \geq 1$ such that $r(k) \in I_n$. On the other hand, according to Proposition \ref{projection}, there exists a combinatorial geodesic between $\epsilon(w)$ and $r(k)$ passing through $\Delta_n$, hence $\Delta_n \leq r(k)$. Therefore, we have proved that, for every $n \geq 1$, $\Delta_n$ is a prefix of $\Xi$: this precisely means that $\Delta$ is a prefix of $\Xi$. The $\leq$-minimality and the uniqueness of $\Delta$ is thus proved. 

\medskip \noindent
Because the $D_k$'s are prefixes of $\Delta$, the inclusion $\bigcup\limits_{D\in \mathcal{D}} \mathcal{M}(D) \subset \mathcal{M} (\Delta)$ is clear. Conversely, let $D \in \mathcal{M}(\Delta)$; let $J$ denote the associated hyperplane. Because $D$ is a prefix of $\Delta$, it is a prefix of $\Delta_k$ for some $k \geq 1$; and because $\Delta_k$ is the combinatorial projection of $\epsilon(w)$ onto $I_k$, we deduce that $J$ separates $\epsilon(w)$ and $I_k$. On the other hand, $I_k= \bigcap\limits_{i=1}^k H_i^+$ so there exists some $1 \leq i \leq k$ such that either $J=H_i$ or $J$ separates $\epsilon(w)$ and $H_i$; equivalently, $D$ is a prefix of $D_i$. We have proved $\mathcal{M}(\Delta) \subset \bigcup\limits_{D\in \mathcal{D}} \mathcal{M}(D)$. $\square$ 

\medskip \noindent
The following lemma states how the operation $\mathrm{sup}$ behaves with respect to the inclusion, the almost-inclusion and the almost-equality.

\begin{lemma}\label{lem:supalmost}
Let $\mathcal{D}_1$ and $\mathcal{D}_2$ be two countable consistent collections of reduced $(w,\ast)$-diagrams. 
\begin{itemize}
	\item[(i)] If $\mathcal{D}_1 \subset \mathcal{D}_2$ then $\mathrm{sup}(\mathcal{D}_1) \leq \mathrm{sup}(\mathcal{D}_2)$;
	\item[(ii)] if $\mathcal{D}_1 \underset{a}{\subset} \mathcal{D}_2$ then $\mathrm{sup}(\mathcal{D}_1)$ is almost a prefix of $\mathrm{sup}(\mathcal{D}_2)$;
	\item[(iii)] if $\mathcal{D}_1 \underset{a}{=} \mathcal{D}_2$ then $\mathrm{sup}(\mathcal{D}_1)$ and $\mathrm{sup}(\mathcal{D}_2)$ are commensurable.
\end{itemize}
\end{lemma}

\noindent
\textbf{Proof.} If $\mathcal{D}_1 \subset \mathcal{D}_2$ then $\mathrm{sup}(\mathcal{D}_2)$ admits any element of $\mathcal{D}_1$ as a prefix, hence $\mathrm{sup}(\mathcal{D}_1) \leq \mathrm{sup}(\mathcal{D}_2)$ by the $\leq$-minimality of $\mathrm{sup}(\mathcal{D}_1)$.

\medskip \noindent
If $\mathcal{D}_1 \underset{a}{\subset} \mathcal{D}_2$ then we can write
\begin{center}
$\mathcal{M}(\mathrm{sup}(\mathcal{D}_1)) = \bigcup\limits_{D \in \mathcal{D}_1} \mathcal{M}(D) = \left( \bigcup\limits_{D \in \mathcal{D}_1 \cap \mathcal{D}_2} \mathcal{M}(D) \right) \sqcup \left( \bigcup\limits_{D \in \mathcal{D}_1 \backslash \mathcal{D}_2} \mathcal{M}(D) \right)$.
\end{center}
On the other hand, $\mathcal{D}_1 \backslash \mathcal{D}_2$ is finite, because $\mathcal{D}_1 \underset{a}{\subset} \mathcal{D}_2$, and $\mathcal{M}(D)$ is finite for any finite diagram $D$, so $\bigcup\limits_{D \in \mathcal{D}_1 \backslash \mathcal{D}_2} \mathcal{M}(D)$ is finite. This proves that $\mathcal{M}(\mathrm{sup}(\mathcal{D}_1)) \underset{a}{\subset} \mathcal{M}(\mathrm{sup}(\mathcal{D}_2))$, ie., $\mathrm{sup}(\mathcal{D}_1)$ is almost a prefix of $\mathrm{sup}(\mathcal{D}_2)$.

\medskip \noindent
If $\mathcal{D}_1 \underset{a}{=} \mathcal{D}_2$, we deduce by applying the previous point twice that $\mathrm{sup}(\mathcal{D}_1)$ is almost a prefix of $\mathrm{sup}(\mathcal{D}_2)$ and vice-versa. Therefore, $\mathrm{sup}(\mathcal{D}_1)$ and $\mathrm{sup} (\mathcal{D}_2)$ are commensurable. $\square$

\medskip \noindent
Finally, Theorem \ref{cb-Farley} will be essentially a consequence of the previous lemma and the following observation. If $r$ is a combinatorial ray in a Farley complex, let $\mathcal{D}(r)$ denote the set of the minimal diagrams which are associated to the hyperplanes of $\mathcal{H}(r)$. Then

\begin{lemma}\label{lem:psisup}
$\Psi(r)= \mathrm{sup}(\mathcal{D}(r))$.
\end{lemma}

\noindent
\textbf{Proof.} Let $A_1,A_2,A_3 \ldots$ be the sequence of atomic diagrams such that 
\begin{center}
$r=(\epsilon(w), \ A_1, \ A_1 \circ A_2, \ A_1 \circ A_2 \circ A_3, \ldots)$.
\end{center}
In particular, $\Psi(r)=A_1 \circ A_2 \circ \cdots$ and $\mathcal{M}(\Psi(r))= \bigcup\limits_{n \geq 1} \mathcal{M}(A_1 \circ \cdots \circ A_n)$. On the other hand, any $D \in \mathcal{D}(r)$ must be a prefix of $A_1 \circ \cdots \circ A_n$ for some $n \geq 1$, and conversely, any minimal prefix of some $A_1 \circ \cdots \circ A_n$ must belong to $\mathcal{D}(r)$. Therefore,
\begin{center}
$\mathcal{M}(\Psi(r)) = \bigcup\limits_{n \geq 1} \mathcal{M}(A_1 \circ \cdots \circ A_n) = \mathcal{M}(\mathrm{sup}(\mathcal{D}(r)))$.
\end{center}
If $\Delta$ is a finite prefix of $\Psi(r)$, then $\mathcal{M}(\Delta) \subset \mathcal{M}(\Psi(r))= \mathcal{M}(\mathrm{sup} (\mathcal{D}(r)))$, so $\Delta= \mathrm{sup}(\mathcal{M}(\Delta))$ must be a prefix of $\mathrm{sup}(\mathcal{D}(r))$ as well. Therefore, $\Psi(r) \leq \mathcal{M}(\mathrm{sup}(\mathcal{D}(r)))$. Similarly, we prove that $\mathcal{M}(\mathrm{sup}(\mathcal{D}(r))) \leq \Psi(r)$, hence $\Psi(r)= \mathcal{M}(\mathrm{sup}(\mathcal{D}(r)))$. $\square$

\medskip \noindent
\textbf{Proof of Theorem \ref{cb-Farley}.} First, we have to verify that, if $r_1$ and $r_2$ are two equivalent combinatorial rays starting from $\epsilon(w)$, then $\Psi(r_1) \underset{a}{=} \Psi(r_2)$. Because $r_1$ and $r_2$ are equivalent, i.e., $\mathcal{H}(r_1) \underset{a}{=} \mathcal{H}(r_2)$, we know that $\mathcal{D}(r_1) \underset{a}{=} \mathcal{D}(r_2)$, which implies 
\begin{center}
$\Psi(r_1)= \mathrm{sup}(\mathcal{D}(r_1)) \underset{a}{=} \mathrm{sup}(\mathcal{D}(r_2)) = \Psi(r_2)$
\end{center}
according to Lemma \ref{lem:supalmost} and Lemma \ref{lem:psisup}. Therefore, $\Psi$ induces to a map $\partial^c X(\mathcal{P},w) \to \partial(\mathcal{P},w)$. For convenience, this map will also be denoted by $\Psi$. 

\medskip \noindent
If $r_1 \prec r_2$, then $\mathcal{H}(r_1) \underset{a}{\subset} \mathcal{H}(r_2)$ hence 
\begin{center}
$\Psi(r_1)= \mathrm{sup}(\mathcal{D}(r_1)) \underset{a}{<} \mathrm{sup}(\mathcal{D}(r_2)) = \Psi(r_2)$.
\end{center}
Therefore, $\Psi$ is a poset-morphism. 

\medskip \noindent
Let $\Delta$ be an infinite reduced $(w,\ast)$-diagram. Let $A_1, A_2, \ldots$ be a sequence of atomic diagrams such that $\Delta=A_1 \circ A_2 \circ \cdots$. Let $\rho$ denote the combinatorial ray $(\epsilon(w), \ A_1, \ A_1 \circ A_2, \ldots)$. Now it is clear that $\Psi(\rho)=A_1 \circ A_2 \circ \cdots = \Delta$, so $\Psi$ is surjective.

\medskip \noindent
Let $r_1$ and $r_2$ be two combinatorial rays starting from $\epsilon(w)$ and satisfying $\Psi(r_1) \underset{a}{=} \Psi(r_2)$. It follows from Lemma \ref{lem:psisup} that $\mathrm{sup}(\mathcal{D}(r_1)) \underset{a}{=} \mathrm{sup}(\mathcal{D}(r_2))$. This means
\begin{center}
$\bigcup\limits_{D \in \mathcal{D}(r_1)} \mathcal{M}(D)= \mathcal{M}(\mathrm{sup}(\mathcal{D}(r_1)) \underset{a}{=} \mathcal{M}(\mathrm{sup}(\mathcal{D}(r_2)))= \bigcup\limits_{D \in \mathcal{D}(r_2)} \mathcal{M}(D)$.
\end{center}
Notice that $\mathcal{D}(r_2)$ is stable under taking a prefix. Indeed, if $P$ is a prefix of some $\Delta \in \mathcal{D}(r_2)$ then the hyperplane $J$ associated to $P$ separates $\epsilon(w)$ and the hyperplane $H$ associated to $\Delta$; because $H \in \mathcal{H}(r_2)$ and $r_2(0)= \epsilon(w)$, necessarily $J \in \mathcal{H}(r_2)$, hence $P \in \mathcal{D}(r_2)$. As a consequence, we deduce that
\begin{center}
$\mathcal{D}(r_1) \backslash \mathcal{D}(r_2) \subset \left( \bigcup\limits_{D \in \mathcal{D}(r_1)} \mathcal{M}(D) \right) \backslash \left( \bigcup\limits_{D \in \mathcal{D}(r_2)} \mathcal{M}(D) \right)$.
\end{center}
We conclude that $\mathcal{D}(r_1) \backslash \mathcal{D}(r_2)$ is finite. We prove similarly that $\mathcal{D}(r_2) \backslash \mathcal{D}(r_1)$ is finite. Therefore, we have $\mathcal{D}(r_1) \underset{a}{=} \mathcal{D}(r_2)$, which implies $\mathcal{H}(r_1) \underset{a}{=} \mathcal{H}(r_2)$, so that $r_1$ and $r_2$ are equivalent. We have proved that $\Psi$ is injective. $\square$

\begin{ex}\label{Thompson}
Let $\mathcal{P}= \langle x \mid x=x^2 \rangle$ be the semigroup presentation usually associated to Thompson's group $F$. An atomic diagram is said \emph{positive} if the associated relation is $x \to x^2$, and \emph{negative} if the associated relation is $x^2 \to x$; by extension, a diagram is \emph{positive} (resp. \emph{negative}) if it is a concatenation of positive (resp. negative) atomic diagrams. Notice that any reduced infinite $(x,\ast)$-diagram has an infinite positive prefix. On the other hand, there exists a reduced infinite postive $(x,\ast)$-diagram containing all the possible reduced infinite positive $(x,\ast)$-diagrams as prefixes: it corresponds to the derivation obtained by replacing at each step each $x$ by $x^2$. Therefore, $\partial(\mathcal{P},x)$ does not contain any isolated vertex. We deduce from Theorem \ref{contracting isometry and boundary} that the automorphism group $\mathrm{Aut}(X(\mathcal{P},x))$ does not contain any contracting isometry.
\end{ex}

\begin{ex}
Let $\mathcal{P}= \langle a,b,p_1,p_2,p_3 \mid a=ap_1, b=p_1b, p_1=p_2, p_2=p_3, p_3=p_1 \rangle$ be the semigroup presentation usually associated to the lamplighter group $\mathbb{Z} \wr \mathbb{Z}$. We leave as an exercice to notice that $\partial (\mathcal{P},ab)$ has no isolated point. Therefore, $\mathrm{Aut}(X(\mathcal{P},ab))$ does not contain any contracting isometry according to Theorem \ref{contracting isometry and boundary}.
\end{ex}

\subsection{Contracting isometries}

\noindent
In this section, we fix a semigroup presentation $\mathcal{P}= \langle \Sigma \mid \mathcal{R} \rangle$ and base word $w \in \Sigma^+$. Our goal is to determine precisely when a spherical diagram $g \in D(\mathcal{P},w)$ induces a contracting isometry on the Farley complex $X(\mathcal{P},w)$. For convenience, we will suppose that $g$ is \emph{absolutely reduced}, ie., for every $n \geq 1$ the concatenation of $n$ copies of $g$ is reduced; this assumption is not really restrictive since, according to \cite[Lemma 15.10]{MR1396957}, any spherical diagram is conjugated to an absolutely reduced spherical diagram (possibly for a different word base). In particular, we may define the infinite reduced $(w,\ast)$-diagram $g^{\infty}$ as $g \circ g \circ \cdots$. Our main criterion is:

\begin{thm}\label{prop1}
Let $g \in D(\mathcal{P},w) - \{ \epsilon(w) \}$ be absolutely reduced. Then $g$ is a contracting isometry of $X(\mathcal{P},w)$ if and only if the following two conditions are satisfied:
\begin{itemize}
	\item[(i)] $g^{\infty}$ does not contain any infinite proper prefix;
	\item[(ii)] if $\Delta$ is a reduced diagram with $g^{\infty}$ as a prefix, then $\Delta$ is commensurable to $g^{\infty}$. 
\end{itemize}
\end{thm}

\noindent
This theorem will be essentially a consequence of the following two lemmas.

\begin{lemma}\label{axis}
Let $g \in D(\mathcal{P},w) - \{ \epsilon(w) \}$ be absolutely reduced. Choose a combinatorial geodesic $[\epsilon(w),g]$ between $\epsilon(w)$ and $g$, and set $\gamma= \bigcup\limits_{n \in \mathbb{Z}} g^n \cdot [\epsilon(w) ,g]$. Then $\gamma$ is a bi-infinite combinatorial geodesic on which $\langle g \rangle$ acts by translation. Moreover, $\gamma(+ \infty)= g^{\infty}$.
\end{lemma}

\noindent
\textbf{Proof.} Notice that $\gamma$ is a bi-infinite combinatorial path passing through $g^n$ for every $n \in \mathbb{Z}$. Therefore, in order to deduce that $\gamma$ is a geodesic, it is sufficient to show that, for every $n,m \in \mathbb{Z}$, the length of $\gamma$ between $g^n$ and $g^{n+m}$ is equal to $d(g^n,g^{n+m})$; but this length is precisely
\begin{center}
$m \cdot \mathrm{length}([\epsilon(w),g])=m \cdot d(\epsilon(w),g) = m \cdot \# (g)$,
\end{center}
and, on the other hand, because $g$ is absolutely reduced, 
\begin{center}
$d(g^n,g^{n+m})= \# (g^{-n}\cdot g^{n+m})= \# (g^m)=m \cdot \# (g)$.
\end{center}
We conclude that $\gamma$ is a combinatorial geodesic. The fact that $\langle g \rangle$ acts on $\gamma$ by translation is clear. Finally, we deduce that $\gamma(+ \infty)= g^{\infty}$ from the fact that, for any $k \geq 1$, we have $\gamma(k \cdot \#(g))=g^k$. $\square$

\begin{definition}
Let $\Delta$ be a possibly infinite $(w,\ast)$-diagram. The \emph{support} of $\Delta$ is the set of maximal subpaths of $w$ whose edges belong to the top path of some cell of $\Delta$. 
\end{definition}

\begin{lemma}\label{cellprec}
Let $\Delta_1,\ldots, \Delta_n$ be copies of an absolutely reduced spherical $(w,w)$-diagram $\Delta$, where $n> |w|$. If $e \subset \Delta_1 \circ \cdots \circ \Delta_n$ which belongs to the bottom path of some cell of $\Delta_1$, then $e$ belongs to the top path of some cell in $\Delta_1 \circ \cdots \circ \Delta_n$.
\end{lemma}

\noindent
\textbf{Proof.} Say that an edge of $\Delta_i$ is a \emph{cutting edge} if it belongs to the bottom path of some cell of $\Delta_i$ but does not belong to the top path of some cell in $\Delta_1 \circ \cdots \circ \Delta_n$. Suppose that there exists a cutting edge $e_1$ in $\Delta_1$. 

\medskip \noindent
In particular, $e_1$ is an edge of the top path of $\Delta_2$, so $e_1$ does not belong to the support of $\Delta_2$. We deduce that $\Delta_2$ decomposes as a sum $\Phi_2 + \epsilon + \Psi_2$, where $\epsilon$ is a trivial diagram with top and bottom paths equal to $e_1$. Let $\Phi_1+ \epsilon+\Psi_1$ denote the same decomposition of $\Delta_1$. Because $e_1$ belongs to the bottom path of a cell of $\Delta_1$, two cases may happen: either $e_1$ belongs to $\mathrm{bot}(\Phi_1)$, so that $\mathrm{top}(\Phi_2) \subsetneq \mathrm{bot}(\Phi_1)$, or $e_1$ belongs to $\mathrm{bot}(\Psi_1)$, so that $\mathrm{top}(\Psi_2) \subsetneq \mathrm{bot}(\Psi_1)$. Say we are in the former case, the latter case being similar. Now, because $\Delta_2$ and $\Delta_1$ are two copies of the same diagram, to the edge $e_1$ of $\Delta_1$ corresponds an edge $e_2$ of $\Delta_2$. Notice that, because $e_1$ belongs to $\mathrm{bot}(\Phi_1)$, necessarily $e_2 \neq e_1$. Moreover, the product $\Delta^{-1} \circ \Delta_1 \circ \cdots \circ \Delta_n$ naturally reduces to a copy of $\Delta_1 \circ \cdots \circ \Delta_{n-1}$, and this process sends the edge $e_2$ to the edge $e_1$. Therefore, because $e_1$ is a cutting edge in $\Delta_1$, we deduce that $e_2$ is a cutting edge in $\Delta_2$. 

\medskip \noindent
By iterating this construction, we find a cutting edge $e_i$ in $\Delta_i$ for every $1 \leq i \leq n$, where $e_i \neq e_j$ provided that $i \neq j$. On the other hand, the path $\mathrm{bot}(\Delta_1 \circ \cdots \circ \Delta_n)$ necessarily contains all these edges, hence $n \leq |w|$. Consequently, $\Delta_1$ cannot contain a cutting edge if $n> |w|$. $\square$

\medskip \noindent
\textbf{Proof of Theorem \ref{prop1}.} Let $\gamma$ be the combinatorial axis of $g$ given by Lemma \ref{axis}. According to Theorem \ref{contracting isometry and boundary}, $g$ is a contracting isometry if and only if $\gamma(+ \infty)$ is isolated in $\partial^c X(\mathcal{P},w)$. Thus, using the isomorphism given by Theorem \ref{cb-Farley}, $g$ is a contracting isometry if and only if $g^{\infty}$ is isolated in $\partial(\mathcal{P},w)$, ie.,
\begin{itemize}
	\item every infinite prefix of $g^{\infty}$ is commensurable to $g^{\infty}$;
	\item if $\Delta$ is a reduced diagram with $g^{\infty}$ as a prefix, then $\Delta$ is commensurable to $g^{\infty}$. 
\end{itemize}
Therefore, to conclude it is sufficient to prove that a proper infinite prefix $\Delta$ of $g^{\infty}$ cannot be commensurable to $g^{\infty}$. Because $\Delta$ is a proper prefix, there exists a cell $\pi_1$ of $g^{\infty}$ which does not belong to $\Delta$. Now, by applying Lemma \ref{cellprec} to a given edge of the bottom path of $\pi_1$, we find a cell $\pi_2$ whose top path intersects the bottom path of $\pi_1$ along at least one edge; by applying Lemma \ref{cellprec} to a given edge of the bottom path of $\pi_2$, we find a cell $\pi_3$ whose top path intersects the bottom path of $\pi_2$ along at least one edge; and so on. Finally, we find an infinite sequence of cells $\pi_2, \pi_3, \ldots$ such that, for every $i \geq 2$, the top path of $\pi_i$ has a common edge with the bottom path of $\pi_{i-1}$. For every $i \geq 1$, let $\Xi_i$ denote the smallest (finite) prefix of $g^{\infty}$ which contains the cell $\pi_i$. This is a minimal diagram, and by construction $\Xi_1, \Xi_2, \ldots \notin \mathcal{M}(\Delta)$, which proves that $\Delta$ and $g^{\infty}$ are not commensurable. $\square$

\begin{ex}
Let $\mathcal{P}= \langle a,b,p \mid a=ap,b=pb \rangle$ and let $g \in D(\mathcal{P},ab)$ be the following spherical diagram:
\begin{center}
\includegraphics[scale=0.6]{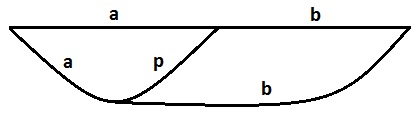}
\end{center}
Then $g^{\infty}$ clearly contains a proper infinite prefix. Therefore, $g$ is not a contracting isometry of $X(\mathcal{P},ab)$.
\begin{center}
\includegraphics[scale=0.6]{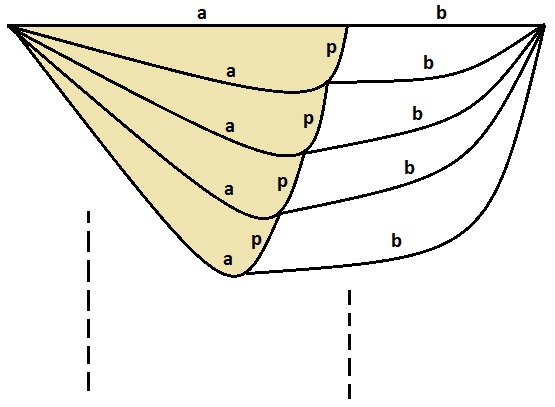}
\end{center}
\end{ex}

\begin{ex}
Let $\mathcal{P}= \langle a,b,c \mid a=b,b=c,c=a \rangle$ and let $g \in D(\mathcal{P},a^2)$ be the following spherical diagram:
\begin{center}
\includegraphics[scale=0.6]{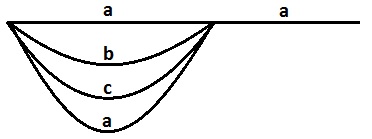}
\end{center}
Then $g^{\infty}$ is a prefix of the diagram $\Delta$ below, but $g^{\infty}$ and $\Delta$ are not commensurable. Therefore, $g$ is not a contracting isometry of $X(\mathcal{P},a^2)$. 
\begin{center}
\includegraphics[scale=0.6]{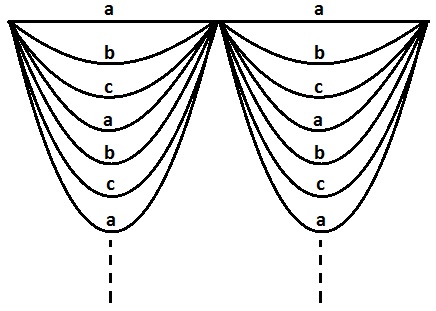}
\end{center}
\end{ex}

\begin{ex}
Let $\mathcal{P}= \langle a,b,c,d \mid ab=ac, cd=bd \rangle$ and let $g \in D(\mathcal{P},ab)$ be the following spherical diagram:
\begin{center}
\includegraphics[scale=0.6]{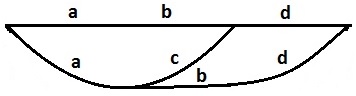}
\end{center}
So the infinite diagram $g^{\infty}$ looks like
\begin{center}
\includegraphics[scale=0.6]{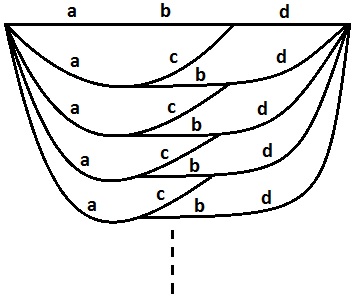}
\end{center}
We can notice that $g^{\infty}$ does not contain a proper infinite prefix that any diagram containing $g^{\infty}$ as a prefix is necessarily equal to $g^{\infty}$.
\end{ex}

\medskip \noindent
The criterion provided by Theorem \ref{prop1} may be considered as unsatisfactory, because we cannot draw an infinite diagram in practice (although it seems to be sufficient for spherical diagrams with only few cells, as suggested by the above examples). We conclude this section by stating and proving equivalent conditions to the assertions $(i)$ and $(ii)$ of Theorem \ref{prop1}. For convenience, we will suppose our spherical diagram absolutely reduced and \emph{normal}, ie., the factors of every decomposition as a sum are spherical; this assumption is not really restrictive since, according to \cite[Lemma 15.13]{MR1396957}, every spherical diagram is conjugated to a normal absolutely reduced spherical diagram (possibly for a different word base), and this new diagram can be found ``effectively'' \cite[Lemma 15.14]{MR1396957}. We begin with the condition $(i)$. The following definition will be needed.

\begin{definition}\label{def:precprec}
Given a diagram $\Delta$, we introduce a partial relation $\prec$ on the set of its cells, defined by: if $\pi, \pi'$ are two cells of $\Delta$, then $\pi \prec \pi'$ holds if the intersection between the bottom path of $\pi$ and the top path of $\pi'$ contains at least one edge. Let $\prec \prec$ denote the transitive closure of $\prec$. 
\end{definition}

\begin{lemma}\label{lem1}
Let $g \in D(\mathcal{P},w) - \{ \epsilon(w) \}$ be absolutely reduced. Suppose that we can write $w=xyz$ where $\mathrm{supp}(g)= \{ y \}$. The following two conditions are equivalent:
\begin{itemize}
	\item[(i)] $\partial(\mathcal{P},x)$ and $\partial(\mathcal{P},z)$ are empty;
	\item[(ii)] if $\Delta$ is a reduced diagram with $g^{\infty}$ as a prefix, then $\Delta$ is commensurable to $g^{\infty}$. 
\end{itemize}
\end{lemma}

\noindent
\textbf{Proof.} First, notice that $g= \epsilon(x)+ g_0 + \epsilon(z)$ for some $(y,y)$-diagram $g_0$. This assertion clearly holds for some $(y,\ast)$-diagram $g_0$, but, because $g$ is a $(w,w)$-diagram, the equality
\begin{center}
$xyz= w = \mathrm{top}(g) = \mathrm{bot}(g) = x \cdot \mathrm{bot}(g_0) \cdot z$
\end{center}
holds in $\Sigma^+$. Therefore, $\mathrm{bot}(g_0)=y$, and we conclude that $g_0$ is indeed a $(y,y)$-diagram.

\medskip \noindent
Now, we prove $(i) \Rightarrow (ii)$. Let $\Delta$ be a reduced diagram containing $g^{\infty}$ as a prefix. Suppose that $\pi$ is a cell of $\Delta$ such that there exists a cell $\pi_0$ of $g^{\infty}$ satisfying $\pi_0 \prec \prec \pi$. So there exists a chain of cells $\pi_0 \prec \pi_1 \prec \cdots \prec \pi_n= \pi$. By definition, $\pi_1$ share an edge with $\pi_0$. On the other hand, it follows from Lemma \ref{cellprec} that this edge belongs to the top path of a cell of $g^{\infty}$, hence $\pi_1 \subset g^{\infty}$. By iterating this argument, we conclude that $\pi$ belongs to $g^{\infty}$. Therefore, if $\Xi$ denotes the smallest prefix of $\Delta$ containing a given cell which does not belong to $g^{\infty}$, then $\mathrm{supp}(\Xi)$ is included into $x$ or $z$. 

\medskip \noindent
As a consequence, there exist a $(x,\ast)$-diagram $\Phi$ and $(z,\ast)$-diagram $\Psi$ such that $\Delta= \Phi + g_0^{\infty} + \Psi$. On the other hand, because $\partial(\mathcal{P},x)$ and $\partial(\mathcal{P},y)$ are empty, necessarily $\# \Phi, \# \Psi <+ \infty$. Therefore,
\begin{center}
$|\mathcal{M}(\Delta) \backslash \mathcal{M}(g^{\infty})|= |\mathcal{M}(\Phi+ \epsilon(yz)) \sqcup \mathcal{M}(\epsilon(xy)+ \Psi)| = \# \Phi + \# \Psi <+ \infty$,
\end{center}
where we used Fact \ref{cardinal}. Therefore, $\Delta$ is commensurable to $g^{\infty}$.

\medskip \noindent 
Conversely, we prove $(ii) \Rightarrow (i)$. Let $\Phi$ be a reduced $(x,\ast)$-diagram and $\Psi$ a reduced $(z,\ast)$-diagram. Set $\Delta= \Phi + g_0^{\infty} + \Psi$. This is a reduced diagram containing $g^{\infty}$ as a prefix. As above, we have
\begin{center}
$|\mathcal{M}(\Delta) \backslash \mathcal{M}(g^{\infty})|= \# \Phi + \# \Psi$.
\end{center}
Because $\Delta$ must be commensurable to $g^{\infty}$, we deduce that $\Phi$ and $\Psi$ are necessarily finite. We conclude that $\partial (\mathcal{P},x)$ and $\partial (\mathcal{P},z)$ are empty. $\square$

\medskip \noindent
Now, we focus on the condition $(ii)$ of Theorem \ref{prop1}. We first introduce the following definition; explicit examples are given at the end of this section.

\begin{definition}
Let $g \in D(\mathcal{P},w)- \{ \epsilon(w) \}$ be absolutely reduced. A subpath $u$ of $w$ is \emph{admissible} if there exists a prefix $\Delta \leq g$ such that $\mathrm{supp}(\Delta)= \{ u \}$ and $\mathrm{bot}(\Delta) \cap \mathrm{bot}(g)$ is a non-empty connected subpath of $\mathrm{bot}(g)$; the subpath of $w$ corresponding to $\mathrm{bot}(\Delta) \cap \mathrm{bot}(g)$ is a \emph{final subpath} associated to $u$; it is \emph{proper} if this is not the whole $\mathrm{bot}(g)$. Given a subpath $u$, we define its \emph{admissible tree} $T_g(u)$ as the rooted tree constructed inductively as follows:
\begin{itemize}
	\item the root of $T_g(u)$ is labeled by $u$;
	\item if $v$ is a non-admissible subpath of $w$ labelling a vertex of $T_g(u)$, then $v$ has only one descendant, labeled by the symbol $\emptyset$;
	\item if $v$ is an admissible subpath of $w$ labelling a vertex of $T_g(u)$, then the descendants of $v$ correspond to the proper final subpaths associated to $u$.
\end{itemize}
If $w$ has length $n$, we say that $T_g(u)$ is \emph{deep} if it contains at least $\alpha(w)+1$ generations, where $\alpha(w)=\frac{n(n+1)}{2}$ is chosen so that $w$ contains at most $\alpha(w)$ subpaths.
%then the \emph{finite admissible tree} $T_g(u)$ is the rooted subtree of $T_g^{\infty}(u)$ corresponding to the $\alpha(w)+1$ first generations, where $\alpha(w)=\frac{n(n-1)}{2}$ is chosen so that $w$ contains at most $\alpha(w)$ subpaths.
\end{definition}

\begin{lemma}\label{lem2}
Let $g \in D(\mathcal{P},w) - \{ \epsilon(w) \}$ be absolutely reduced. The following two conditions are equivalent:
\begin{itemize}
	\item[(i)] $g^{\infty}$ contain an infinite proper prefix;
	\item[(ii)] there exists a subpath $u$ of $w$ such that $T_g(u)$ is deep.
	%contains twice a proper subpath of $\mathrm{supp}(g)$ in a same lineage.
\end{itemize}
\end{lemma}

\noindent
\textbf{Proof.} For convenience, let $g^{\infty}= g_1 \circ g_2 \circ \cdots$ where each $g_i$ is a copy of $g$.

\medskip \noindent
We first prove $(i) \Rightarrow (ii)$. Suppose that $g^{\infty}$ contains an infinite propre prefix $\Xi$. For every $i \geq 1$, $\Xi$ induces a prefix $\Xi_i$ of $g_i$. Up to taking an infinite prefix of $\Xi$, we can suppose that $\mathrm{supp}(\Xi_i)$ and $\mathrm{bot}(\Xi_i) \cap \mathrm{bot}(g_i)$ are connected subpaths of $\mathrm{top}(g_i)$ and $\mathrm{bot}(g_i)$ respectively. Thus, if $\xi_i$ denotes the word labeling $\mathrm{supp}(\Xi_i)$ for every $i\geq 1$, then the sequence $(\xi_i)$ defines an infinite descending ray in the rooted tree $T_g(\xi_1)$. A fortiori, $T_g(\xi_1)$ is deep.

%\medskip \noindent
%Notice that, if there exists some $j \geq 1$ such that $\mathrm{supp}(g)=\{\xi_j \}$, then $\mathrm{supp}(\xi_i)$ for every $i \leq j$. Therefore, by replacing $\Xi$ with $g^{-n} \cdot \Xi$ for some $n \geq 1$ if necessary, we may suppose without loss of generality that $\xi_i$ is a proper subword of $\mathrm{supp}(g)$ for every $i \geq 1$. Because $w$ contains at most $\alpha(w)$ subwords, necessarily $\{ \xi_i \mid 1 \leq i \leq \alpha(w)+1 \}$ contains twice the same word. Finally, $T_g(\xi_1)$ contains twice a proper subword of $\mathrm{supp}(g)$ in a same lineage.

\medskip \noindent
Conversely, we prove $(ii) \Rightarrow (i)$. 
Suppose that there exists a subword $u$ of $w$ such that $T_g(u)$ is deep. Because $w$ contains at most $\alpha(w)$ subwords, necessarily
%Suppose that there exists a subword $u$ of $w$ such that 
$T_g(u)$ contains twice a proper subword $\xi$ of $\mathrm{supp}(g)$ in a same lineage. Let $\xi=\xi_0, \ldots, \xi_n=\xi$ be the geodesic between these two vertices in $T_g(u)$. For every $0 \leq i \leq n-1$, there exists a prefix $\Xi_i$ of $g$ such that $\mathrm{supp}(\Xi_i)=\{ \xi_i \}$ and $\mathrm{bot}(X_i) \cap \mathrm{bot}(g)$ is a connected path labeled by $\xi_{i+1}$. In particular, $\Xi_i= \epsilon(x_i)+ \widetilde{\Xi}_i + \epsilon(y_i)$ for some words $x_i,y_i \in \Sigma^+$ and some $(\xi_i,p_i \xi_{i+1}q_i)$-diagram $\widetilde{\Xi}_i$. Now, define
\begin{center}
$\Xi_{i+1}'= \epsilon(x_0p_0 \cdots p_i) + \widetilde{\Xi}_{i+1} + \epsilon(q_i \cdots q_0y_0)$ for $i \geq 0$.
\end{center}
Notice that $\Xi_{i+1}'$ is a $(x_0p_0 \cdots p_i \xi_i q_i \cdots q_0y_0,x_0p_0 \cdots p_{i+1} \xi_{i+1}q_{i+1} \cdots q_0y_0)$-diagram, so that the concatenation $\Xi= \Xi_0 \circ \Xi_1' \circ \cdots \circ \Xi_{n-1}'$ is well-defined. 

\medskip \noindent
We have proved that $g^n$ contains a prefix $\Xi$ such that $\mathrm{supp}(\Xi)=\{ \xi \}$ and $\mathrm{bot}(\Xi) \cap \mathrm{bot}(g^n)$ is the same subpath $\xi$. 

\medskip \noindent
Say that $\Xi$ is a $(x \xi y, xp \xi qy)$-diagram, for some words $x,y,p,q \in \Sigma^+$. In particular, $\Xi$ decomposes as a sum $\epsilon(x)+ \widetilde{\Xi}+ \epsilon(y)$ where $\widetilde{\Xi}$ is a $(\xi,p \xi q)$-diagram. For every $i \geq 0$, set
\begin{center}
$\Delta_i = \epsilon(xp^i) + \widetilde{\Xi} + \epsilon(q^iy)$.
\end{center}
Finally, the concatenation $\Delta= \Delta_0 \circ \Delta_1 \circ \cdots$ defines an infinite prefix of $g^{\infty}$. Moreover,
\begin{center}
$\mathrm{supp}(\Delta)= \mathrm{supp}(\Xi)= \{ \xi \} \subsetneq \mathrm{supp}(g) \subset \mathrm{supp}(g^{\infty})$,
\end{center}
so $\Delta$ is a proper prefix of $g^{\infty}$. $\square$

\medskip \noindent
By combining Lemma \ref{lem1} and Lemma \ref{lem2} with Theorem \ref{prop1}, we obtain the following new criterion.

\begin{prop}\label{criterion2}
Let $g \in D(\mathcal{P},w) - \{ \epsilon(w) \}$ be normal and absolutely reduced. Then $g$ is a contracting isometry of $X(\mathcal{P},w)$ if and only if the following two conditions are satisfied:
\begin{itemize}
	\item the support of $g$ is reduced to $\{ w_1 \}$ for some subword $w_1$ of $w$, and, if we write $w=xw_1y$, then $\partial(\mathcal{P},x)$ and $\partial(\mathcal{P},y)$ are empty;
	\item for every subpath $u$ of $w$, $T_g(u)$ is not deep.
	%contains twice a proper subpath of $\mathrm{supp}(g)$ in a same lineage.
\end{itemize}
\end{prop}

\noindent
\textbf{Proof.} According to Theorem \ref{prop1}, $g$ is contracting if and only if
\begin{itemize}
	\item[(i)] $g^{\infty}$ does not contain any infinite proper prefix;
	\item[(ii)] if $\Delta$ is a reduced diagram with $g^{\infty}$ as a prefix, then $\Delta$ is commensurable to $g^{\infty}$. 
\end{itemize}
Now, Lemma \ref{lem2} states that the condition $(i)$ is equivalent to
\begin{itemize}
	\item[(iii)] for every subpath $u$ of $w$, $T_g(u)$ is not deep;
	%contains twice a proper subpath of $\mathrm{supp}(g)$ in a same lineage.
\end{itemize}
and Lemma \ref{lem1} implies that the condition $(ii)$ is a consequence of
\begin{itemize}
	\item[(iv)] the support of $g$ is reduced to $\{ w_1 \}$ for some subword $w_1$ of $w$, and, if we write $w=xw_1y$, then $\partial(\mathcal{P},x)$ and $\partial(\mathcal{P},y)$ are empty.
\end{itemize}
To conclude the proof, it is sufficient to show that $(iv)$ is a consequence of $(i)$ and $(ii)$. So suppose that $(i)$ and $(ii)$ hold. If $\mathrm{supp}(g)$ is not reduced to a single subpath, then $g$ decomposes as a sum of non-trivial diagrams, and because $g$ is normal, it decomposes as a sum of non-trivial spherical diagrams, say $g=g_1+g_2$. Now, since $g_1$ and $g_2$ are spherical, $g^{\infty} = g_1^{\infty}+ g_2^{\infty}$, so $g^{\infty}$ clearly contains an infinite proper prefix, contradicting $(i)$. Therefore, $\mathrm{supp}(g)$ reduces to a single subpath, and Lemma \ref{lem1} implies that $(iv)$ holds. $\square$ 

\begin{remark}
The first condition of Proposition \ref{criterion2} is obviously satisfied if $\mathrm{supp}(g)= \{ w \}$. Because the case often happen, we will say that spherical diagrams satisfying this property are \emph{full}.
\end{remark}

\begin{ex}
The group $\mathbb{Z} \bullet \mathbb{Z}= \langle a,b,t \mid [a,b^{t^n}]=1, n \geq 0 \rangle$, introduced in \cite[Section 8]{MR1396957}, is isomorphic to the diagram group $D(\mathcal{P},ab)$, where
\begin{center}
$\mathcal{P}= \left\langle a_1,a_2,a_3, p, b_1,b_2,b_3 \left| \begin{array}{l} a_1=a_1p \\ b_1=pb_1 \end{array}, \begin{array}{l} a_1=a_2,a_2=a_3,a_3=a_1 \\ b_1=b_2, b_2=b_3, b_3=b_1 \end{array} \right. \right\rangle$.
\end{center}
Let $g \in D(\mathcal{P},ab)$ be the following spherical diagram group:
\begin{center}
\includegraphics[scale=0.6]{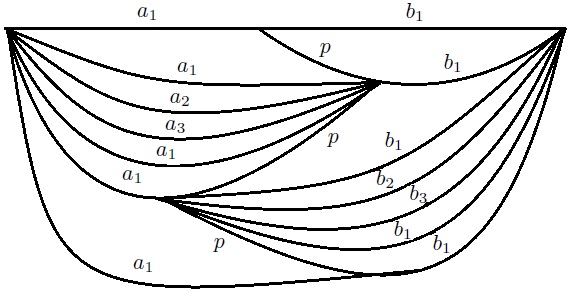}
\end{center}
This is a full, normal, and absolutely reduced diagram. Moreover, $a_1b_1$ is the only admissible subpath, and $T_g(a_1b_1)$ is
\begin{displaymath}
\xymatrix{ \ar@{-}[d] a_1b_1 \\ b_1 }
\end{displaymath}
Therefore, $g$ is a contracting isometry of $X(\mathcal{P},ab)$. In particular, $\mathbb{Z} \bullet \mathbb{Z}$ is acylindrically hyperbolic.
\end{ex}

\begin{ex}
Let $\mathcal{P}= \langle a,b \mid ab=ba, ab^2=b^2a \rangle$ and let $g \in D(\mathcal{P},ab^4)$ be the following spherical diagram:
\begin{center}
\includegraphics[scale=0.6]{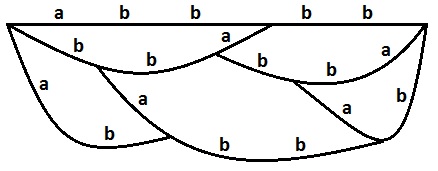}
\end{center}
This is a full, normal, and absolutely reduced diagram. Morever, $ab^4$ is the only admissible subpath, and $T_g(ab^4)$ is
\begin{displaymath}
\xymatrix{& \ar@{-}[dl] \ar@{-}[dr] ab^4 & \\  b^4 \ar@{-}[d] & & b \ar@{-}[d] \\ \emptyset & & \emptyset}
\end{displaymath}
We conclude that $g$ is a contracting isometry of $X(\mathcal{P},ab^4)$. 
\end{ex}

\subsection{Some examples}\label{section:examples}

\noindent
In this section, we exhibit some interesting classes of acylindrically hyperbolic diagram groups by applying the results proved in the previous section. In Example \ref{ex:braid}, we consider a family of groups $U(m_1, \ldots, m_n)$ already studied in \cite{MR1396957} and prove that they are acylindrically hyperbolic. Morever, we notice that each $L_n=U(1, \ldots, 1)$ turns out to be naturally a subgroup of a right-angled Coxeter group. In Example \ref{ex:bullet}, we show that the $\bullet$-product, as defined in \cite{MR1396957}, of two non trivial diagram groups is acylindrically hyperbolic but not relatively hyperbolic. Finally, in Example \ref{nonRAAG}, we prove that some diagram products, as defined in \cite{MR1725439}, of diagram groups turn out to be acylindrically hyperbolic. As a by-product, we also exhibit a cocompact diagram group which is not isomorphic to a right-angled Artin group.

\begin{ex}\label{ex:braid}
Let $\mathcal{P}_n = \langle x_1, \ldots, x_n \mid x_ix_j=x_jx_i, 1 \leq i < j \leq n \rangle$ and let $U(m_1, \ldots, m_n)$ denote the diagram group $D(\mathcal{P},x_1^{m_1} \cdots x_n^{m_n})$, where $m_i \geq 1$ for every $1 \leq i \leq n$. Notice that, for every permutation $\sigma \in S_n$, $U(m_{\sigma(1)},\ldots, m_{\sigma(n)})$ is isomorphic to $U(m_1, \ldots, m_n)$; therefore, we will always suppose that $m_1 \leq \cdots \leq m_n$. For example, the groups $U(p,q,r)$ were considered in \cite[Example 10.2]{MR1396957} (see also \cite[Example 5.3]{arXiv:1505.02053} and \cite[Example 5.11]{arXiv:1507.01667}). In particular, we know that 
\begin{itemize}
	\item $U(m_1, \ldots, m_n)$ is trivial if and only if $n \leq 2$; 
	\item $U(m_1, \ldots, m_n)$ is infinite cyclic if and only if $n=3$ and $m_1=m_2=m_3=1$;
	\item $U(m_1,\ldots, m_n)$ is a non-abelian free group if and only if $n=3$ and $m_1=1$, or $n=4$ and $m_1=m_2=m_3=1$, or $n=5$ and $m_1=m_2=m_3=m_4=m_5=1$.
\end{itemize}
We claim that, whenever it is not cyclic, the group $U(m_1, \ldots, m_n)$ is acylindrically hyperbolic. For instance, if $g \in U(1,1,2)$ denotes the spherical diagram associated to the following derivation
\begin{center}
$x_1x_2x_3x_3 \to x_2x_1x_3x_3 \to x_2x_3x_1x_3 \to x_3x_2x_1x_3 \to x_3x_1x_2x_3 \to x_1x_3x_2x_3 \to x_1x_2x_3x_3$,
\end{center}
then $g$ is a contracting isometry of the associated Farley cube complex. (The example generalizes easily to other cases although writing an argument in full generality is laborious.) In particular, we deduce that $U(m_1,\ldots, m_n)$ does not split as non trivial direct product.

\medskip \noindent
An interesting particular case is the group $L_n$, corresponding to $U(1, \ldots, 1)$ with $n$ ones. A first observation is that $L_n$ is naturally a subgroup of $U(m_1, \ldots, m_n)$, and conversely, $U(m_1, \ldots, m_n)$ is naturally a subgroup of $L_{m_1+ \cdots + m_n}$. Secondly, $L_n$ can be described as a group of \emph{pseudo-braids}; this description is made explicit by the following example of a diagram with its associated pseudo-braid:
\begin{center}
\includegraphics[scale=0.6]{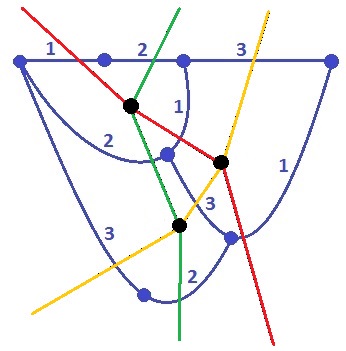}
\end{center}
Of course, we look at our pseudo-braids up to the following move, corresponding to the reduction of a dipole in the associated diagram:
\begin{center}
\includegraphics[scale=0.6]{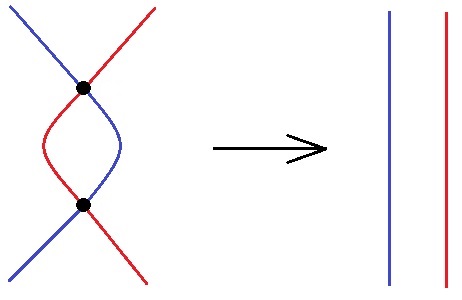}
\end{center}
Finally, if $C_n$ denotes the group of all pictures of pseudo-braids, endowed with the obvious concatenation, then $L_n$ corresponds to ``pure subgroup'' of $C_n$. It is clear that, if $\sigma_i$ corresponds to the element of $C_n$ switching the $i$-th and $(i+1)$-th braids, then a presentation of $C_n$ is 
\begin{center}
$\langle \sigma_1, \ldots, \sigma_{n-1} \mid \sigma_i^2=1, [\sigma_i,\sigma_j]=1 \ \text{if} \ |i-j| \geq 2 \rangle$.
\end{center}
Alternatively, $C_n$ is isomorphic to the right-angled Coxeter group $C( \overline{P_{n-1}})$, where $\overline{P_{n-1}}$ is the complement of the graph $P_{n-1}$ which is a segment with $n-1$ vertices, i.e., $\overline{P_{n-1}}$ is the graph with the same vertices as $P_{n-1}$ but whose edges link two vertices precisely when they are not adjacent in $P_{n-1}$. In particular, $L_n$ is naturally a finite-index subgroup of $C( \overline{P_{n-1}})$.

\medskip \noindent
With this interpretation, the group $U(m_1, \ldots, m_n)$ corresponds to the subgroup of $L_{m}$, where $m=m_1+ \cdots +m_n$, defined as follows: Color the $m$ braids with $n$ colors, such that, for every $1 \leq i \leq n$, there are $m_i$ braids with the color $i$. Now, $U(m_1, \ldots, m_n)$ is the subgroup of the pure pseudo-braids where two braids of the same color are never switched.

\medskip \noindent
Therefore, the groups $U(m_1, \ldots, m_n)$ turn out to be natural subgroups of right-angled Coxeter groups.
\end{ex}

\begin{ex}\label{ex:bullet}
Let $G$ and $H$ be two groups. We define the product $G \bullet H$ by the relative presentation
\begin{center}
$\langle G,H,t \mid [g,h^t]=1, g \in G,h \in H \rangle$.
\end{center}
As proved in \cite[Theorem 8.6]{MR1396957}, the class of diagram groups is stable under the $\bullet$-product. More precisely, let $\mathcal{P}_1= \langle \Sigma_1 \mid \mathcal{R}_1 \rangle$, $\mathcal{P}_2= \langle \Sigma_2 \mid \mathcal{R}_2 \rangle$ be two semigroup presentations and $w_1,w_2$ two base words. For $i=1,2$, suppose that there does not exist two words $x$ and $y$ such that $xy$ is non-empty and $w_i=xw_iy$; as explained in during the proof of \cite[Theorem 8.6]{MR1396957}, we may suppose that this condition holds up to a small modification of $\mathcal{P}_i$ which does not modify the diagram group. Then, if
\begin{center}
$\mathcal{P}= \langle \Sigma_1 \sqcup \Sigma_2 \sqcup \{p \} \mid \mathcal{R}_1 \sqcup \mathcal{R}_2 \sqcup \{w_1=w_1p, w_2=pw_2 \} \rangle$,
\end{center}
then $D(\mathcal{P},w_1w_2)$ is isomorphic to $D(\mathcal{P}_1,w_1) \bullet D(\mathcal{P}_2,w_2)$ \cite[Lemma 8.5]{MR1396957}. Now, we claim that, whenever $D(\mathcal{P}_1,w_1)$ and $D(\mathcal{P}_2,w_2)$ are non-trivial, the product $D(\mathcal{P}_1,w_1) \bullet D(\mathcal{P}_2,w_2)$ is acylindrically hyperbolic. Indeed, if $\Gamma \in D(\mathcal{P}_1,w_1)$ and $\Delta \in D(\mathcal{P}_2,w_2)$ are non-trivial, then the spherical diagram
\begin{center}
\includegraphics[scale=0.6]{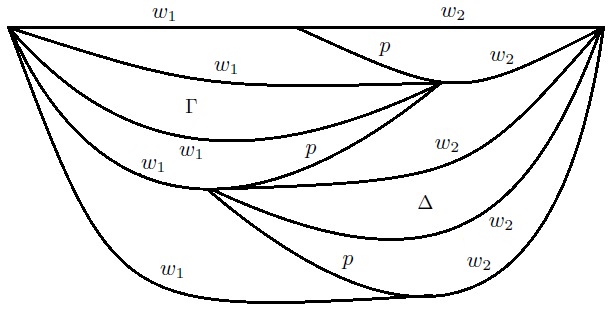}
\end{center}
of $D(\mathcal{P},w_1w_2)$ is a contracting isometry. On the other hand, $D(\mathcal{P}_1,w_1) \bullet D(\mathcal{P}_2,w_2)$ contains a copy of $\mathbb{Z} \bullet \mathbb{Z}$, so it cannot be cyclic. We conclude that $D(\mathcal{P}_1,w_1) \bullet D(\mathcal{P}_2,w_2)$ is indeed acylindrically hyperbolic.

\medskip \noindent
In fact, the product $G \bullet H$ can be described as an HNN extension: if $H^{\infty}$ denote the free product of infinitely many copies of $H$, then $G \bullet H$ is isomorphic to the HNN extension of $G \times H^{\infty}$ over $H^{\infty}$ with respect to two monomorphisms $H^{\infty} \hookrightarrow G \times H^{\infty}$ associated to $H_i \mapsto H_i$ and $H_i \mapsto H_{i+1}$. Now, if $t$ denotes the stable letter of the HNN extension and if $h \in H^{\infty}$ is a non-trivial element of the first copy of $H$ in $H^{\infty}$, then it follows from Britton's lemma that $H^{\infty} \cap (H^{\infty})^{t^{-1}ht} = \{1 \}$. Therefore, $H^{\infty}$ is a weakly malnormal subgroup, and it follows from \cite[Corollary 2.3]{arXiv:1310.6289} that, if $G$ and $H$ are both non-trivial, the product $G \bullet H$ is acylindrically hyperbolic. It is interesting that:

\begin{fact}
If $G$ and $H$ are torsion-free, then any proper malnormal subgroup of $G \bullet H$ is a free group.
\end{fact}

\noindent
Indeed, let $K$ be a malnormal subgroup of $G \bullet H$. Suppose that $K$ contains a non-trivial element of $G \times H^{\infty}$, say $gh$. Then $h \in K$ and $g \in K$ since 
\begin{center}
$\langle gh \rangle = \langle gh \rangle \cap \langle gh \rangle^h \subset K \cap K^h$ and $\langle gh \rangle = \langle gh \rangle \cap \langle gh \rangle^g \subset K \cap K^g$.
\end{center}
If $h$ is not trivial, then $G \subset K$ since
\begin{center}
$\langle h \rangle = \langle h \rangle \cap \langle h \rangle^G \subset K \cap K^G$;
\end{center}
then, we deduce that $H^{\infty} \subset K$ since
\begin{center}
$G=G \cap G^{H^{\infty}} \subset K \cap K^{H^{\infty}}$.
\end{center}
Therefore, $G \times H^{\infty} \subset K$. Similarly, if $h$ is trivial then $g$ is non-trivial, and we deduce that $H^{\infty} \subset K$ from
\begin{center}
$\langle g \rangle = \langle g \rangle \cap \langle g \rangle^{H^{\infty}} \subset K \cap K^{H^{\infty}}$.
\end{center}
We also notice that $G \subset K$ since
\begin{center}
$H^{\infty} = H^{\infty} \cap (H^{\infty})^{G} \subset K \cap K^G$.
\end{center}
Consequently, $G \times H^{\infty} \subset K$. Thus, we have proved that a malnormal subgroup of $G \bullet H$ either intersects trivially $G \times H^{\infty}$ or contains $G \times H^{\infty}$. If $K$ intersects trivially $G \times H^{\infty}$, the action of $K$ on the Bass-Serre tree of $G \bullet H$, associated to its decomposition as an HNN extension we mentionned above, is free, so $K$ is necessarily free. Otherwise, $K$ contains $G \times H^{\infty}$. If $t$ denotes the stable letter of our HNN extension, then 
\begin{center}
$H_2 \ast H_3 \ast \cdots = H^{\infty} \cap (H^{\infty})^t \subset K \cap K^t$,
\end{center}
where we used the notation $H^{\infty} = H_1 \ast H_2 \ast H_3 \ast \cdots$, such that each $H_i$ is a copy of $H$. Therefore, $t \in K$, and we conclude that $K$ is not a proper subgroup. The fact is proved.

\medskip \noindent
As a consequence, $\bullet$-products of two non-trivial diagram groups are good examples of acylindrically hyperbolic groups in the sense that they are not relatively hyperbolic. Indeed, if a group is hyperbolic relatively to a collection of subgroups, then this collection must be malnormal (see for instance \cite[Theorem 1.4]{OsinRelHyp}), so that a relatively hyperbolic $\bullet$-product of two torsion-free groups must hyperbolic relatively to a collection of free groups according to the previous fact. On the other hand, such a group must be hyperbolic (see for instance \cite[Corollary 2.41]{OsinRelHyp}), which is impossible since we know that a $\bullet$-product of two non trivial torsion-free groups contains a subgroup isomorphic to $\mathbb{Z}^2$.
\end{ex}

\begin{ex}\label{nonRAAG}
Let $\mathcal{P}= \langle \Sigma \mid \mathcal{R} \rangle$ be a semigroup presentation and $w \in \Sigma^+$ a base word. Fixing a family of groups $\mathcal{G}_{\Sigma} = \{ G_s \mid s \in \Sigma \}$, we define the \textit{diagram product} $\mathcal{D}(\mathcal{G}_{\Sigma}, \mathcal{P},w)$ as the fundamental group of the following 2-complex of groups:
\begin{itemize}
	\item the underlying 2-complex is the 2-skeleton of the Squier complex $S(\mathcal{P},w)$;
	\item to any vertex $u=s_1\cdots s_r \in \Sigma^+$ is associated the group $G_u=G_{s_1} \times \cdots \times G_{s_r}$;
	\item to any edge $e=(a, u \to v,b)$ is associated the group $G_{e}=G_a \times G_b$;
	\item to any square is associated the trivial group;
	\item for every edge $e=(a,u \to v,b)$, the monomorphisms $G_e \to G_{aub}$ and $G_e \to G_{avb}$ are the canonical maps $G_a \times G_b \to G_a \times G_u \times G_b$ and $G_a \times G_b \to G_a \times G_v \times G_b$.
\end{itemize}
Guba and Sapir proved that the class of diagram groups is stable under diagram product \cite[Theorem 4]{MR1725439}. Explicitely,

\begin{thm}
If $G_s= D(\mathcal{P}_s,w_s)$ where $\mathcal{P}_s= \langle \Sigma_s \mid \mathcal{R}_s \rangle$ is a semigroup presentation and $w_s \in \Sigma_s^+$ a base word, then the diagram product $\mathcal{D}(\mathcal{G}_{\Sigma},\mathcal{P},w)$ is the diagram group $D(\bar{\mathcal{P}},w)$, where $\bar{\mathcal{P}}$ is the semigroup presentation 
\begin{center}
$\langle \Sigma \sqcup \bar{\Sigma} \sqcup \Sigma_0 \mid \mathcal{R} \sqcup \bar{\mathcal{R}} \sqcup \mathcal{S} \rangle$ 
\end{center}
where $\Sigma_0= \{ a_s \mid s \in \Sigma \}$ is a copy of $\Sigma$, $\bar{\Sigma} = \bigsqcup\limits_{s \in \Sigma} \Sigma_s$, $\bar{\mathcal{R}}= \bigsqcup\limits_{s \in \Sigma} \mathcal{R}_s$ and finally $\mathcal{S}= \bigsqcup\limits_{s \in \Sigma} \{ s=a_sw_sas \}$. \end{thm}

\noindent
Of course, any semigroup diagram with respect to $\mathcal{P}$ is a semigroup diagram with respect to $\bar{\mathcal{P}}$. Thus, if $D(\mathcal{P},w)$ contains a normal, absolutely reduced, spherical diagram which does not contains a proper infinite prefix and whose support is $\{w\}$, then any diagram product over $(\mathcal{P},w)$ will be acylindrically hyperbolic or cyclic. (Notice that such a diagram product contains $D(\mathcal{P},w)$ as a subgroup, so it cannot be cyclic if $D(\mathcal{P},w)$ is not cyclic itself.)

\medskip \noindent
For instance, let $\mathcal{P}= \langle x,y,z,a,b \mid yz=az, xa=xb, b=y \rangle$. Then $D(\mathcal{P},xyz)$ contains the following spherical diagram, which is a contracting isometry of $X(\mathcal{P},xyz)$:
\begin{center}
\includegraphics[scale=0.6]{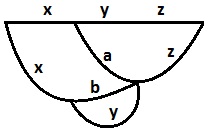}
\end{center}
By the above observation, if $G$ denotes the diagram product over $(\mathcal{P},apb)$ with $G_x=G_z=\mathbb{Z}$ and $G_y=G_a=G_b=\{1 \}$, then $G$ is acylindrically hyperbolic. In this case, the Squier complex $S(\mathcal{P},xyz)$ is just a cycle of length three. By computing a presentation of the associated graph of groups, we find
\begin{center}
$G = \langle x,y,t \mid [x,y]=[x,y^t]=1 \rangle$.
\end{center}
This example is interesting because this is a cocompact diagram group which is not a right-angled Artin group. Let us prove this assertion. From the presentation of $G$, it is clear that its abelianization is $\mathbb{Z}^3$ and we deduce from \cite[Exercice II.5.5]{BrownCoho} that the rank of $H_2(G)$ is at most two (the number of relations of our presentation). On the other hand, if $A(\Gamma)$ denotes the right-angled Artin group associated to a given graph $\Gamma$, then $H_1(A(\Gamma))= \mathbb{Z}^{\# V(\Gamma)}$ and $H_2(A(\Gamma))= \mathbb{Z}^{\# E(\Gamma)}$, where $V(\Gamma)$ and $E(\Gamma)$ denote the set of vertices and edges of $\Gamma$ respectively. Therefore, the only right-angled Artin groups which might be isomorphic to $G$ must correspond to a graph with three vertices and at most two edges. We distinguish two cases: either this graph is not connected or it is a segment of length two. In the former case, the right-angled Artin group decomposes as a free product, and in the latter, it is isomorphic to $\mathbb{F}_2 \times \mathbb{Z}$. Thus, it is sufficient to show that $G$ is freely irreducible and is not isomorphic to $\mathbb{F}_2 \times \mathbb{Z}$ in order to deduce that $G$ cannot be isomorphic to a right-angled Artin group. First, notice that, because $G$ is acylindrically hyperbolic, its center must be finite (see \cite[Corollary 7.3]{arXiv:1304.1246}), so $G$ cannot be isomorphic to $\mathbb{F}_2 \times \mathbb{Z}$. Next, suppose that $G$ decomposes as a free product. Because the centralizer of $x$ is not virtually cyclic, we deduce that $x$ belongs to a conjugate $A$ of some free factor; in fact, $A$ must contain the whole centralizer of $x$, hence $y,y^t \in A$. As a consequence, $y^t \in A \cap A^t$, so that $A \cap A^t$ must be infinite. Since a free factor is a malnormal subgroup, we deduce that $t \in A$. Finally, $x,y,t \in A$ hence $A=G$. Therefore, $G$ does not contain any proper free factor, i.e., $G$ is freely irreducible. This concludes the proof of our assertion.

\medskip \noindent
To conclude, it is worth noticing that the whole diagram product may contain contracting isometries even if the underlying diagram group (i.e., the diagram group obtained from the semigroup presentation and the word base, forgetting the factor-groups) does not contain such isometries on its own Farley cube complex. For instance, if $\mathcal{P}= \langle a,b,p \mid a=ap,b=pb \rangle$, $w=ab$, $G_p=\{1 \}$ and $G_a=D(\mathcal{P}_1,w_1)$, $G_b= D(\mathcal{P}_2,w_2)$ are non-trivial, then our diagram product is isomorphic to $G_a \bullet G_b$ (see \cite[Example 7]{MR1725439}), and the description of $G_a \bullet G_b$ as a diagram group is precisely the one we mentionned in the previous example, where we saw that the action on the associated Farley cube complex admits contracting isometries. On the other hand, $X(\mathcal{P},ab)$ is a ``diagonal half-plane'', so that its combinatorial boundary does not contain any isolated point: a fortiori, $X(\mathcal{P},ab)$ cannot admit a contracting isometry.
\end{ex}

\subsection{Cocompact case}

\noindent
In this section, we focus on cocompact diagram groups, ie., we will consider a semigroup presentation $\mathcal{P}= \langle \Sigma \mid \mathcal{R} \rangle$ and a word base $w \in \Sigma^+$ whose class modulo $\mathcal{P}$ is finite. Roughly speaking, we want to prove that $D(\mathcal{P},w)$ contains a contracting isometry if and only if it does not split as a direct product in a natural way, which we describe now.

\medskip \noindent
If there exists a $(w,u_1v_1\cdots u_nv_nu_{n+1})$-diagram $\Gamma$, for some words $u_1,v_1, \ldots, u_n,v_n,u_{n+1} \in \Sigma^+$, then we have the natural map $$\left\{ \begin{array}{ccc} D(\mathcal{P},v_1) \times \cdots \times D(\mathcal{P},v_n) & \to & D(\mathcal{P},w) \\ (V_1, \ldots, V_n) & \mapsto & \Gamma \cdot (\epsilon(u_1)+ V_1 + \cdots + \epsilon(u_n)+V_n+\epsilon(u_{n+1})) \cdot \Gamma^{-1} \end{array} \right.$$ This is a monomorphism, and we denote by $\Gamma \cdot D(\mathcal{P},v_1) \times \dots \times D(\mathcal{P},v_n) \cdot \Gamma^{-1}$ its image. A subgroup of this form will be referred to as a \emph{canonical subgroup}; if furthermore at least two factors are non trivial, then it will be a \emph{large canonical subgroup}. 

\medskip \noindent
The main result of this section is the following decomposition theorem (cf. Theorem \ref{main7} in the introduction).

\begin{thm}\label{decompositioncor}
Let $\mathcal{P}= \langle \Sigma \mid \mathcal{R} \rangle$ be a semigroup presentation and $w \in \Sigma^+$ a base word whose class modulo $\mathcal{P}$ is finite. Then there exist some words $u_0,u_1, \ldots, u_m \in \Sigma^+$ and a $(w,u_0u_1 \cdots u_m)$-diagram $\Gamma$ such that
\begin{center}
$D(\mathcal{P},w)= \Gamma \cdot ( D(\mathcal{P},u_0) \times D(\mathcal{P},u_1) \times \cdots \times D(\mathcal{P},u_m)) \cdot \Gamma^{-1}$,
\end{center}
where $D(\mathcal{P},u_i)$ is either infinite cyclic or acylindrically hyperbolic for every $1 \leq i \leq m$, and $m \leq \dim_{\mathrm{alg}} D(\mathcal{P},w) \leq \dim X(\mathcal{P},w)$.
\end{thm}

\noindent
The acylindrical hyperbolicity will follow from the existence of a contracting isometry in the corresponding Farley cube complex. In order to understand what happens when such an isometry does not exist, we will use the following dichotomy, which is a slight modification, in the cocompact case, of the Rank Rigidity Theorem proved in \cite{MR2827012}. 

\begin{thm}\label{thm:rankrigidity}
Let $G$ be a group acting geometrically on an unbounded CAT(0) cube complex $X$. Either $G$ contains a contracting isometry of $X$ or there exists a $G$-invariant convex subcomplex $Y \subset X$ which splits as the cartesian product of two unbounded subcomplexes.
\end{thm}

\noindent
\textbf{Proof.} Recall that a hyperplane of $X$ is \emph{essential} if no $G$-orbit lies in a neighborhood of some halfspace delimited by this hyperplane; we will denote by $\mathcal{E}(X)$ the set of the essential hyperplanes of $X$. Let $Y \subset X$ denote the \emph{essential core} associated to the action $G \curvearrowright X$, which is defined as the \emph{restriction quotient} of $X$ with respect to $\mathcal{E}(X)$. See \cite[Section 3.3]{MR2827012} for more details. According to \cite[Proposition 3.5]{MR2827012}, $Y$ is a $G$-invariant convex subcomplex of $X$. 

\medskip \noindent
It is worth noticing that, because the action $G \curvearrowright X$ is cocompact, a hyperplane $J$ of $X$ does not belong to $\mathcal{E}(X)$ if and only if one of the halfspaces it delimits lies in the $R(J)$-neighborhood of $N(J)$ for some $R(J) \geq 1$. Set $R= \max\limits_{J \in \mathcal{H}(X)} R(J)$. Because there exist only finitely many orbits of hyperplanes, $R$ is finite.

\medskip \noindent
We claim that two hyperplanes of $Y$ are well-separated in $Y$ if and only if they are well-separated in $X$. 

\medskip \noindent
Of course, two hyperplanes of $Y$ which are well-separated in $X$ are well-separated in $Y$. Conversely, let $J_1,J_2 \in \mathcal{H}(Y)$ be two hyperplanes and fix some finite collection $\mathcal{H} \subset \mathcal{H}(X)$ of hyperplanes intersecting both $J_1,J_2$ which does not contain any facing triple. We write $\mathcal{H}= \{ H_0, \ldots, H_n \}$ and we fix a halfspace $H_i^+$ delimited by $H_i$ for every $0 \leq i \leq n$, so that $H_0^+ \supset H_1^+ \supset \cdots \supset H_n^+$. If $n \leq 2R$, there is nothing to prove. Otherwise, by our choice of $R$, it is clear that the hyperplanes $H_R, \ldots, H_{n-R}$ belong to $\mathcal{E}(X)$, and a fortiori to $\mathcal{H}(Y)$. Consequently, if $J_1$ and $J_2$ are $L$-well-separated in $Y$, then they are $(L+2R)$-well-separated in $X$. 

\medskip \noindent
Now we are ready to prove our theorem. Because the action $G \curvearrowright Y$ is geometric and \emph{essential} (ie., every hyperplane of $Y$ is essential), it follows from \cite[Theorem 6.3]{MR2827012} that either $Y$ splits non trivially as a cartesian product of two subcomplexes or $G$ contains a contracting isometry of $Y$. In the former case, notice that the factors cannot be bounded because the action $G \curvearrowright Y$ is essential; in the latter case, we deduce that the contracting isometry of $Y$ induces a contracting isometry of $X$ by combining our previous claim with Theorem \ref{equivalence - contracting isometry}. This concludes the proof. $\square$

\medskip \noindent
Therefore, the problem reduces to understand the subproducts of $X(\mathcal{P},w)$. This purpose is achieved by the next proposition.

\begin{prop}\label{prop:subproducts}
Let $\mathcal{P}= \langle \Sigma \mid \mathcal{R} \rangle$ be a semigroup presentation and $w \in \Sigma^+$ a base word. Let $Y \subset X(\mathcal{P},w)$ be a convex subcomplex which splits as the cartesian product of two subcomplexes $A,B \subset Y$. Then there exist a word $w' \in \Sigma^+$, a $(w,w')$-diagram $\Gamma$ and words $u_1,v_1, \ldots, u_k,v_k \in \Sigma^+$ for some $k \geq 1$ such that $w'=u_1v_1 \cdots u_kv_k$ in $\Sigma^+$ and $$A \subset \Gamma \cdot X(\mathcal{P},u_1) \times \cdots \times X(\mathcal{P},u_k) \cdot \Gamma^{-1} \ \text{and} \ B \subset \Gamma \cdot X(\mathcal{P},v_1) \times \cdots \times X(\mathcal{P},v_k) \cdot \Gamma^{-1}.$$ Furthermore, $$Y \subset \Gamma \cdot X(\mathcal{P},u_1) \times X( \mathcal{P},v_1) \times \cdots \times X(\mathcal{P},u_k) \times X ( \mathcal{P},v_k) \cdot \Gamma^{-1}.$$
\end{prop}

\noindent
In the previous statement, $Y$ is a Cartesian product of the subcomplexes $A$ and $B$ in the sense that there exists an isomorphism $\varphi : A \times B \to Y$ which commutes with the inclusions $A,B,Y \hookrightarrow X$. More explicitely, we have the following commutative diagram:
\begin{displaymath}
\xymatrix{ & A \ar@{_{(}->}[dl] \ar@{^{(}->}[dr] & \\ A \times B \ar[r]^{\varphi} & Y \ar@{^{(}->}[r] & X \\ & B \ar@{_{(}->}[ul] \ar@{^{(}->}[ur] & }
\end{displaymath}

\noindent
\textbf{Proof of Proposition \ref{prop:subproducts}.} Up to conjugating by a diagram $\Gamma$ of $A \cap B$, we may assume without loss of generality that $\epsilon(w) \in A \cap B$. First, we claim that, for every $\Delta_1 \in A$ and $\Delta_2 \in B$, $\mathrm{supp}(\Delta_1) \cap \mathrm{supp}(\Delta_2)= \emptyset$. 

\medskip \noindent
Suppose by contradiction that there exist some $\Delta_1 \in A$ and $\Delta_2 \in B$ satisfying $\mathrm{supp}(\Delta_1) \cap \mathrm{supp}(\Delta_2) \neq \emptyset$. Without loss of generality, we may suppose that we chose a counterexample minimizing $\# \Delta_1+ \# \Delta_2$. Because $\mathrm{supp}(\Delta_1) \cap \mathrm{supp}(\Delta_2) \neq \emptyset$, there exist two cells $\pi_1 \subset \Delta_1$, $\pi_2 \subset \Delta_2$ and an edge $e \subset w$, where $w$ is thought of as a segment representing both $\mathrm{top}(\Delta_1)$ and $\mathrm{top}(\Delta_2)$, such that $e \subset \mathrm{top}(\pi_1) \cap \mathrm{top}(\pi_2)$. Notice that $\pi_i$ is the only cell of the maximal thin suffix of $\Delta_i$ for $i=1,2$. Indeed, if this was not the case, taking the minimal prefixes $P_1,P_2$ of $\Delta_1, \Delta_2$ containing $\pi_1, \pi_2$ respectively, we would produce two diagrams satisfying $\mathrm{supp}(P_1) \cap \mathrm{supp}(P_2) \neq \emptyset$ and $\# P_1 + \# P_2 < \# \Delta_1 + \# \Delta_2$. Moreover, because $Y$ is convex necessarily the factors $A,B$ have to be convex as well, so that $A$ and $B$, as sets of diagrams, have to be closed by taking prefixes since $\epsilon(w) \in A \cap B$ and according to the description of the combinatorial geodesics given by Lemma \ref{geodesic}; this implies that $P_1 \in A$ and $P_2 \in B$, contradicting the choice of our counterexample $\Delta_1, \Delta_2$. For $i=1,2$, let $\overline{\Delta}_i$ denote the diagram obtained from $\Delta_i$ by removing the cell $\pi_i$. We claim that $$\mathrm{supp}(\overline{\Delta}_1) \cap \mathrm{suppp}( \Delta_2)=\mathrm{supp}(\Delta_1) \cap \mathrm{suppp}( \overline{\Delta}_2)= \emptyset.$$
Indeed, since we know that $A$ and $B$ are stable by taking prefixes, $\overline{\Delta}_1$ belongs to $A$ and $\overline{\Delta}_2$ belongs to $B$. Therefore, because 
$$\# \overline{\Delta}_1 + \# \Delta_2 = \# \Delta_1 + \# \overline{\Delta}_2 < \# \Delta_1+ \# \Delta_2,$$
having $\mathrm{supp}(\overline{\Delta}_1) \cap \mathrm{suppp}( \Delta_2) \neq \emptyset$ or $\mathrm{supp}(\Delta_1) \cap \mathrm{suppp}( \overline{\Delta}_2) \neq \emptyset$ would contradict our choice of $\Delta_1$ and $\Delta_2$. 

%Otherwise, there would exist two cells $\pi_1' \subset \overline{\Delta}_1$, $\pi_2' \subset \overline{\Delta}_2$ and an edge $e' \subset \mathrm{supp}(\overline{\Delta}_1) \cap \mathrm{supp}( \overline{\Delta}_2)$ such that $e' \subset \mathrm{top}(\pi_1') \cap \mathrm{top}(\pi_2')$. Now, since we know that $A$ and $B$ are stable by taking prefixes, replacing $\Delta_1,\Delta_2$ with their minimal prefixes containing $\pi_1',\pi_2'$ respectively would produce a smaller counterexample, contradiction our choice of $\Delta_1, \Delta_2$. 

\medskip \noindent
Because $\mathrm{supp}(\overline{\Delta}_1) \cap \mathrm{suppp}( \overline{\Delta}_2)= \emptyset$, it is possible to write $w=a_1b_1 \cdots a_pb_p$ in $\Sigma^+$ so that $$\left\{ \begin{array}{l} \overline{\Delta}_1 = A_1+ \epsilon(b_1) + \cdots + A_p+ \epsilon(b_p) \\ \overline{\Delta}_2 = \epsilon(a_1) + B_1 +  \cdots + \epsilon(a_p) + B_p \end{array} \right.$$ for some $(a_i,\ast)$-diagram $A_i$ and $(b_i,\ast)$-diagram $B_i$. Set $\Delta_0= A_1+B_1+ \cdots +A_p+B_p$. Notice that $\overline{\Delta}_1$ and $\overline{\Delta}_2$ are prefixes of $\Delta_0$. 

\medskip \noindent
For $i=1,2$, let $E_i$ denote the atomic diagram such that the concatenation $\Delta_0 \circ E_i$ corresponds to gluing $\pi_i$ on $\overline{\Delta}_i$ as a prefix of $\Delta_0$. Notice that the diagrams $E_1,E_2$ exist precisely because
$$\mathrm{supp}(\overline{\Delta}_1) \cap \mathrm{suppp}( \Delta_2)=\mathrm{supp}(\Delta_1) \cap \mathrm{suppp}( \overline{\Delta}_2)= \emptyset.$$
As $\mathrm{top}(\pi_1) \cap \mathrm{top}(\pi_2) \neq \emptyset$, since the intersection contains $e$, the two edges $(\Delta_0, \Delta_0 \circ E_1)$ and $(\Delta_0, \Delta_0 \circ E_2)$ of $X(\mathcal{P},w)$ do not span a square. Therefore, the hyperplanes $J_1,J_2$ dual to these two edges respectively are tangent. For $i=1,2$, noticing that $\Delta_i$ is a prefix of $\Delta_0 \circ E_i$ whereas it is not a prefix of $\Delta_0$, we deduce from Proposition \ref{prop:hypFarley} that the minimal diagram associated to $J_i$ is precisely $\Delta_i$. On the other hand, the hyperplanes corresponding to $\Delta_1, \Delta_2$ are clearly dual edges of $A,B$ respectively. A fortiori, $J_1$ and $J_2$ must be transverse (in $Y= A \times B$). Therefore, we have found an inter-osculation in the CAT(0) cube complex $X(\mathcal{P},w)$, which is impossible.

\medskip \noindent
Thus, we have proved that, for every $\Delta_1 \in A$ and $\Delta_2 \in B$, $\mathrm{supp}(\Delta_1) \cap \mathrm{supp}(\Delta_2)= \emptyset$. In particular, we can write $w=u_1v_1 \cdots u_kv_k$ in $\Sigma^+$, for some $u_1,v_1, \ldots, u_k,v_k \in \Sigma^+$, so that $\mathrm{supp}(\Delta_1) \subset u_1 \cup \cdots \cup u_k$ and $\mathrm{supp}(\Delta_2) \subset v_1 \cup \cdots \cup v_k$ for every $\Delta_1 \in A$ and $\Delta_2 \in B$. By construction, it is clear that $$A \subset X(\mathcal{P},u_1) \times \cdots \times X(\mathcal{P},u_k) \ \text{and} \ B \subset X(\mathcal{P},v_1) \times \cdots \times X(\mathcal{P},v_k).$$ Now, let $\Delta \in Y=A \times B$ be vertex and let $\Delta_1 \in A$, $\Delta_2 \in B$ denote its coordinates; they are also the projections of $\Delta$ onto $A,B$ respectively. By the previous remark, we can write $$\left\{ \begin{array}{l} \Delta_1= U_1 + \epsilon(v_1)+ \cdots + U_k+ \epsilon(v_k) \\ \Delta_2= \epsilon(u_1)+ V_1+ \cdots + \epsilon(u_k)+ V_k \end{array} \right.$$ for some $(u_i, \ast)$-diagram $U_i$ and some $(v_i, \ast)$-diagram $V_i$. Set $\Xi = U_1+V_1+ \cdots U_k+V_k$. Noticing that $d(\Xi,A) \geq \# V_1+ \cdots + \# V_k$ and $d(\Xi,\Delta_1)= \# V_1+ \cdots + \# V_k$, we deduce that $\Delta_1$ is the projection of $\Xi$ onto $A$. Similarly, we show that $\Delta_2$ is the projection of $\Xi$ onto $B$. Consequently, $$\Delta=\Xi \in X(\mathcal{P},u_1) \times X(\mathcal{P},v_1) \times \cdots \times X(\mathcal{P},u_k) \times X(\mathcal{P},v_k),$$ which concludes the proof. $\square$

\medskip \noindent
By combining Theorem \ref{thm:rankrigidity} and Proposition \ref{prop:subproducts}, we are finally able to prove:

\begin{cor}\label{cor:dichotomy}
Let $\mathcal{P}= \langle \Sigma \mid \mathcal{R} \rangle$ be a semigroup presentation and $w \in \Sigma^+$ a base word whose class modulo $\mathcal{P}$ is finite. Then either $D(\mathcal{P},w)$ contains a contracting isometry of $X(\mathcal{P},w)$ or there exists a $(w,uv)$-diagram $\Gamma$ such that $$D(\mathcal{P},w) = \Gamma \cdot D(\mathcal{P},u) \times D(\mathcal{P},v) \cdot \Gamma^{-1},$$ where $D(\mathcal{P},u)$ and $D(\mathcal{P},v)$ are non trivial. 
\end{cor}

\noindent
\textbf{Proof.} According to Theorem \ref{thm:rankrigidity}, if $D(\mathcal{P},w)$ does not contain a contracting isometry of $X(\mathcal{P},w)$, then there exists a convex $D(\mathcal{P},w)$-invariant subcomplex $Y \subset X(\mathcal{P},w)$ which splits as the cartesian product of two unbounded subcomplexes $A,B$. Up to conjugating by the smallest diagram of $Y$, we may assume without loss of generality that $Y$ contains $\epsilon(w)$. Let $\Gamma$ be the $(w,u_1v_1 \cdots u_kv_k)$-diagram given by Proposition \ref{prop:subproducts}. Notice that, because $\epsilon(w) \in Y$ and that $Y$ is $D(\mathcal{P},w)$-invariant, necessarily any spherical diagram belongs to $Y$. Thus, we deduce from Lemma \ref{somme} below that $$D(\mathcal{P},w)= \Gamma \cdot D(\mathcal{P},u_1) \times D(\mathcal{P},v_1) \times \cdots \times D(\mathcal{P},u_k) \times D(\mathcal{P},v_k) \cdot \Gamma^{-1}.$$ Because $A \subset \Gamma \cdot X(\mathcal{P},u_1) \times \cdots \times X(\mathcal{P},u_k) \cdot \Gamma^{-1}$, we deduce from the unboundedness of $A$ that $X(\mathcal{P},u_1) \times \cdots \times X(\mathcal{P},u_k)$ is infinite; since the class of $w$ modulo $\mathcal{P}$ is finite, this implies that there exists some $1 \leq i \leq k$ such that $D(\mathcal{P},u_i)$ is non trivial. Similarly, we show that there exists some $1 \leq j \leq k$ such that $D(\mathcal{P},v_j)$ is non trivial. Now, write $u_1v_1 \cdots u_kv_k=uv$ in $\Sigma^+$ where $u,v \in \Sigma^+$ are two subwords such that one contains $u_i$ and the other $v_j$. Then $$D(\mathcal{P},w) = \Gamma \cdot D(\mathcal{P},u) \times D(\mathcal{P},v) \cdot \Gamma^{-1},$$ where $D(\mathcal{P},u)$ and $D(\mathcal{P},v)$ are non trivial. $\square$

\begin{lemma}\label{somme}
Let $\mathcal{P}= \langle \Sigma \mid \mathcal{R} \rangle$ be a semigroup presentation and $w \in \Sigma^+$ a base word whose class modulo $\mathcal{P}$ is finite. If a spherical diagram $\Delta$ decomposes as a sum $\Delta_1 + \cdots + \Delta_n$, then each $\Delta_i$ is a spherical diagram.
\end{lemma}

\noindent
\textbf{Proof.} Suppose by contradiction that one factor is not spherical. Let $\Delta_k$ be the leftmost factor which is not spherical, ie., $\Delta_{\ell}$ is spherical, say a $(x_{\ell},x_{\ell})$-diagram, for every $\ell<k$. For $\ell \geq k$, say that $\Delta_{\ell}$ is a $(x_{\ell},y_{\ell})$-diagram. In $\Sigma^+$, we have the equality
\begin{center}
$x_1 \cdots x_{k-1} x_k x_{k+1} \cdots x_{n} = \mathrm{top}(\Delta)= \mathrm{bot}(\Delta)=x_1 \cdots x_{k-1} y_k y_{k+1} \cdots y_n$,
\end{center}
hence $x_k \cdots x_n = y_k \cdots y_n$. Because the equality holds in $\Sigma^+$, notice that $\mathrm{lg}(x_k)= \mathrm{lg}(y_k)$ implies $x_k=y_k$, which is impossible since we supposed that $\Delta_k$ is not spherical. Therefore, two cases may happen: either $\mathrm{lg}(x_k)> \mathrm{lg}(y_k)$ or $\mathrm{lg}(x_k)< \mathrm{lg}(y_k)$. Replacing $\Delta$ with $\Delta^{-1}$ if necessary, we may suppose without loss of generality that we are in the latter case. Thus, $x_k$ is a prefix of $y_k$: there exists some $y \in \Sigma^+$ such that $y_k=x_ky$. On the other hand, because $\Delta_k$ is a $(x_k,y_k)$-diagram, the equality $x_k=y_k$ holds modulo $\mathcal{P}$. A fortiori, $x_k=x_ky$ modulo $\mathcal{P}$. We deduce that
\begin{center}
$x_1 \cdots x_{k-1} x_k y^n x_{k+1} \cdots x_n \in [w]_{\mathcal{P}}$
\end{center}
for every $n \geq 1$. This implies that $[w]_{\mathcal{P}}$ is infinite, contradicting our hypothesis. $\square$

\medskip \noindent
\textbf{Proof of Theorem \ref{decompositioncor}.} We argue by induction on the algebraic dimension of the considered diagram group, ie., the maximal rank of a free abelian subgroup. If the algebraic dimension is zero, there is nothing to prove since the diagram group is trivial in this case. From now on, suppose that our diagram group $D(\mathcal{P},w)$ has algebraic dimension at least one. By applying Corollary \ref{cor:dichotomy}, we deduce that either $D(\mathcal{P},w)$ contains a contracting isometry of $X(\mathcal{P},w)$, and we conclude that $D(\mathcal{P},w)$ is either cyclic or acylindrically hyperbolic, or there exist two words $u,v \in \Sigma^+$ and a $(w,uv)$-diagram $\Gamma$ such that $D(\mathcal{P},w)= \Gamma \cdot ( D(\mathcal{P},u) \times D(\mathcal{P},v)) \cdot \Gamma^{-1}$ where $D(\mathcal{P},u)$ and $D(\mathcal{P},v)$ are non-trivial. In the first case, we are done; in the latter case, because $D(\mathcal{P},u)$ and $D(\mathcal{P},v)$ are non-trivial (and torsion-free, as any diagram group), their algebraic dimensions are strictly less than the one of $D(\mathcal{P},w)$, so that we can apply our induction hypothesis to find a $(u,u_1\cdots u_r)$-diagram $\Phi$ and a $(v,v_1 \cdots v_s)$-diagram $\Psi$ such that $$D(\mathcal{P},u)= \Phi \cdot D(\mathcal{P},u_1) \times \cdots \times D(\mathcal{P},u_r) \cdot \Phi^{-1}$$ and similarly $$D(\mathcal{P},v)= \Psi \cdot D(\mathcal{P},v_1) \times \cdots \times D(\mathcal{P},v_s) \cdot \Psi^{-1},$$ where, for every $1 \leq i \leq r$ and $1 \leq j \leq s$, $D(\mathcal{P},u_i)$ and $D(\mathcal{P},v_j)$ are either infinite cyclic or acylindrically hyperbolic. Now, if we set $\Xi= \Gamma \cdot (\Phi+ \Psi)$, then $$D(\mathcal{P},w)= \Xi \cdot D(\mathcal{P},u_1) \times \cdots \times D(\mathcal{P},u_r) \times D(\mathcal{P},v_1) \times \cdots \times D(\mathcal{P},v_s) \cdot \Xi^{-1}$$ is the decomposition we are looking for.

\medskip \noindent
Finally, the inequality $n \leq \dim_{\mathrm{alg}} D(\mathcal{P},w)$ is clear since diagram groups are torsion-free, and the inequality $\dim_{\mathrm{alg}} D(\mathcal{P},w) \leq \dim X(\mathcal{P},w)$ is a direct consequence of \cite[Corollary 5.2]{arXiv:1505.02053}. $\square$

\medskip \noindent
We conclude this section by showing that, in the cocompact case, being a contracting isometry can be characterized algebraically. For convenience, we introduce the following definition.

\begin{prop}\label{decompositionprop}
Let $\mathcal{P}= \langle \Sigma \mid \mathcal{R} \rangle$ be a semigroup presentation and $w \in \Sigma^+$ a baseword whose class modulo $\mathcal{P}$ is finite. If $g \in D(\mathcal{P},w) - \{ \epsilon(w) \}$, the following statements are equivalent: 
\begin{itemize}
	\item[(i)] $g$ is a contracting isometry of $X(\mathcal{P},w)$;
	\item[(ii)] the centraliser of $g$ is infinite cyclic;
	\item[(iii)] $g$ does not belong to a large canonical subgroup.
\end{itemize}
\end{prop}

\noindent
\textbf{Proof.} It is well-known that the centraliser of a contracting isometry of a group acting properly is virtually cyclic. Because diagram groups are torsion-free, we deduce the implication $(i) \Rightarrow (ii)$.

\medskip \noindent
Suppose that there exists a $(w,uv)$-diagram $\Gamma$ such that $g \in \Gamma \cdot D(\mathcal{P},u) \times D(\mathcal{P},v) \cdot \Gamma^{-1}$, where $D(\mathcal{P},u)$ and $D(\mathcal{P},v)$ are non trivial. Let $g_1 \in D(\mathcal{P},u)$ and $g_2 \in D(\mathcal{P},v)$ be two spherical diagrams such that $g= \Gamma \cdot (g_1+g_2) \cdot \Gamma^{-1}$. In particular, the centraliser of $g$ contains: $\Gamma \cdot ( \langle g_1 \rangle \times \langle g_2 \rangle) \cdot \Gamma^{-1}$ if $g_1$ and $g_2$ are non trivial; $\Gamma \cdot ( D(\mathcal{P},u) \times \langle g_2 \rangle) \cdot \Gamma^{-1}$ if $g_1$ trivial and $g_2$ is not; $\Gamma \cdot ( \langle g_1 \rangle \times D(\mathcal{P},v) ) \cdot \Gamma^{-1}$ if $g_2$ is trivial and $g_1$ is not. We deduce that $(ii) \Rightarrow (iii)$.   

\medskip \noindent
Finally, suppose that $g$ is not a contracting isometry. Up to taking a conjugate of $g$, we may suppose without loss of generality that $g$ is absolutely reduced. According to Theorem \ref{prop1} and Lemma \ref{lem1}, three cases may happen: 
\begin{itemize}
	\item[(a)] $\mathrm{supp}(g)$ has cardinality at least two;
	\item[(b)] $w=xyz$ in $\Sigma^+$ with $\mathrm{supp}(g)= \{y\}$ and $X(\mathcal{P},x)$ or $X(\mathcal{P},z)$ infinite;
	\item[(c)] $g^{\infty}$ contains a proper infinite prefix.
\end{itemize}
In case $(a)$, $g$ decomposes as a sum with at least two non-trivial factors, and we deduce from Lemma \ref{somme} that $g$ belongs to a large canonical subgroup. In case $(b)$, if $X(\mathcal{P},x)$ is infinite then $D(\mathcal{P},x)$ is non-trivial since $[x]_{\mathcal{P}}$ is finite, and we deduce that $g$ belongs to the large canonical subgroup $D(\mathcal{P},x) \times D(\mathcal{P},yz)$; we argue similarly if $X(\mathcal{P},z)$ is infinite. In case $(c)$, suppose that $g^{\infty}$ contains a proper infinite prefix $\Xi$. For convenience, write $g^{\infty}=g_1 \circ g_2 \circ \cdots$, where each $g_i$ is a copy of $g$. Suppose by contradiction that $g$ cannot be decomposed as a sum with at least two non-trivial factors. As a consequence, for each $i \geq 1$, the subdiagram $g_i$ either is contained into $\Xi$ or it contains a cell which does not belong to $\Xi$ but whose top path intersects $\Xi$. Because $\Xi$ is a proper prefix of $g^{\infty}$, there exists some index $j \geq 1$ such that $g_j$ is not included into $\Xi$, and it follows from Lemma \ref{cellprec} that, for every $i \geq j$, $g_i$ satisfies the same property. Let $\Xi_{j+r}$ denote the prefix of $g_1 \circ \cdots \circ g_{j+r}$ induced by $\Xi$. We know that, for every $0 \leq s \leq r$, the subdiagram $g_{j+s}$ contains a cell which does not belong to $\Xi$ but whose top path intersects $\Xi$. This implies that $\mathrm{bot}(\Xi_{j+r})$ has length at least $r$. On the other hand, the finiteness of $[w]_{\mathcal{P}}$ implies that the cardinality of $[\mathrm{bot}(\Xi_{j+r})]_{\mathcal{P}}$ is bounded by a constant which does not depend on $r$. Thus, we get a contradiction if $r$ is chosen sufficiently large. We have proved that $g$ decomposes as a sum with at least two non-trivial factors. As in case $(a)$, we deduce from Lemma \ref{somme} that $g$ belongs to a large canonical subgroup. $\square$

\begin{remark}
The implication $(iii) \Rightarrow (ii)$ also follows from the description of centralisers in diagram groups \cite[Theorem 15.35]{MR1396957}, since they are canonical subgroups.
\end{remark}

\appendix

\section{Combinatorial boundary vs. other boundaries}

\noindent
In this appendix, we compare the combinatorial boundary of a CAT(0) cube complex with three other boundaries. Namely, the \emph{simplicial boundary}, introduced in \cite{simplicialboundary}; the \emph{Roller boundary}, introduced in \cite{Roller}; and its variation studied in \cite{Guralnik}. In what follows, we will always assume that our CAT(0) cube complexes have countably many hyperplanes. For instance, this happens when they are locally finite.

\subsection{Simplicial boundary}

\noindent
For convenience, we will say that a CAT(0) cube complex is \emph{$\omega$-dimensional} if it does not contain an infinite collection of pairwise transverse hyperplanes. 

\medskip \noindent
Let $X$ be an $\omega$-dimensional CAT(0) cube complex. A \emph{Unidirection Boundary Set} (or \emph{UBS} for short) is an infinite collection of hyperplanes $\mathcal{U}$ not containing any facing triple which is:
\begin{itemize}
	\item \emph{inseparable}, ie., each hyperplane separating two elements of $\mathcal{U}$ belongs to $\mathcal{U}$;
	\item \emph{unidirectional}, ie., any hyperplane of $\mathcal{U}$ bounds a halfspace containing only finitely many elements of $\mathcal{U}$.
\end{itemize}
According to \cite[Theorem 3.10]{simplicialboundary}, any UBS is commensurable to a disjoint union of \emph{minimal} UBS, where a UBS $\mathcal{U}$ is \emph{minimal} if any UBS $\mathcal{V}$ satisfying $\mathcal{V} \subset \mathcal{U}$ must be commensurable to $\mathcal{U}$; the number of these minimal UBS is the dimension of the UBS we are considering. Then, to any commensurability class of a $k$-dimensional UBS is associated a \emph{$k$-simplex at infinity} whose faces correspond to its subUBS. All together, these simplices at infinity define a simplicial complex, called the \emph{simplicial boundary} $\partial_{\triangle}X$ of $X$. See \cite{simplicialboundary} for more details.

\medskip \noindent
The typical example of a UBS is the set of the hyperplans intersecting a given combinatorial ray. A simplex arising in this way is said \emph{visible}. Therefore, to any point in the combinatorial boundary naturally corresponds a simplex of the simplicial boundary. The following proposition is clear from the definitions:

\begin{prop}
Let $X$ be an $\omega$-dimensional CAT(0) cube complex. Then $(\partial^c X, \leq)$ is isomorphic to the face-poset of the visible part of $\partial_{\triangle} X$.
\end{prop}

\noindent
In particular, when $X$ is \emph{fully visible}, ie., when every simplex of the simplicial boundary is visible, the two boundaries essentially define the same objet. Notice however that Farley cube complexes are not necessarily fully visible. For example, if $\mathcal{P}= \langle a,b,p \mid a=ap, b=pb \rangle$, then $X(\mathcal{P},ab)$ is a ``diagonal half-plane'' and is not fully visible. Moreover, although Farley cube complexes are complete, they are not necessarily $\omega$-dimensional, so that the simplicial boundary may not be well-defined. This is the case for the cube complex associated to Thompson's group $F$, namely $X(\mathcal{P},x)$ where $\mathcal{P}= \langle x \mid x=x^2 \rangle$; or the cube complex associated to the lamplighter group $\mathbb{Z} \wr \mathbb{Z}$, namely $X(\mathcal{P},ab)$ where $\mathcal{P}= \langle a,b,p,q,r \mid a=ap,b=pb,p=q,q=r,r=p \rangle$.

\subsection{Roller boundary}

\noindent
A \emph{pocset} $(\Sigma, \leq, ^{\ast})$ is a partially-ordered set $(\Sigma,\leq)$ with an order-reversing involution $^{\ast}$ such that $a$ and $a^*$ are not $\leq$-comparable for every $a \in \Sigma$. A subset $\alpha \subset \mathfrak{P}(\Sigma)$ is an \emph{ultrafilter} if:
\begin{itemize}
	\item for every $a \in \Sigma$, exactly one element of $\{ a , a^* \}$ belongs to $\alpha$;
	\item if $a,b \in \Sigma$ satisfy $a \leq b$ and $a \in \alpha$, then $b \in \alpha$.
\end{itemize}
Naturally, to every CAT(0) cube complex $X$ is associated a pocset: the set of half-spaces $\mathcal{U}(X)$ with the inclusion and the complementary operation. In this way, every vertex can be thought of as an ultrafilter, called a \emph{principal ultrafilter}; see for example \cite{SageevCAT(0)} and references therein. For convenience, let $X^{\circ}$ denote the set of ultrafilters of the pocset associated to our cube complex. Notice that $X^{\circ}$ is naturally a subset of $2^{\mathcal{U}(X)}$, and can be endowed with the Tykhonov topology. 

\medskip \noindent
Now, the \emph{Roller boundary} $\mathfrak{R}X$ of $X$ is defined as the space of non principal ultrafilters of $X^{\circ}$ endowed with the Tykhonov topology. Below, we give an alternative description of $\mathfrak{R}X$. 

\medskip \noindent
Given a CAT(0) cube complex $X$, fix a base vertex $x \in X$, and let $\mathfrak{S}_x X$ denote the set of combinatorial rays starting from $x$ up to the following equivalence relation: two rays $r_1,r_2$ are equivalent if $\mathcal{H}(r_1)= \mathcal{H}(r_2)$. Now, if $\xi_1, \xi_2 \in \mathfrak{S}_x X$, we define the distance between $\xi_1$ and $\xi_2$ in $\mathfrak{S}_x X$ by
\begin{center}
$d(\xi_1, \xi_2)= 2^{- \min \{ d(x,J) \ \mid \ J \in \mathcal{H}(r_1) \oplus \mathcal{H}(r_2) \} }$,
\end{center}
where $\oplus$ denotes the symmetric difference. In fact, as a consequence of the inclusion $A \oplus B \subset (A \oplus B) \cup (B \oplus C)$ for every sets $A,B,C$, the distance $d$ turns out to be \emph{ultrametric}, i.e., 
$$d(\xi_1,\xi_3) \leq \max \left( d(\xi_1,\xi_2),d(\xi_2,\xi_3) \right)$$
for every $\xi_1,\xi_2, \xi_3 \in \mathfrak{S}_x X$. Also, notice that this distance is well-defined since its expression does not depend on the representatives we chose. 

\medskip \noindent
Given a combinatorial ray $r$, let $\alpha(r)$ denote the set of half-spaces $U$ such that $r$ lies eventually in $U$. 

\begin{prop}\label{RollerBoundary}
The map $r \mapsto \alpha(r)$ induces a homeomorphism $\mathfrak{S}_x X \to \mathfrak{R} X$.
\end{prop}

\noindent
\textbf{Proof.} If $r$ is a combinatorial ray starting from $x$, it is clear that $\alpha(r)$ is a non-principal ultrafilter, ie., $\alpha(r) \in \mathfrak{R} X$. Now, if for every hyperplane $J$ we denote by $J^+$ (resp. $J^-$) the halfspace delimited by $J$ which does not contain $x$ (resp. which contains $x$), then we have
\begin{center}
$\alpha(r)= \{ J^+ \mid J \in \mathcal{H}(r) \} \sqcup \{ J^- \mid J \notin \mathcal{H}(r) \}$.
\end{center}
We first deduce that $\alpha(r_1)= \alpha(r_2)$ for any pair of equivalent combinatorial rays $r_1,r_2$ starting from $x$, so that we get a well-defined induced map $\alpha : \mathfrak{S}_x X \to \mathfrak{R} X$; secondly, the injectivity of $\alpha$ follows, since $\alpha(r_1)=\alpha(r_2)$ clearly implies $\mathcal{H}(r_1) = \mathcal{H}(r_2)$.

\medskip \noindent
Now, let $\eta \in \mathfrak{R} X$ be a non-principal ultrafilter, and let $\mathcal{H}^+(\eta)$ denote the set of the hyperplanes $J$ satisfying $J^+ \in \eta$. Notice that, if $\mathcal{H}^+(\eta)$ is finite, then all but finitely many halfspaces of $\eta$ contain $x$, which implies that $\eta$ is principal; therefore, $\mathcal{H}^+(\eta)$ is infinite, say $\mathcal{H}^+(\eta)= \{ J_1,J_2, \ldots \}$. On the other hand, $J_i^+ \cap J_j^+$ is non-empty for every $i,j \geq 1$, so, for every $k \geq 1$, the intersection $C_k= \bigcap\limits_{i=1}^k J_i^+$ is non-empty; let $x_k$ denote the combinatorial projection of $x$ onto $C_k$. According to Lemma \ref{inclusion}, $x_{k+1}$ is the combinatorial projection of $x_k$ onto $C_{k+1}$ for every $k \geq 1$. Therefore, by applying Proposition \ref{projection}, we know that there exists a combinatorial ray $\rho$ starting from $x$ and passing through all the $x_k$'s. Now, if $J \in \mathcal{H}(\rho)$, there exists some $k \geq 1$ such that $J$ separates $x$ and $x_k$, and a fortiori $x$ and $C_k$ according to Lemma \ref{hyperplan séparant}. We deduce that either $J$ belongs to $\{J_1, \ldots, J_k \}$ or $J$ separates $x$ and some hyperplane of $\{ J_1, \ldots, J_k \}$; equivalently, there exists some $1 \leq \ell \leq k$ such that $J_{\ell}^+ \subset J^+$. Because $J_{\ell}^+ \in \eta$, we deduce that $J^+ \in \eta$. Thus, we have proved that $\mathcal{H}(\rho)= \mathcal{H}^+(\eta)$. It follows that $\alpha(\rho)= \eta$. The surjectivity of $\alpha$ is proved.

\medskip \noindent
Finally, it is an easy exercice of general topology to verify that $\alpha$ is indeed a homeomorphism. The details are left to the reader.
$\square$

\medskip \noindent
Therefore, the combinatorial boundary $\partial^c X$ may be thought of as a quotient of the Roller boundary. This is the point of view explained in the next section. 

\medskip \noindent
It is worth noticing that this description of the Roller boundary allows us to give a precise description of the Roller boundary of Farley cube complexes:

\begin{prop}
Let $\mathcal{P}= \langle \Sigma \mid \mathcal{R} \rangle$ be a semigroup presentation and $w \in \Sigma^+$ a base word. The Roller boundary of $X(\mathcal{P},w)$ is homeomorphic to the space of the infinite reduced $w$-diagrams endowed with the metric
\begin{center}
$(\Delta_1, \Delta_2) \mapsto 2^{- \min \{ \# \Delta \ \mid \ \Delta \ \text{is not a prefix of both} \ \Delta_1 ~ \text{and} \ \Delta_2 \} }$,
\end{center}
where $\# ( \cdot)$ denotes the number of cells of a diagram.
\end{prop}

\subsection{Guralnik boundary}

\noindent
Given a pocset $(\Sigma, \leq, ^*)$, Guralnik defines the boundary $\mathfrak{R}^*\Sigma$ as the set of almost-equality classes of ultrafilters, partially ordered by the relation: $\Sigma_1 \leq \Sigma_2$ if $\Sigma_2 \cap \overline{\Sigma_1}\neq \emptyset$, where $\overline{\cdot}$ denotes the closure in $\mathfrak{R}X$ with respect to the Tykhonov topology. See \cite{Guralnik} for more details.

\medskip \noindent
In particular, if $X$ is a CAT(0) cube complex, it makes sense to define the boundary $\mathfrak{R}^*X$ as the previous boundary of the associated pocset. Notice that $\mathfrak{R}^*X$ contains a particular point, corresponding to the almost-equality class of the principal ultrafilters; let $\mathfrak{p}$ denote this point. Then $\mathfrak{R}^*X \backslash \{ \mathfrak{p} \}$ is naturally the quotient of the Roller boundary $\mathfrak{R}X$ by the almost-equality relation. 

\begin{remark}
Although Guralnik considers only $\omega$-dimensional pocsets in \cite{Guralnik}, the definition makes sense for every pocset, and in particular is well-defined for every CAT(0) cube complex. Notice also that the boundary $\mathfrak{R}^*$ is called the \emph{Roller boundary} there, whereas the terminology now usually refers to the boundary defined in the previous section.
\end{remark}

\begin{prop}
The posets $(\mathfrak{R}^*X \backslash \{ \mathfrak{p} \}, \leq)$ and $(\partial^c X,\leq)$ are isomorphic.
\end{prop}

\noindent
\textbf{Proof.} If $r$ is an arbitrary combinatorial ray, choose a second ray $\rho \in \mathfrak{S}_xX$ equivalent to it, and set $\varphi(r) \in \mathfrak{R}^* X$ as the almost-equality class of $\alpha(\rho)$, where $\alpha : \mathfrak{S}_x X \to \mathfrak{R}(X)$ is the isomorphism given by Proposition \ref{RollerBoundary}. Notice that, if $\rho' \in \mathfrak{S}_xX$ is equivalent to $r$, then $\alpha(\rho)$ and $\alpha(\rho')$ are almost-equal, so we deduce that $\varphi(r)$ does not depend on our choice of $\rho$. As another consequence, $\varphi$ must be constant on the equivalence classes of rays, so that we get a well-defined map 
$$\varphi : \partial^c X \to \mathfrak{R}^* X \backslash \{ \mathfrak{p} \}.$$
Clearly, $\varphi$ is surjective. Now, if $r_1,r_2 \in \partial^c X$ satisfy $\varphi(r_1)=\varphi(r_2)$, then $\alpha(r_1) \underset{a}{=} \alpha(r_2)$; using the notation of the proof of Proposition \ref{RollerBoundary}, this implies that 
\begin{center}
$\mathcal{H}(r_1)= \mathcal{H}^+(\alpha(r_1)) \underset{a}{=} \mathcal{H}^+ (\alpha(r_2))= \mathcal{H}(r_2)$,
\end{center}
hence $r_1=r_2$ in $\partial^cX$. Thus, $\varphi$ is injective.

\medskip \noindent
Let $r_1,r_2 \in \partial^c X$ be two combinatorial rays satisfying $r_1 \prec r_2$, ie., $\mathcal{H}(r_1) \underset{a}{\subset} \mathcal{H}(r_2)$. According to Lemma \ref{basepoint} and Lemma \ref{basepointbis}, we may suppose without loss of generality that $r_1(0)=r_2(0)=x$ and $\mathcal{H}(r_1) \subset \mathcal{H}(r_2)$. Let $\mathcal{H}(r_2)= \{J_1, J_2, \ldots \}$. For every $k \geq 1$, let $C_k= \bigcap\limits_{i=1}^k J_i^+$ and let $x_k$ denote the combinatorial projection of $x$ onto $C_k$. Let $\mathcal{H}(r_1)= \{ H_1, H_2, \ldots \}$, and, for every $i \geq 1$, let $K_i= \bigcap\limits_{j=1}^i H_j \cap C_k$; notice that, because $\mathcal{H}(r_1) \subset \mathcal{H}(r_2)$, $K_i$ is non-empty. Fix some $k \geq 1$. Finally, for every $i \geq 1$, let $y_i$ denote the combinatorial projection of $x_k$ onto $K_i$. By combining Lemma \ref{inclusion} and Proposition \ref{projection}, we know that there exists a combinatorial ray starting from $x_k$ and passing through all the $y_i$'s; let $\rho_k$ denote the concatenation of a combinatorial geodesic between $x$ and $x_k$ with the previous ray. 

\medskip \noindent
If $\rho_k$ is not a combinatorial ray, then there exists a hyperplane $J$ separating $x$ and $x_k$, and $x_k$ and some $y_j$. On the other hand, it follows from Lemma \ref{hyperplan séparant} that $J$ must separate $x$ and $C_k$, so that $J$ cannot separate $x_k$ and $y_j$ since $x_k,y_j \in C_k$. Therefore, $\rho_k$ is a combinatorial ray. 

\medskip \noindent
If $\mathcal{H}_k$ denotes the set of the hyperplanes separating $x$ and $x_k$, then 
\begin{center}
$\mathcal{H}(\rho_k) = \mathcal{H}_k \cup \mathcal{H}(r_1)$.
\end{center}
Indeed, if $J$ is a hyperplane intersecting $\rho_k$ which does not separate $x$ and $x_k$, then $J$ must separate $x_k$ and some $y_j$; a fortiori, $J$ must be disjoint from $K_j$ according to Lemma \ref{hyperplan séparant}. Thus, if $z$ is any vertex of $r_1 \cap K_j$, we know that $J$ separates $x_k$ and $y_j$, but cannot separate $x$ and $x_k$, or $y_j$ and $z$; therefore, necessarily $J$ separates $x$ and $z$, hence $J \in \mathcal{H}(r_1)$. Conversely, let $J$ be a hyperplane of $\mathcal{H}(r_1)$. In particular, $J \in \mathcal{H}(r_2)$, so $\rho_k$ eventually lies in $J^+$. Because $\rho_k(0)=x \notin J^+$, we conclude that $J \in \mathcal{H}(\rho_k)$. 

\medskip \noindent
In particular, because $\mathcal{H}_k$ is finite, we deduce that $\mathcal{H}(\rho_k) \underset{a}{=} \mathcal{H}(r_1)$, ie., $r_1= \rho_k$ in $\partial^cX$. On the other hand, it is clear that, for any finite subset $\mathcal{H} \subset \mathcal{H}(r_2)$, there exists some $k \geq 1$ sufficiently large so that $\mathcal{H} \subset \mathcal{H}_k$. Therefore, the sequence of ultrafilters $(\alpha(\rho_k))$ converges to $\alpha(r_2)$ in the Roller boundary. 

\medskip \noindent
Consequently, $\alpha(r_2)$ belongs to the closure of the almost-equality class of $\alpha(r_1)$. This precisely means that $\varphi(r_1) \leq \varphi(r_2)$ in $\mathfrak{R}^* X$. We have proved that $\varphi$ is a poset-morphism.

\medskip \noindent
Now, let $r_1,r_2 \in \partial^c X$ such that $\varphi(r_1) \leq \varphi(r_2)$ in $\mathfrak{R}^*X$. This means that there exists a combinatorial ray $\rho$ such that $\alpha(\rho)$ belongs to the intersection between the almost-equality class of $\alpha(r_2)$ and the closure of the almost-equality class of $\alpha(r_1)$ in the Roller boundary. So there exists a sequence of combinatorial rays $(\rho_k)$ such that $\rho_k=r_1$ in $\partial^cX$ for every $k \geq 1$ and $\alpha(\rho_k) \to \alpha(\rho)$. Let $r_0$ be the minimal element of the class of $r_1$ in $\mathfrak{S}_x X$ which is given by Lemma \ref{minimal} below. Then $\mathcal{H}(r_0) \subset \mathcal{H}(\rho_k)$ for every $k \geq 1$, hence $\mathcal{H}(r_0) \subset \mathcal{H}(\rho)$. Therefore,
\begin{center}
$\mathcal{H}(r_1) \underset{a}{\subset} \mathcal{H}(r_0) \subset \mathcal{H}(\rho) \underset{a}{=} \mathcal{H}(r_2)$,
\end{center}
hence $r_1 \prec r_2$ in $\partial^c X$. We have proved that $\varphi^{-1}$ is a poset-morphism as well. $\square$

\begin{lemma}\label{minimal}
Let $r \in \mathfrak{S}_xX$ be a combinatorial ray. There exists $r_0 \in \mathfrak{S}_x X$ equivalent to $r$ such that, for any combinatorial ray $\rho \in \mathfrak{S}_xX$ equivalent to $r$, we have $\mathcal{H}(r_0) \subset \mathcal{H}(\rho)$.
\end{lemma}

\noindent
\textbf{Proof.} Let $\mathcal{H} \subset \mathcal{H}(r)$ be the set of the hyperplanes $J \in \mathcal{H}(r)$ such that $r$ does not stay in a neighborhood of $J$; say $\mathcal{H}=\{ J_1, J_2, \ldots \}$. For every $k \geq 1$, let $C_k= \bigcap\limits_{i=1}^k J_i^+$ and let $x_k$ denote the combinatorial projection of $x$ onto $C_k$. By combining Lemma \ref{inclusion} and Proposition \ref{projection}, we know that there exists a combinatorial ray $r_0$ starting from $x$ and passing through all the $x_k$'s. We claim that $r_0$ is the ray we are looking for.

\medskip \noindent
Let $\rho \in \mathfrak{S}_xX$ be a combinatorial ray equivalent to $r$, and let $J \in \mathcal{H}(r_0)$ be a hyperplane. By construction, $J$ separates $x$ and some $x_k$; a fortiori, it follows from Lemma \ref{hyperplan séparant} that $J$ separates $x$ and $C_k$. On the other hand, given some $1 \leq \ell \leq k$, because $r$ does not stay in neighborhood of $J_{\ell}$, there exist infinitely many hyperplanes of $\mathcal{H}(r)$ which are included into $J_{\ell}^+$; therefore, since $\rho$ is equivalent to $r$, we deduce that $\rho$ eventually lies in $J_{\ell}^+$. A fortiori, $\rho$ eventually lies in $C_k$. Consequently, because $J$ separates $x=\rho(0)$ and $C_k$, we conclude that $J \in \mathcal{H}(\rho)$. $\square$

\addcontentsline{toc}{section}{References}

\bibliographystyle{alpha}
\bibliography{CICCC}

\end{document}